\documentclass{article} 
\usepackage{url}
\usepackage{wrapfig} 
\usepackage{floatrow}
\title{Lazily Adapted Constant Kinky Inference for Nonparametric Regression and Model-Reference Adaptive Control}

\author{
Jan-Peter Calliess\\ 
Dept. of Engineering Science\\
University of Oxford, UK\footnote{At the time of the submission of the original version, the author was a member of the University of Cambridge.}
}
\usepackage[margin=1.2in]{geometry}
\usepackage{subfigure}
\usepackage{amsmath}
\usepackage{amssymb}
\usepackage{amsthm}
\usepackage{tikz}
\newtheorem{thm}{Theorem}[section]
\newtheorem{cor}[thm]{Corollary}
\newtheorem{lem}[thm]{Lemma}

\theoremstyle{definition}
\newtheorem{defn}[thm]{Definition}
\theoremstyle{remark}
\newtheorem{rem}[thm]{Remark}
\newtheorem{remark}[thm]{Remark}

\newtheorem{ex}[thm]{\textbf{Example}}

\newcommand{\matnorm}[1]{{\left\vert\kern-0.25ex\left\vert\kern-0.25ex\left\vert #1 
    \right\vert\kern-0.25ex\right\vert\kern-0.25ex\right\vert}}
\newcommand{\opnorm}[1]{{\left\vert\kern-0.25ex\left\vert\kern-0.25ex\left\vert #1 
    \right\vert\kern-0.25ex\right\vert\kern-0.25ex\right\vert}}		
\newcommand{\specnorm}[1]{\matnorm{#1}_2}
\newcommand{\norm}[1]{\left\Vert#1\right\Vert}

\newcommand{\abs}[1]{\left\vert#1\right\vert}

\newcommand{\Real}{\mathbb R}

\newcommand{\nat}{\mathbb N}

\newcommand{\argmin}{\text{argmin}}
\newcommand{\argmax}{\text{argmax}}

\newcommand{\vc}[1]{#1}

\newcommand{\SP}[2]{\ensuremath{\mathbf{\langle} \vc{#1} \mathbf{,} \vc{#2} \mathbf{\rangle}}}

\newcommand{\ball}[2]{\ensuremath{\mathfrak B_{\vc{#1}} \bigl( #2\bigr) }}

\newcommand{\dif}[2]{{\operatorname{d}\over\operatorname{d}#2}#1}

\renewcommand{\d}[1]{\text{ d}#1}

\newcommand{\floor}[1]{ \left\lfloor #1 \right\rfloor }
\newcommand{\ceil}[1]{ \left\lceil #1 \right\rceil }

\newcommand{\specrad}{\rho}

\newcommand{\data}{\ensuremath{ \mathcal D} }

\newcommand{\param}{\ensuremath{\Xi}}

\newcommand{\state}{\ensuremath{ \vc x}} 

\newcommand{\statespace}{\ensuremath{ \mathcal X}}

\newcommand{\inspace}{\ensuremath{ \mathcal X}}
\newcommand{\outspace}{\ensuremath{ \mathcal Y}}

\newcommand{\fctspace}{\ensuremath{ \mathcal F}}

\newcommand{\iaspace}{\ensuremath{ \mathcal S}} 

\newcommand{\grid}{\ensuremath{  G}}

\newcommand{\indset}{\ensuremath{ \mathcal I}}
\newcommand{\indsett}{\ensuremath{ {\mathcal I_{t}}}}

\newcommand{\maxerrn}{\bar{\mathfrak N}} 

\newcommand{\metric}{\, \mathfrak{d}} 
\newcommand{\metricp}{\,  \tilde {d}} 

\newcommand{\Metrici}[2]{\metric\bigl(#1,#2\bigr) }
\newcommand{\Metrico}[2]{\metric_\outspace\bigl(#1,#2\bigr) }

\newcommand{\dist}{\ensuremath{ \text{dist}}} 

\newcommand{\predf}{\, \mathfrak{  \hat f}} 
\newcommand{\predfn}{\, \mathfrak{  \hat f_n}} 
\newcommand{\predfnj}{\, \mathfrak{  \hat f_{n,j}}} 
\newcommand{\prederr}{\, \mathfrak{\hat v}} 
\newcommand{\prederrn}{\, \mathfrak{\hat v_n}} 
\newcommand{\prederrnj}{\, \mathfrak{\hat v_{n,j}}} 

\newcommand{\noise}{{\nu}} 

\newcommand{\hexp}{{ \alpha }}
\newcommand{\hestthresh}{\ensuremath{ \lambda}}
\newcommand{\ubf}{\, {\bar B}} 
\newcommand{\lbf}{\, {\underline{B}}} 

\newcommand{\hoelset}[3]{\mathfrak{H}_{#2}(#1,#3)}
\newcommand{\Hoelset}{\mathfrak{H}}

\newcommand{\obserrpar}{e} 
\newcommand{\obserr}{\mathfrak e} 
\newcommand{\obserrbnd}{\bar{\mathfrak e}}

\newcommand{\decke}{\ensuremath{\mathfrak u}}
\newcommand{\boden}{\ensuremath{\mathfrak l}}

\newcommand{\seq}[2]{\ensuremath{\bigl(#1\bigr)_{#2}}}

\renewcommand{\Pr}{\mathrm{Pr}}

\renewcommand{\d}{\ensuremath{\text{ d}}}

\newcommand{\beq}{\begin{equation}}
\newcommand{\eeq}{\end{equation}}

\newcommand{\eqn}[2]{\begin{equation} #2 \label{#1}\end{equation}}
\newcommand{\dt}{\ensuremath{ \, dt \,}}
\newcommand{\tinc}{\ensuremath{ \Delta}}

\newcommand{\convton}{\stackrel{n \to \infty}\longrightarrow}
\newcommand{\convto}{\longrightarrow}
\newcommand{\setconvto}{\longrightarrow}

\newcommand{\bd}{\leadsto} 
\newcommand{\bdu}{\twoheadrightarrow} 

\newcommand{\errmetric}{\mathcal E} 
\newcommand{\queryset}{I} 

\begin{document}

\maketitle
\begin{abstract}
Techniques known as \emph{Nonlinear Set Membership} prediction, \emph{Lipschitz Interpolation} or \emph{Kinky Inference} are  approaches to machine learning that utilise \emph{presupposed} Lipschitz properties to compute inferences over unobserved function values. Provided a bound on the true best Lipschitz constant of the target function is known a priori they offer convergence guarantees as well as bounds around the predictions.
Considering a more general setting that builds on H\"older continuity relative to pseudo-metrics, we propose an online method for estimating the H\"older constant online from function value observations that possibly are corrupted by bounded observational errors. Utilising this to compute  adaptive parameters within a kinky inference rule gives rise to a nonparametric machine learning method, for which we establish strong \emph{universal approximation} guarantees. That is, we show that our prediction rule can learn any continuous function in the limit of increasingly dense data to within a worst-case error bound that depends on the level of observational uncertainty.
We apply our method in the context of \emph{nonparametric model-reference adaptive control (MRAC)}. 
Across a range of simulated aircraft roll-dynamics and performance metrics our approach outperforms recently proposed alternatives that were based on \emph{Gaussian processes} and \emph{RBF-neural networks}. For discrete-time systems, we provide guarantees on the tracking success of our learning-based controllers both for the batch and the online learning setting. 
\end{abstract}

\section{Introduction}

Typically, a controller is designed on the basis of a dynamical model of the system. When little is known about these dynamics a priori or, if the dynamics may be subject to unexpected change, machine learning methods can be employed to learn such a model online on the basis of measurements. 

Supervised machine learning methods are algorithms for inductive inference. On the basis of a sample, they construct (learn) a computable model of a data generating process that facilitates inference over the underlying ground truth function and aims to predict its function values at unobserved inputs.
Among supervised learning methods, nonparametric algorithms tend to offer greater flexibility to learn rich 
function classes (e.g. rich classes of nonlinear dynamics). 
 
Perhaps the most popular nonparametric machine learning method is Bayesian inference with \textit{Gaussian processes (GPs)} \cite{GPbook:2006}. GPs offer a flexible and principled probabilistic method for nonparametric regression and have evolved into one of the primary work-horses for learning dynamic systems \cite{Deisenroth2009,Deisenroth2011,Deisenroth2015,Tuongmodellearningsurvey2011,deisenrothsurvey2013,McHutchonthesis2014} in the research communities related to artificial intelligence. 
However, they suffer from several limitations, including scalability to large data sets, a lack of understanding of how to bound the closed-loop dynamics resulting from controlling on the basis of a GP state-space model and the question of how to choose a good prior in a principled yet computable manner. To alleviate the last problem, it is common practice to tune hyper-parameters of a chosen (typically universal) kernel to explain the data via the marginal log-likelihood. While often successful on many data sets, the result can be highly sensitive to the choice of optimiser, initialisations, data sets and computational budget. Unfortunately, little theoretical understanding of the important interplay between these components in the resulting inference mechanism seems to exist.

In contrast to such Bayesian methods this work considers an extension of ideas existing in applied mathematics  (e.g. \cite{Shubert:72,Zabinsky2003,Cooper1995,Cooper2006,Baran2008,Beliakov2006}) as well as in control \cite{Milanese2004,Canale2014} that harnesses Lipschitz regularity of the target function to provide bounds on the predictions of the target function at unobserved inputs. Applied to machine learning, the basic idea is that Lipschitz continuity constrains the set of possible function values of a target function at a query input as a function of the distance between the query and the previously observed training examples. A prediction is then made by choosing a function value in the middle of the set of possible function values. This idea, which at least goes back to \cite{Sukharev1978}, has been leveraged in different fields under different headlines including \emph{Lipschitz Interpolation} \cite{Zabinsky2003,Beliakov2006} and \emph{Nonlinear Set Membership (NSM)}  methods \cite{Milanese2004}. If the Lipschitz constant is known, Lipschitz interpolation provides uncertainties around the predictions of function values at unobserved inputs. The uncertainties are maximally tight if no other knowledge is known other than the presupposed Lipschitz regularity \cite{Sukharev1978,calliess2014_thesis}.
The presupposed Lipschitz constant is a crucial parameter of the inference rule, quite similar to the choice of a prior (e.g. via kernel and hyperparameters)  in Bayesian inference.

We would argue that one of the advantages of these methods is their computational simplicity. That is, they are numerically robust and only involve basic computational steps that could be more efficiently computed even on a simple embedded RISC micro-controller.
A practical concern is that the predictions and bounds of these methods hinge on the a priori knowledge of the presupposed Lipschitz parameter of the underlying target function. Some previous works remark that, in the absence of such knowledge the constants might in practice be estimated from the data (e.g. via estimators discussed in  \cite{Strongin1973,Wood1996,Beliakovsmoothing2007,Milanese2004,Beliakov2006}). In fact, \cite{Milanese2004} suggest fitting a parametric regressor to the data first and utilise the Lipschitz constant of the fitted model for the NSM approach. Unfortunately, there is no theoretical analysis given anywhere what the impact of that estimate is to the predictive performance of the regression method. (And, consequently control methods that rely on the resulting predictions would be unable to assert closed-loop convergence or robustness guarantees).

Our work addresses this deficiency by proposing an approach that allows us to provide learning and tracking guarantees even in the absence of a priori knowledge of this constant.

 As a first step, we rehearse \emph{Kinky Inference (KI)} \cite{calliess2014_thesis} which generalises the Lipschitz interpolation and NSM frameworks in several ways.
%
We then combine the KI machine learning approach with a simple online parameter estimator that allows us to incrementally compute an estimate of the H\"older constant on the basis of incrementally arriving (noisy) data. This merger yields a new inference rule we refer to as \emph{Lazily Adapted Constant Kinky Inference (LACKI)}. We prove that this rule is sample-consistent (up to the level observational error in the data) and that the LACKI predictors themselves are H\"older continuous. This allows us to establish strong universal approximation properties: That is, in the limit of infinite data, the LACKI rule is capable of approximating not only any H\"older continuous target function but also any non-H\"older continuous function with arbitrarily low error (up to a bound dependent on the level of observational error in the training data). Since our LACKI rule can be seen as an extension of Lipschitz interpolation with empirical Lipschitz constant estimation as considered in \cite{Beliakov2006}, our results also provide new theoretical guarantees for this previously proposed interpolation rule where the Lipschitz parameters are estimated online from the noisy training data with our proposed modified constant estimation method.

In addition to these learning-theoretic considerations, we apply our LACKI approach to online-learning based model-reference adaptive control.

 As a testbed, we replicate simulations of  the roll dynamics of an F-4 fighter aircraft under uncertain wing rock previously considered by other authors in model-reference adaptive control \cite{Chowdhary2013,Monahemi1996,chowdharyacc2013} and compare our controller against their methods. Here our LACKI-based controller outperforms its competitors across a range of metrics including computational speed, prediction and tracking accuracy. 

For discrete-time feedback-linearisable systems with uncertain nonlinear drift, we provide theoretical guarantees on the tracking success of our LACKI- model-reference adaptive controller both in the batch and in the online learning setting. 

In contrast to much of the standard literature of probabilistic nonparametric regression (e.g. \cite{Gyoerfi2002,Tsybakov2009}), our analysis focusses on the derivation of deterministic  worst-case error bounds. While possibly being more conservative, we would argue that this type of analysis has the benefit of being more meaningful in a control setting where the learner receives training examples and queries that will typically violate statistical assumptions typically presupposed  in the statistical literature.

The remainder of the paper is structured as follows: 

Sec. \ref{sec:LACKI_all} contains the core of the LACKI regression methods. Following a rehearsal of the kinky inference (KI) framework for nonparametric learning, Sec. \ref{sec:lacki} describes the our LACKI approach for setting the KI parameters. Sec. \ref{sec:properties_lacki} is dedicated the derivation of several properties of the resulting LACKI approach, including our consistency guarantees.

Sec. \ref{sec:MRAC_lacki} contains the control part of the paper. We first introduce the framework of model-reference adaptive control in which we propose a controller based on our LACKI learning method.  For illustration purposes, we closely follow the setting of wing-rock control considered in \cite{Chowdhary2013,ChowdharyCDC2013} and compare our LACKI-based controller to other MRAC controllers consdered and proposed by previous work. The section concludes by giving convergence guarantees for LACKI-MRAC in discrete-time systems.

The paper concludes with Sec. \ref{sec:concl}, summarising our findings and containing various suggestions for future work. The appendix contains a variety of background material on 
H\"older continuity and various supplementary derivations referred to at various points of the main body of the paper. 

In comparison to the original 2016 version of this preprint, this updated version adds some experiments in Sec. \ref{sec:regr_benchmarks} and corrects a typo that had existed in Thm. \ref{thm:trackerrconv}.
This paper has originally been released in 2016. Since then alternative versions have been presented at the  European control conference \cite{calliess2018ECC} and been published in Automatica \cite{calliess2020Automatica}. The present update of this preprint contains a correction of the exposition of the proof of Lemma 2.16.

\section{Kinky Inference with lazily adapted constants}

\label{sec:LACKI_all}

\subsection{Kinky Inference}
\label{sec:KI_core}
In this section, we will introduce the class of learning rules we refer to as \emph{Kinky Inference}. They encompass a host of other methods such as Lipschitz Interpolation and Nonlinear Set Interpolation.  
%

The rules possess a parameter $L(n)$ that needs to be specified by any KI algorithm. In this paper, we are most concerned with studying LACKI, a KI rule algorithm where $L(n)$ coincides with a noise-robust and multi-variate generalisation of Strongin's estimate \cite{Strongin1973} of a H\"older constant computed from the data set $\data_n$ available at time step $n$. 

\textbf{Setting.}  Let $\inspace$, $\outspace$ be two spaces endowed with (pseudo-) metrics $\metric: \inspace^2 \to \Real_{\geq 0}, \metric_\outspace:\outspace^2 \to \Real_{\geq 0}$, respectively. Spaces $\inspace, \outspace$ will be referred to as \textit{input space} and \textit{output space}, respectively. 
It will be convenient to restrict our attention to input and output spaces that are additive abelian groups and which are \emph{translation-invariant} with respect to their (pseudo-) metrics. That is, for the input space $\inspace$, we assume $ \metric(x+x',x''+x') = \metric(x,x''),\forall x,x',x'' \in \inspace$. 

For simplicity, throughout the remainder of this work, we will assume the output space is the canonical Hilbert space $\outspace = \Real^m$ $(m \in \nat)$ endowed with the $\metric_{\outspace}(y,y') = \norm{y-y'}_\infty, \forall y,y' \in \outspace$. 

Let $f: \inspace \to \outspace$ be a \emph{target} or \emph{ground-truth} function we desire to learn. For our purposes, learning means regression. That is, we utilise the data to construct a computable function that allows us predict values of the target function at any given input.


Assume that, at time step $n$, we have access to a \textit{sample} or \textit{data set} $\data_n:= \{\bigl( s_i, \tilde f_i \bigr) \, \vert \, i=1,\ldots, N_n \} $ containing $N_n \in \nat$ sample vectors $\tilde f_i = \tilde f(s_i) \in \outspace$ of an \emph{observable}   function $\tilde f$ at sample input $s_i \in \inspace$. Here, the observable $\tilde f : \inspace \to \outspace$ is a ``noise-corrupted'' version of the true \emph{target function} $f: \inspace \to \outspace$ we would like to make inferences about on the basis of the available sample. In this work, we will typically assume that the observable $\tilde f$ coincides with the target $f$ up to a level of interval-bounded observational noise: $\forall x \in \inspace: \metric_\outspace(\tilde f(x),f(x) ) \leq \obserr(x) $ where $\obserr : \inspace \to \Real_{\geq 0}$ is the error bound function whose values we assume to be bounded by some (known or unknown) bound $\obserrbnd \in \Real_{\geq 0}$. We can model the situation by the presence of a bounded additive observational error (or ``noise'') function $\noise:\inspace \to \outspace$ with $\tilde f = f + \noise$.
The interpretation of these errors depends on the given application and this ``noise'' may be deterministic or stochastic. For instance, in the context of system identification, the sample might be based on noisy measurements of velocities and it may be due to sensor noise or, the noise might model systematic error such as numerical approximation errors. 

Furthermore, $\obserr$ may also accommodate \textit{input} uncertainty (that is when predicting $f(x)$, $x$ is uncertain) (for details refer to \cite{calliess2014_thesis}). In the course of our theoretical considerations below the error will also serve to absorb the discrepancy between a H\"older and a non-H\"older function.
%

\textbf{Learning.}
It is our aim to learn target function $f$ in the sense that, combining prior knowledge about $f$ with the observed data $\data_n$, we infer \textit{predictions} $\predfn(x)$ of $f(x)$ at unobserved \textit{query inputs} $x \notin \grid_n$. Here, $\grid_n =\{s_i | i =1\ldots,N_n\} \subset \inspace$ refers to the (not necessarily regular) \emph{grid} of sample inputs. The entire function $\predfn$ that is learned to facilitate predictions is referred to as the \textit{predictor}. Since the computation of the predictor is based on the available data and utilised to make inferences over unobserved inputs, we can view the learning process as an instance of (inductive) inference. Therefore, the formula to compute the predictor $\predfn$ will also be referred to as an inference rule.
%
%
In our context, we will understand a \textit{machine learning algorithm} to implement a such an inference rule. That is, it is a computable function that maps a data set $\data_n $ to a 
 predictor $\predf_n$ (and possibly an uncertainty estimate function $\prederr_n$). 
A typical \textit{desideratum} of a good predictor is that it is \textit{efficiently computable}. Its learning performance is measured in terms of the degree and rapidity it  \textit{converges to the target} (up to the observational error given by $\obserr$) in the limit of increasingly dense data. Of course there are many different metrics with respect to which one can assess convergence. Perhaps, the most convenient one is mean-square convergence. However, inspired by our control applications, we desire to investigate worst-case convergence rates which hold independently from distributional assumptions and will yield performance guarantees even in zero-measure events.
\textbf{The Kinky Inference Learning rules.}
In this work we will expand on the basis of the following class of predictors to perform learning as inference over unobserved function values:

\begin{defn}[Kinky Inference (KI) rule ] \label{def:KIL}
Let  $\Real_\infty := \Real \cup\{- \infty, \infty\}$ and $\inspace$ be some space endowed with a pseudo-metric $\metric$. Let $\lbf,\ubf: \inspace \to \outspace \subseteq \Real_\infty^m$ denote \textit{lower- and upper bound functions}, that can be specified in advance and assume $\lbf(x) \leq \ubf(x), \forall x \in I \subset \inspace$ component-wise. 
Furthermore, let $\obserrpar$ denote a parameter that specifies a deterministic belief about the true observational error bound $\obserr$.
Given sample set $\data_n$, we define the predictive functions $\predfn: \inspace \to \outspace, \prederrn: \inspace \to \Real^m_{\geq 0}$ to perform inference over function values.
For $j=1,\ldots,m$, their $j$th output components are given by:
	\begin{align*}
   \predfnj(x; \param(n)) :=& \frac{1}{2} \min\{ \ubf_j(x), \decke_{n,j}(x;\param(n)\bigr)\} \\
    &+ \frac{1}{2} \max\{ \lbf_j(x), \boden_{n,j}(x;\param(n)\bigr) \},\\
	\prederrnj(x; \param(n)) :=& \frac{1}{2} \min\{ \ubf_j(x), \decke_{n,j}\bigl(x;\param(n)\bigr)\} \\
	&- \frac{1}{2} \max\{ \lbf_j, \boden_{n,j}\bigl(x;\param(n)\bigr) \}.
	\end{align*}
	Here, $\decke_n\bigl(\cdot;\param(n)\bigr), \boden_n\bigl(\cdot;\param(n)\bigr): \inspace \to \Real^m$ are called ceiling and floor functions, respectively. Their $j$th component functions are given by
	\[\decke_{n,j}\bigl(x; \param(n)\bigr) := \min_{i=1,\ldots,N_n}   \tilde f_{i,j} +  \tilde \metric(x,s_i;\param(n)) \] and 
	\[\boden_{n,j}\bigl(x; \param(n)\bigr) := \max_{i=1,\ldots,N_n}   \tilde f_{i,j} -  \tilde  \metric(x,s_i;\param(n)),\] respectively.
  Here, $\tilde \metric (\cdot,\cdot;\param(n))$ is a mapping parameterised by $ \param(n)$. While there are many conceivable parameterisations, we restrict our attention to the case where, for some pseudo-metric $\metric$ on the input space $\inspace$, we have $\param(n) = (L(n), \hexp,\obserrpar)$ with $$\tilde \metric(\cdot,\cdot; \param(n)) = L(n) \metric^\hexp(\cdot,\cdot) + \obserrpar(x).$$
  As will be seen below,  parameter $L(n)$ has the interpretation of a H\"older constant of the predictor relative to pseudo-metric $\metric$ while $\hexp \in (0,1]$ can be interpreted as a H\"older exponent (cf. Thm. \ref{lem:LACKIpredHoelder}). That is, we will show that $\predfn(\cdot;\param(n))$ belongs to the class $$\Hoelset \bigl(L(n),\hexp\bigr) = \{ \phi: \inspace \to \outspace | \forall x,x' \in \inspace : \metric_\outspace(\phi(x),\phi(x') ) \leq L(n) \metric(x,x')^\hexp\}$$ of $L(n)-\hexp$- H\"older continuous functions. Note, we could alternatively re-express this H\"older class as a generalised class of Lipschitz functions $\text{Lip}(L(n))= \{ \phi: \forall x,x' \in \inspace: \Metrico{\phi(x)}{\phi(x')} \leq L(n) \tilde \metric (x,x')\} $ where, for any $\hexp \in (0,1]$, $\tilde \metric = \metric^\alpha$ is a pseudo-metric provided  $\metric$ is (refer to Lem. \ref{lem:hoeldererror_metric}). However, as it is often customary to define Lipschitz and H\"older continuity in a more restricted sense relative to standard norm-induced metrics (in which case the H\"older class is strictly more general than the Lipschitz class) we chose to refer to H\"older continuity rather than Lipschitz continuity to highlight that such H\"older functions can be learned as well.
  
  As insinuated by our notation, we consider parameter $L(n)$ to be adaptive, i.e.  data-dependent while the other parameters are assumed to be set in advance.
 Function $\obserrpar$ can be utilised to accommodate observational noise. 
That is, if the noise level in the data is assumed to be contained in $[-\obserrbnd, \obserrbnd]$ then one would choose $\obserrpar(x) = \obserrbnd,\forall x$.
In addition,  functions $\lbf,\ubf:  \inspace \to \Real_\infty^m$ are parameters that have to be specified in advance and can impose a priori knowledge of bounds on the target function. For example, if we know the target function to exclusively map to nonnegative values, then one can set $\lbf(x) = 0,\forall x$.
To disable restrictions of boundedness, it is allowed to specify the upper and lower bound functions to constants $\infty$ or $-\infty$, respectively.

Function $\predfn$ is the predictor that is to be utilised for predicting/inferring function values at unseen inputs. Function $\prederrn(x;\param(n))$ is meant to quantify the uncertainty of prediction $\predfn(x;\param(n))$. 
For ease of notation, we shall often omit explicit mention of the parameter, e.g. we may write $\predfn(x)$ instead of $\predfn(x;\param(n))$.
\end{defn}

To provide an intuition of the inference rule, consider the following special case where we have access to a noise-free sample $\data_n$ and suppose 
the target $f$ is a real-valued $L^*-\hexp$ H\"older continuous function. Observing the noise-free sample point $(s_i,f_i)$ constrains the set of function values $f(x)$ to the set $\mathbb S_i(x) =\{ \phi \in \outspace | \metric_\outspace(\phi, f_i) \leq L^* \metric(s_i,x) \}$. Considering a set of sample points $\data_n$, target value $f(x)$ is constrained to lie in the intersection $\mathbb S(x)=\cap_{i=1}^{N_n} \mathbb S_i(x)$. It is easy to see that the floor and ceiling functions are tight lower and upper bounds of $\mathbb S(x)$ with $\mathbb S(x) := \{ \phi \in \outspace | \boden_n(x; L^*) \leq \phi \leq \decke_n(x;L^*)\}$. In other words, setting parameter $L(n)$ to the best H\"older constant $L^*$ and bounds $\lbf =-\infty,\ubf=+\infty$ yields a predictor $\predfn(x)$ that for every query $x$ chooses the mid-point of the set $\mathbb S(x)$ of those function values that can possibly be assumed by a H\"older continuous function that interpolates the observed sample. Prediction error $\prederrn(x)$ simply is the radius of the set.

For the case of $\hexp=1$, this approach is known as \emph{Lipschitz interpolation} \cite{Beliakov2006,Zabinsky2003}. Since a set is utilised for interpolation, the approach is also known as \emph{Nonlinear Set Interpolation} \cite{Milanese2004,Canale2014} in control. Specification of $\ubf,\lbf$ allows us to incorporate additional knowledge and constrain our set $\mathbb S(x)$ further. 
For instance, when estimating densities we might incorporate the knowledge of dealing with nonnegative functions. In this case, it makes sense to set $\lbf$ to a constant value of zero yielding $\mathbb S(x) = \{\phi | \phi \geq 0 \} \cap_{i=1}^{N_n} \mathbb S_i(x) $.

When choosing $L(n)$ to coincide with the best H\"older constant, one can give strong guarantees of convergence to the target as on the tightness of the prediction bounds  \cite{Sukharev1978,calliess2014_thesis} showing that bounds are as tight as possible without imposing additional assumptions and that the predictor minimises the worst-case risk.

Unfortunately, this requires us to know at least an upper bound of $L^*$ and therefore, several authors have proposed different approaches of how to estimate the constant from the data (e.g. \cite{Strongin1973,Wood1996}). However, it appears to be largely unknown how to do so in the presence of bounded observational noise $\obserr >0$ in a principled manner. Furthermore, when we replace $L(n)$ by the empirical estimates, nothing seems to be known about the convergence properties of the resulting kinky inference rule that is based on such  estimates.  

In the remainder of the paper, we shall address this gap. Firstly, we propose an estimator to be utilised in place of $L(n)$ that can be set to be robust to noise (i.e. does not grow unbounded). Referring to the resulting KI rule as LACKI, we then prove universal approximation properties of the LACKI rule before considering its performance in a control application.

\subsection{Lazily Adapted Constant Kinky Inference (LACKI)}
\label{sec:lacki}
Above we explained the benefits of choosing parameter $L(n)$ to coincide with a H\"older constant of the target. 
However, if such a constant is unavailable a priori, we desire to compute $L(n)$ as a data-dependent estimate of the H\"older constant. Our proposal for such an estimator will be introduced next. 

For notational convenience, for two sets $S,S' \subset \inspace$ of inputs we define  $$U(S,S') := \{(s,s') \in S \times S' | \metric(s,s') >0\} \text{ and  let } U_n := U(\grid_n,\grid_n) $$ be the set of all pairs grid inputs deemed disparate under the pseudo-metric $\metric$.

The \emph{best} H\"older constant of a function $f$ is the smallest nonnegative number $L^*$ such that $f$ is contained in the set 
$\hoelset {L^*}{}{p} = \{\phi: \inspace \to \outspace \, \vert \, \forall x,x' \in \inspace : \metric_\outspace \bigl(\phi(x),\phi(x')\bigr) \leq L^* \, \bigl( \metric (x,x') \bigr)^\hexp\} $ of $L^*-\hexp$- H\"older continuous functions. So, this best H\"older constant is given by  
$$L^* = \sup_{(x,x') \in U(\inspace,\inspace)} \frac{\metric_\outspace \bigl(f(x) - f(x')\bigr)}{\metric^\hexp(x,x') }.$$

Given the noisy data $\data_n = \{(s_i,\tilde f_i) | i=1,\ldots,N_n\}$ a natural estimate of the best H\"older constant might be to compute $\hat \ell_n^* := \max_{(s,s') \in U_n} \frac{\metric_\outspace(\tilde f(s),\tilde f(s')) }{\metric^\hexp(s,s')}$ \cite{Strongin1973}. In the absence of \emph{noise} (may it be stochastic or deterministic), that is, if $\tilde f_i = f(s_i),\forall i$, $\hat \ell_n^*$ never overestimates the true best H\"older constant. That is, $\hat \ell_n^* \leq L^*$. However, in the presence of noise $\noise: \inspace \to \outspace$ (such that $\tilde f = f+ \noise$) this boundedness assumption of the estimates no longer holds true. In particular if the noise is not H\"older continuous, we expect $\hat \ell_n^*$ to grow unbounded with increasingly dense data.
For practical reasons and for the sake of our theoretical arguments presented below, we desire the parameters $L(n)$ to remain bounded. Thus, without further modifications $\hat \ell_n^*$ is not a suitable candidate for $L(n)$.

To ensure bounded estimates even in the presence of noise, we propose the following estimator : %
\begin{equation}\label{eq:lazyconstupdaterule_batch_main}
\ell(\data_n;\hestthresh,\underline L)  := 
 \max \Bigl\{ \underline L, \max_{(s,s') \in U_n} \frac{\metric_\outspace(\tilde f(s),\tilde f(s')) - \hestthresh}{\metric^\hexp(s,s')} \Bigr\}.
\end{equation}
\begin{figure}
        \centering
    \includegraphics[width = 8.8cm,height = 4cm]
								{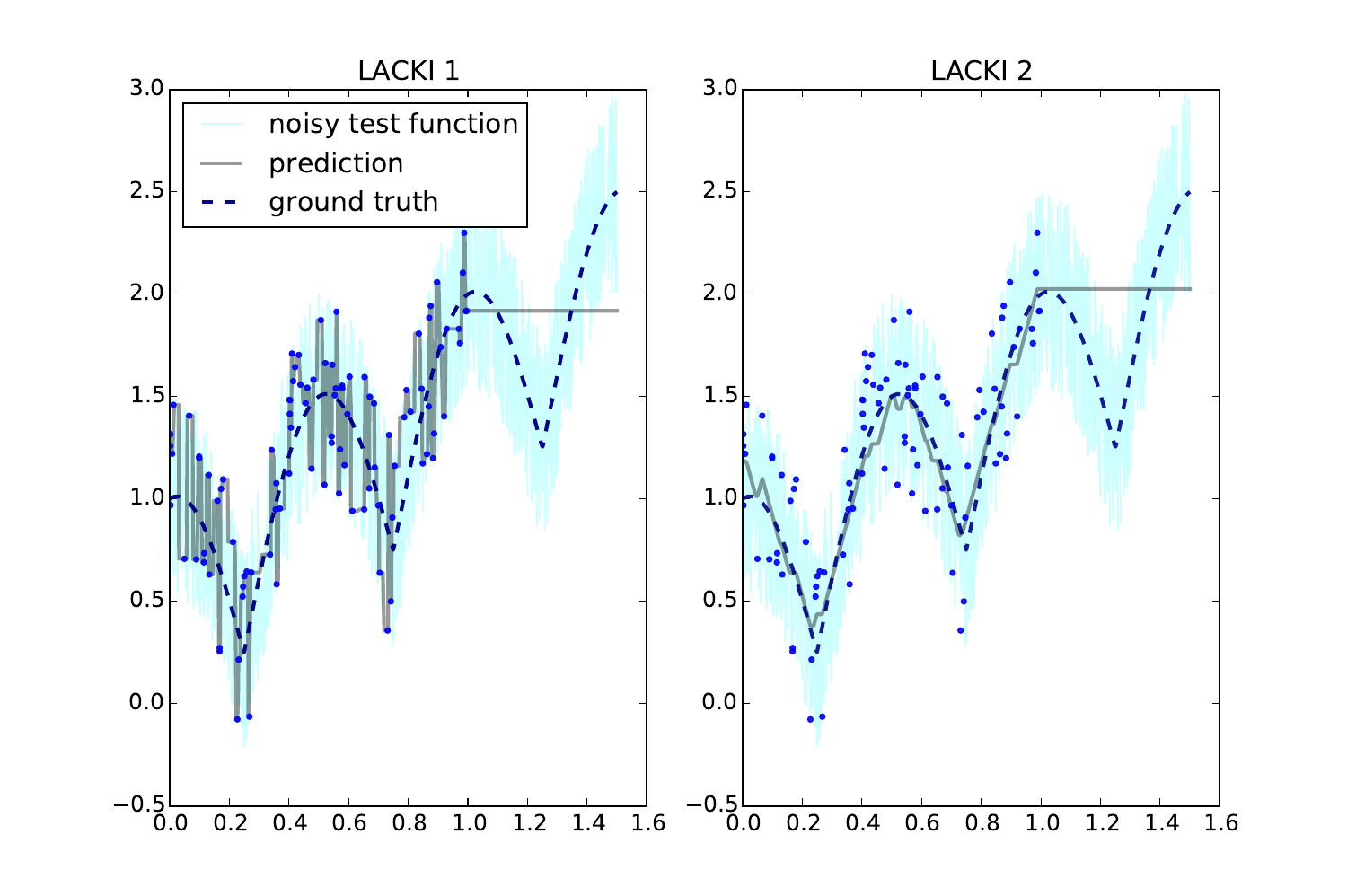}
   \caption{Two LACKI inferences over function values of the target $f: x \mapsto \abs{\cos(2\pi x)}+x$ (dashed line) on the basis of a noisy sample (plotted as dots). The predictors are plotted in grey, the noisy observational function $\tilde f (\cdot) = f(x) + \nu_x$ is plotted in cyan. Here the $\nu_x$ were drawn i.i.d. at random from a uniform distribution on the interval $[-0.5,0.5]$. In both cases we chose the parameters $\ubf \equiv \infty,\lbf \equiv -\infty,\underline L = 0$ and $\hexp=1$.
   The left plot shows the LACKI predictor $\predfn$ (grey line) for parameter choice $\hestthresh =0$, falsely assuming absence of observational noise. As a result, the prediction overfitted to the noise. The right plot depicts the prediction $\predfn$ (grey curve) for correct parameter choice $\hestthresh = 2 \obserrbnd =1$, causing the noise to be smoothed out and resulting in more accurate prediction of the underlying ground truth $f$.  }
			\label{fig:LACKInoise}
\end{figure}

Here $\underline L$ is a parameter that can be used to specify a priori knowledge of a lower bound on the best Lipschitz constant. In the absence of particular domain-specific knowledge, one can of course always set $\underline L =0$.
\begin{rem} \label{rem:bndedlipconstestimates}
By setting parameter $\hestthresh$ at least as large as twice the maximum level of observational noise, i.e. $\hestthresh = 2 \obserrbnd + q$ for some $q\geq0$, it is easy to see that the $\ell(\data_n;\hestthresh,0)$ are bounded from above by $\bar L  =\sup_{x,x', \metric(x,x') >0} \frac{\metric_\outspace(f(x), f(x'))   -q}{\metric^\hexp(x,x')} \leq L^* $ ( and,  $L^* < \infty$ if the target is H\"older continuous). 
\end{rem}
Next, consider an online learning situation where the available data increases over time. 
That is, $\data_n \subseteq \data_{n+1}$ for all time steps $n \in \nat$. 
For time step $n \in \nat$, let $S_{n+1} := G_{n+1} \backslash \grid_n$ be the set of new sample inputs.
We can define an incremental update rule recursively as follows: 
\begin{align} \label{eq:Hoelconstlazyupdateincr_main}
\ell_{n+1} := \max\Bigl\{ & \ell_n, \max_{(s,s') \in U(\grid_n, S_{n+1})} \frac{\metric_\outspace\bigl(\tilde f(s),\tilde f(s')\bigr) - \hestthresh}{\metric^\hexp(s,s')},\\
&\max_{(s,s') \in U(S_{n+1}, S_{n+1})} \frac{\metric_\outspace\bigl(\tilde f(s),\tilde f(s')\bigr) - \hestthresh}{\metric^\hexp(s,s')} \Bigr\},
\end{align} for $n \in \nat$ 
and where 
$\ell_0 := \underline L$. 
The effort of computing $\ell_{n+1}$ is in the order of $\mathcal O\bigl(M (\abs{S_{n+1}} N_n+ \abs{S_{n+1}}^2)\bigr)$ where $M$ denotes the effort for evaluating the pseudo-metrics.
By construction, we have $\ell_n =\ell(\mathcal \data_n;\hestthresh,\underline L), \forall n \in \nat.$ 

\begin{rem}\label{rem:convlipconstestimates}
Remember that $\seq{\ell_n}{n\in \nat}$ is bounded. Since it is also growing monotonically we can appeal to the monotone convergence theorem to show that the sequence is convergent to some number $\bar L \leq \max\{L^*,\underline L\}$.
\end{rem}

So far, we have defined a rule of how to update noise-robust and convergent estimates $\ell_n$ of the H\"older constant. Using these data-dependent estimates in place of $L(n)$ in our kinky inference framework as per Def. \ref{def:KIL} yields an inference rule that shall henceforth be referred to as \textit{Lazily Adapted Constant Kinky Inference (LACKI)}.

\begin{defn}[LACKI rule] \label{def:LACKI} For each output component $j \in \{1,\ldots,\text{dim} \, \outspace \}$ 
define $\predfn(\cdot)_{j}$  as per Def. \ref{def:KIL} but assume we choose the parameters $L(n) := \ell(\data_n;\hestthresh,\underline L)$ (according to Eq. \ref{eq:lazyconstupdaterule_batch_main}). We refer to the resulting predictor $\predfn\bigl(\cdot \bigr)$ as a \emph{Lazily Adapted Constant Kinky Inference (LACKI)} rule. Here, the free parameters are $\hexp \in (0,1], \hestthresh \in \Real_{\geq 0}$ and $\underline L \in \Real_{\geq 0}$.  
\end{defn}

To develop a first feel for our inference rule, refer to Fig. \ref{fig:LACKInoise}. Here, we depicted the predictors for an underlying ground-truth function on the basis of a sample but with different parameter choices $\hestthresh$. When setting this parameter to two times the observational noise level, the predictor accurately smoothes out the noise. In contrast, when the parameter is set to zero, the resulting predictor will perfectly interpolate through the noisy observations, thereby limiting the approximation quality to the level of observational noise. 

Furthermore, we notice that the predictors are H\"older continuous but non-differentiable. Informally speaking, the inference exhibits ``kinks'', motivating the term ``kinky inference''.

Finally note, the estimator $\ell_n$ determining  parameter $L(n)$ is ``lazy'' in the sense that it only increases the estimate of the H\"older constant just enough to be consistent with the observed data. That is, it chooses $L(n)$ to coincide with the smallest H\"older constant of a conceivable target function $f$ that could have generated the data under the given noise assumption. Below, we will see that the predictor $\predfn$ has H\"older constant $L(n)$. Therefore, the ``laziness'' of the estimator of $L(n)$ implements \emph{Occam's razor}: it \emph{regularises} the hypothesis space of continuous functions to prefer simple explanations of the data (i.e. functions with low H\"older constants) over complex ones (i.e. functions with higher H\"older constants). Here, $\hestthresh$ can serve as a parameter that can be utilised to regularise the predictor further in order to compensate for (bounded) noise in the data.

\begin{figure}
        \centering
    \includegraphics[width = 8.8cm,height = 4cm]
								{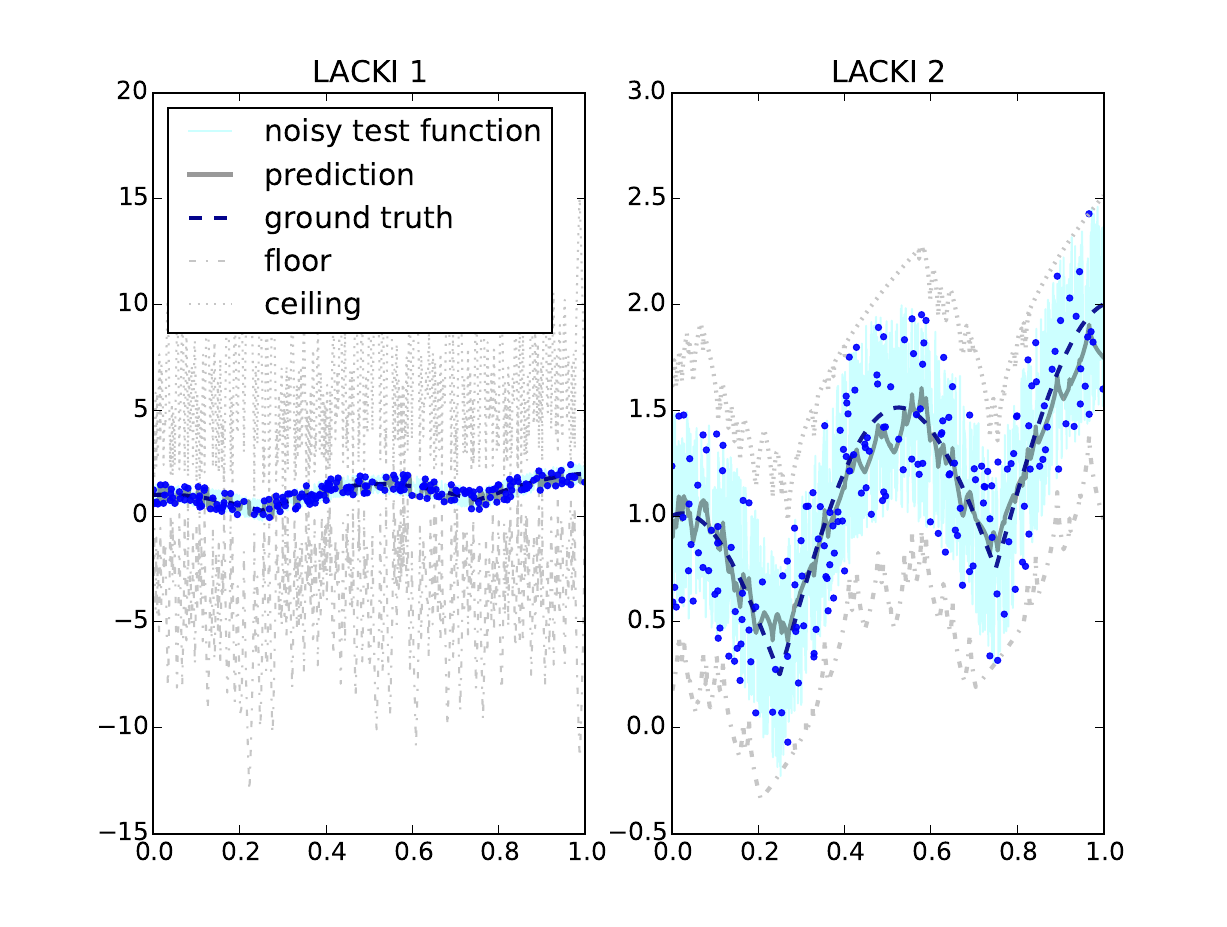}
   \caption{Repetition of the experiment but with $p=0.5$. This time, we also plotted floor and ceiling functions (grey dotted and dashed-dotted curves) delimiting the uncertainty bounds. Note, when setting $\lambda =0$ the estimate $\ell(\mathcal \data_n;\hestthresh,\underline L)$ was found to be 128 resulting in extremely conservative uncertainty estimates (left figure). In contrast, choosing $\lambda =2 \obserrbnd $ gave a parameter estimate $\ell(\mathcal \data_n;\hestthresh,\underline L) = 2.4$ yielding sensible uncertainty bounds (right figure). }
			\label{fig:LACKInoise2}
\end{figure}

\subsection{Properties}
\label{sec:properties_lacki}
We will now establish several properties of the LACKI rules including boundedness of the predictors, sample-consistency and H\"older continuity. Most importantly however, we will show that the LACKI rules are universal approximators, in the sense that they can be set to learn any continuous function with arbitrary worst-case error.

The core idea behind this can be sketched as follows: 
First, we establish H\"older continuity and sample-consistency. This allows us to prove that LACKI can learn any H\"older function.
Note, some universal approximators, such as RBFNs with Gaussian kernels, are provably Lipschitz. Therefore, learning any continuous function can be interpreted as learning some Gaussian RBFN with an observational error level that absorbs the discrepancy between the RBFN and the ground truth. Since a finite RBFN with smooth, bounded-derivative kernel is provably H\"older and since we can learn any H\"older function with LACKI up to the level of observational error, we can learn the continuous ground-truth up to the approximation error of the RBFN. 

Following this outline, we will now proceed to establish the desired properties formally.

\begin{lem}[Boundedness of the predictor]
Irrespective of the boundedness of input space $\inspace$ and assuming finite sample size  $N_n=\abs{\data_n} < \infty$, the predictor $\predfn:\inspace \to \outspace$ is bounded. In particular,
for $\hexp =1$, we have  
$\forall x \in \inspace: \norm{\predfn(x)}_\infty \leq \max_{i=1,...,N_n} \norm{\tilde f_i}_\infty + \frac{L(n)}{2}  \max_{i,j = 1,...,N_n} \metric^\alpha(s_i,s_j) <\infty$.
\begin{proof}
Let $D = \max_{i,j = 1,...,N_n} \metric^\hexp(s_i,s_j)$ and for the $k$th output dimension let  $F_k= \max_{i=1,...,N_n} \abs{\tilde f_{i,k}}$. As shown in Sec. \ref{sec:Hoelder_brief}, $\metric^\hexp$ is a pseudo-metric too and hence, adheres to the triangle inequality. Utilising the definition of the predictor and the triangle inequality we see that, for any  $x \in \inspace$ and any output dimension $k$, there are some $i,j \in \{1,...,N_n\}$ such that we have: 
$\predf_{n,k}(x) =\frac{\tilde f_{j,k} + \tilde f_{i,k}}{2} + \frac{L(n)}{2} \bigl( \metric^\hexp(x,s_i) - \metric^\hexp(x,s_j) \bigr) \leq \frac{\tilde f_{j,k} + \tilde f_{i,k}}{2} +  \frac{L(n)}{2} \metric^\hexp(s_j,s_i) \leq F_k + \frac{L(n)}{2}  D < \infty$.
\end{proof}
\label{lem:LACKIpredbounded}
\end{lem}

As promised, we establish that the predictors of the LACKI inference rule are H\"older continuous:
\begin{lem}[H\"older regularity of LACKI]
With definitions as before, let $(\outspace,\metric_\outspace) = (\Real^m,(x,y) \mapsto \norm{x-y}_\infty)$.
Provided that the bounding functions $\lbf,\ubf$ are H\"older continuous (or set to $-\infty,\infty$, respectively),
the predictors $\predfn$ are H\"older continuous $(n \in \nat)$ with constant $L(n)$ and exponent $\hexp$. That is, $\forall n \in \nat: \predfn \in \Hoelset \bigl(L(n),\hexp\bigr)$.

\begin{proof}
It is easy to show that the one-dimensional mappings of the form $x \mapsto \ell \Metrici{x}{x'}^\hexp$ are $\ell-\hexp-$ H\"older continuous for any choices of $\ell,\hexp$ and inputs $x'$. Furthermore, taking point-wise $\max$, $\min$ as well as averages of H\"older continuous functions is known to not change their H\"older properties (e.g. cf. \cite{calliess2014_thesis}). Therefore, the output-component predictors  $\predf_{n,j}$ $(j=1,...,m)$ are $L(n)$-$\hexp$- H\"older. 
Hence, $\forall x,x': \Metrico{\predfn(x)}{\predfn(x')} = \norm{\predfn(x) - \predfn(x')}_\infty = \max_{j=1,...,m} \abs{\predf_{n,j}(x) - \predf_{n,j}(x')} \leq \max_{j=1,...,m} L(n) \metric(x,x')^\hexp = L(n) \metric(x,x')^\hexp$.
\end{proof} 
\label{lem:LACKIpredHoelder}
\end{lem} 

We now establish how well our LACKI rule can interpolate the training data as function of the noise bound and regularisation parameter $\hestthresh$: 
\begin{lem}[Sample-consistency of LACKI] \label{lem:LACKIsampleconsistency} If for each output dimension $j \in \{1,...,d\}$ and some $\hestthresh \geq 0$ we have $L(n)  \geq  \max_{(s,s') \in U_n} \frac{\abs{\tilde f_j(s)- \tilde f_j(s')} -\hestthresh }{\metric^\hexp(s,s')}
$ then the LACKI rule is sample-consistent (up to $\frac \hestthresh 2$). 
That is,  \[\forall q \in \{1,\ldots,N_n \}:  \predfn(s_q) \in \ball{\frac \hestthresh 2}{\tilde f_q} \] where $\ball{\frac \hestthresh 2}{\tilde f_q} =\{ x \in \outspace |  \norm{x-\tilde f_q}_\infty \leq \frac \hestthresh 2\}$ denotes the $\frac \hestthresh 2$-ball around the observation $\tilde f_q$.\\
Thus, we also have $\norm{ f(s_q) - \predfn(s_q)}_\infty \leq \frac \hestthresh 2 + \norm{\obserr(s_q)}_\infty \leq \frac \hestthresh 2 + \obserrbnd$.
\begin{proof}
Remember, our output-space metric is given by $\Metrico{y}{y'} = \norm{y-y'}_\infty$.
For ease of notation, we will confine our proof to the case of one-dimensional outputs ($d=1$). The multi-dimensional case follows trivially from the one-dimensional result by applying it to each output component function. 
Let $n \in \nat  $ be fixed and, for ease of notation, write $L := L(n)$. Let $j,k \in   \{1,\ldots,N_n \} $ such that $j \in \argmin_i \tilde f_i + L \metric^\hexp(s_i,s_q) $ and 
$k \in \argmax_i \tilde f_i - L \metric^\hexp(s_i,s_q) $. 
By definition of $\predfn$ we have:
\begin{align}
\predfn(s_q) 
&= \frac 1 2 \bigl(\underbrace{ \tilde f_j + L \metric^\hexp(s_j,s_q)}_{:=B}   \bigr) + \frac 1 2 \bigl( \underbrace{\tilde  f_k - L \metric^\hexp(s_k,s_q) }_{=:A}\bigr). \label{eq:hru582jsokbbbn}
\end{align}
(i) Firstly, we show  \underline{$ A \in [\tilde f_q, \tilde f_q +\hestthresh]$}:
If $k = q$, this holds trivially true since then $A= \tilde f_q$. 
So, assume $k \neq q$. 
We have $\tilde f_k \geq \tilde f_k - L \metric^\hexp(s_k,s_q)  \geq \tilde f_q - L \metric^\hexp(s_q,s_q)  = \tilde f_q$ where the second inequality holds due to $k \in \argmax_i \tilde f_i - L \metric^\hexp(s_i,s_q) $. That is,
\begin{equation} 
A=\tilde f_k - L \metric^\hexp(s_k,s_q) \geq \tilde f_q. \label{eq:9p4t37ru8ewihk3hjtu}
\end{equation} 
On the other hand, since $L  \geq  \max_{(s,s') \in U_n} \frac{\abs{\tilde f(s)-\tilde f(s')} - \hestthresh}{\metric^\hexp(s,s')}$  we have in particular: $L  \geq  \frac{ \abs{\tilde f_k- \tilde f_q}-\hestthresh }{\metric^\hexp(s_k,s_q)}$. Thus,  
 $L \metric^\hexp(s_k,s_q) + \hestthresh \geq \abs{\tilde f_k-\tilde f_q} = \tilde f_k - \tilde f_q$. Hence, $\tilde f_q + \hestthresh \geq \tilde f_k - L \metric^\hexp(s_k,s_q) = A$. 
 Together with (\ref{eq:9p4t37ru8ewihk3hjtu}) we have shown $A \in [\tilde f_q , \tilde f_q +\hestthresh].$
%

(ii) The proof of \underline{$B \in [ \tilde f_q - \hestthresh, \tilde f_q]$} is completely analogous to that of (i) and hence, is omitted.

(iii) Together, the statements in (i) and (ii)  entail  $\predfn(s_q) = \frac 1 2 A + \frac 1 2 B  \in [\tilde f_q - \frac \hestthresh 2, \tilde f_q + \frac \hestthresh 2]$.

Hence, $\Metrico{\predfn(s_q)}{\tilde f(s_q)} \leq \frac \hestthresh 2$.

Moreover, for any sample input $s_q$ we have $\predfn(s_q) = f(s_q) + \phi_q + \psi_q$ with $\metric_{\outspace}(0,\psi_q) \leq \frac \hestthresh 2,  \metric_{\outspace}(0,\phi_q) \leq \Metrico 0 {\obserr(s_q)} \leq \obserrbnd$. 
Our output-space metric is translation-invariant and hence, $\Metrico{f(s_q)}{\predfn(s_q)} = \Metrico{0}{\predfn(s_q)- f(s_q)} = \Metrico{0}{  \phi_q + \psi_q} \leq \frac \hestthresh 2 +\Metrico 0 {\obserr(s_q)} \leq \frac \hestthresh 2 + \obserrbnd$.
\end{proof}
\end{lem}
%

\subsubsection{Prediction error bounds and consistency}

To asses our learning rule, we might be interested the discrepancy $\metric_\fctspace(\predfn,f)$ between the predictor $\predfn$ and the target function $f$ relative to some metric $\metric_{\fctspace}$ between functions in the space $\fctspace$ of continuous functions. In statistics, a typical choice is the mean-square error metric assessed with respect to some distribution over inputs, the function space and the noise. However, in many safety-critical applications, often arising in control, worst-case error considerations are of greater value, leading to a worst-case metric 
$\metric_\fctspace(f,g) = \sup_{x \in I} \Metrico{ f(x)}{g(x)} $ for some subset $I \subseteq \inspace$ of queries on finds interesting to take into consideration. 

Therefore, we will now establish worst-case consistency guarantees of our LACKI inference rules. That is, we shall study the worst-case error  
sequence $\errmetric^\infty :=\seq{\errmetric^\infty_n}{n \in \nat}$, 
\begin{equation}
\errmetric^\infty_n := \sup_{x \in \queryset } \metric_\outspace\bigl(\predfn(x), f(x)\bigr)\end{equation}
for data $\data_n$ that becomes increasingly dense over time relative to a set of query inputs $\queryset \subseteq \inspace$. 
To clarify the latter concept, 
consider the sequence of grids $\seq{\grid_n}{n \in \nat}$. 
We say this sequence converges to a set that \emph{becomes dense relative to a set $I$ in the limit of large $n$ } if we can use points in the sequence to approximate any points in $I$ with increasing accuracy. That is, if $\forall \epsilon >0,x \in I \exists n_0 \forall n \geq n_0 \exists g \in \grid_n: \metric(x,g) < \epsilon$. If the rate at which this happens is independent of $x$ then we say that the grid sequence becomes dense \emph{uniformly}. This is the case iff $\forall \epsilon >0 \exists n_0 \forall n \geq n_0, x \in I \exists g \in \grid_n: \metric(x,g) < \epsilon$.

To make the rates explicit in our notation, we list the following general definitions: 


\begin{defn} [Becoming dense, rates, $\stackrel{r}{\convto}, \stackrel{r}{\bd}, \stackrel{r}{\bdu}$]
Let $\inspace$ be a space endowed with a pseudo-metric $\metric$. 
Let $r:\nat \to \Real$ be a ``rate'' function. 
that vanishes, that is, with $\lim_{n \to \infty} r(n) = 0$ (i.e. $r \in o(1)$).
\begin{itemize}
\item The sequence $s =\seq{s_n}{n \in \nat}$ of points in $\inspace$ is said to converge to a point $x \in \inspace$ \emph{with rate} $r$ (denoted by $s \stackrel{r} \convto x$) iff  $\forall n \in \nat : \metric(x,s_n) \leq r(n)$ and $r(n) \convton 0$. 
\begin{itemize}
\item The sequence $s $ is said to converge to a set $\mathbb S \subset \inspace$  with rate $r:\nat \to \Real$ (denoted by $s \stackrel{r} \setconvto \mathbb S$) iff $\forall n\in \nat : \inf_{x \in \mathbb S}\metric(x,s_n) \leq r(n)$ and $r(n) \convton 0$. 
\end{itemize}
\item A sequence of sets $S = \seq{S_n}{n \in \nat}$ is said to \emph{become dense} relative to $x \in \inspace$  \emph{with rate} $r$ (denoted by $S \stackrel{r}{\bd} x$) iff $S$ contains a point sequence that converges to $x$ with that rate. That is, iff $\exists s= \seq{s_n}{n\in \nat}: s \stackrel{r}\convto x \wedge \forall n: s_n \in S_n$. 
\begin{itemize}
\item Similarly, the sequence of sets $S$ is said to become dense relative to a \emph{set}  of points $\mathbb S \subset \inspace $ (denoted by $S \bd \mathbb S$) iff it becomes dense relative to all points of $\mathbb S$, i.e. iff $\forall x \in \mathbb S : S \stackrel{r_x}{\bd} x$ for some vanishing rate $r_x : \nat \to \Real$. 
\item The sequence is becoming dense relative to $\mathbb S$ \emph{uniformly} (denoted by $S \bdu \mathbb S$) iff there is a single vanishing rate for all $x \in \mathbb S$. That is, if $\exists r:\nat \to \Real:$ $\lim_{n\to \infty} r(n) =0 \wedge \sup_{x \in \mathbb S} \inf{s_n \in S_n} \metric(s_n,x) \leq r(n), \forall n$. Function $r$ is referred to as the convergence rate and we write $S \stackrel{r}{\bdu}\mathbb S$ to denote that $S$ becomes dense relative to $\mathbb S$ with uniform rate $r$.
\end{itemize}
\end{itemize}
\end{defn}

\begin{thm}[LACKI can learn any H\"older function] \label{thm:convergenceifboundedconstandsamplecons_LACKI}
Assume the following holds true:
\begin{enumerate}
\item The observational errors given by $\obserr$ are bounded from above by $\obserrbnd= \sup_x \metric_\outspace\bigl(0,\obserr(x)\bigr) \in \Real_{\geq 0}$.
\item  The target $f: \inspace \to \outspace$ is H\"older continuous, i.e. $\exists L^* \in \Real: f \in \hoelset {L^*}{}  p$.
\end{enumerate}
%

Under these assumptions we can give the following guarantees: 

\textbf{(A)} If the grid becomes dense (pointwise), the point-wise worst-case error vanishes up to $\frac \hestthresh 2 + \obserr$:
$$
\text{If }\forall x \in I \subset \inspace \exists r_x \in o(1) : L(\cdot) r_x^\hexp(\cdot) \in o(1) \wedge \seq{G_n}{n \in \nat} \stackrel{r_x}{\bd} x  \text{ then }\forall x \in I : \seq{\Metrico{\predfn(x)}{f(x)}}{n \in \nat} \stackrel{\varrho_x} \convto [0,\frac \hestthresh 2 + \obserr]$$
where for the error convergence rate $\varrho_x$ we have $\varrho_x(n) \leq (L(n)+ L^*) r_x^\hexp(n) , \forall n \in \nat$. 

\textbf{(B)} If the grid becomes dense in $I \subset \inspace$ uniformly, then the worst-case prediction error vanishes uniformly (up to $\frac \lambda 2 + \obserrbnd$):
$$\text{If } \exists r \in o(1) :  L(\cdot) r^\hexp(\cdot) \in o(1) \wedge (G_n) \stackrel{r}{\bdu} I  \text{ then }  \errmetric^\infty \stackrel{\varrho}{\convto} [0,\frac \hestthresh 2 + \obserrbnd]$$
where for the uniform error convergence rate $\varrho$ we have $\varrho(n) \leq (L(n)+ L^*) r^\hexp(n) , \forall n \in \nat$. 

%
%
%


\begin{proof}
We have established that the predictors $\predfn(\cdot)$ of the LACKI rule are $L(n)$-$\hexp$- H\"older (Lem. \ref{lem:LACKIpredHoelder}) and sample-consistent up to level $\frac \hestthresh 2$ (Lem. \ref{lem:LACKIsampleconsistency}). 

For any input $x \in \inspace$ let $\xi_n^x$ denote a nearest neighbour of $x$ in grid $\grid_n$. That is, $\xi_n^x \in \arg\inf_{s \in \grid_n} \metric(x,s)$.
Since $\grid_n$ is assumed to become dense in the input domain $\inspace$, for any input $x$ there is a rate function $r_x : \nat \to \Real_{\geq 0}$ such that $r_x(n) \stackrel{n \to \infty}{\longrightarrow} 0$ and $ \metric(x,\xi_n^x )^\hexp \leq r_x(n) 
, \forall n \in \nat$. In the case of uniform convergence a rate function can be chosen independently of $x$ and will be denoted by $r$ rather than $r_x$.

\textbf{(A)} 
For all $n \in \nat$ and $x \in \inspace$ we have:
\begin{align} \metric_\outspace\bigl(\predfn(x) , f(\xi_n^x)  \bigr)
&\stackrel{(i)}{\leq} 
\metric_\outspace\bigl(\predfn(x) , \predfn(\xi_n^x)  \bigr) + \metric_\outspace\bigl(\predfn(\xi_n^x) , f(\xi_n^x)  \bigr)\\
&\stackrel{(ii)}{\leq} 
\metric_\outspace\bigl(\predfn(x) , \predfn(\xi_n^x)  \bigr) + \frac \lambda 2 + \metric_\outspace\bigl(0, \obserr(\xi_n^x) \bigr)
\stackrel{ }{=} 
\norm{\predfn(x) - \predfn(\xi_n^x) }_\infty + \frac \lambda 2 + \obserrbnd\\
&\stackrel{(iii) }{\leq} L(n)
\Metrici{x}{\xi_n^x }^\hexp + \frac \lambda 2 + \obserrbnd \label{ineq:kjshkhd8e}
\end{align}
Here, (i) follows from the triangle inequality, (ii) leverages Lem. \ref{lem:LACKIsampleconsistency} and (iii) follows by H\"older continuity of the predictors (Lem. \ref{lem:LACKIpredHoelder}).

Thus, for $x \in \inspace, n \in \nat$:
\begin{align}
0 \leq \metric_\outspace\bigl(\predfn(x) , f(x)\bigr) &\leq \metric_\outspace\bigl(\predfn(x) ,  f(\xi_n^x)  \bigr) + \metric_\outspace\bigl(f(\xi_n^x) , f(x)\bigr)\\
&\stackrel{(\dagger)}{\leq} (L(n)+ L^*) \metric(x,\xi_n^x )^\hexp +  \frac \lambda 2 +\obserrbnd
\end{align}
where ($\dagger$) follows from (\ref{ineq:kjshkhd8e}) and the presupposed H\"older continuity of $f$. \\

Since by assumption, $ \metric(x,\xi_n^x )^\hexp \leq r_x(n) 
, \forall n$ this implies:
\begin{equation}
\metric_\outspace\bigl(\predfn(x) , f(x)\bigr)  \in \bigl[0, (L(n)+ L^*) r_x(n)^\hexp +  \frac \lambda 2 +\obserrbnd \bigr], \forall n. \end{equation} 
By assumption $ r_x(n), L(n) r_x^\hexp(n) \stackrel{n \to \infty}{\to} 0, \forall x$ and hence, $\metric_\outspace\bigl(\predfn(x) , f(x)\bigr)$ converges to $  [0,\frac \lambda 2 +\obserrbnd], \forall x$ with rate $\varrho_x \leq(L(n)+ L^*) r_x(n)^\hexp $.
%
%
%

\textbf{(B)} Proceeding analogously as before, but utilising uniform convergence with rate $r$, we obtain:  \begin{equation} 
\metric_\outspace\bigl(\predfn(x) , f(x)\bigr)  \in \bigl[0, (L(n)+ L^*) r(n)^\hexp +  \frac \lambda 2 +\obserrbnd \bigr], \forall x \forall n. \end{equation} 
By assumption,  $L(n) r^\hexp(n) \in o(1)$ and thus, $\lim_{n \to \infty} L(n) r(n)^\hexp = 0$. Hence,  $$\errmetric^\infty = \seq{\sup_{x \in I} \metric_\outspace\bigl(\predfn(x) , f(x)\bigr) }{n \in \nat} \stackrel{\varrho}{\convto} [0, \frac \lambda 2 +\obserrbnd] $$ with rate $\varrho $ such that $\varrho(n) \leq (L(n)+ L^*) r(n)^\hexp, \forall n$.

\end{proof}
\end{thm}

Note a necessary condition was that the product of $L(n)$ and the rate was in $o(1)$, that is, vanishing in the limit of $n \to \infty$. A sufficient condition for this to hold is if $L(n)$ is guaranteed to be bounded (assuming the rate is vanishing). Above, we have established a sufficent condition for this (cf. Rem. \ref{rem:convlipconstestimates}): $L(n)$ is bounded as long as parameter $\hestthresh \geq 2 \obserrbnd + q$ for any $q \geq 0$. This yields the following result:

\begin{cor}
With definitions and assumption as in Thm. \ref{rem:convlipconstestimates}, if parameter $\hestthresh $ is chosen to be $2 \obserrbnd + q $ for any $q \geq 0$ then convergence to the ground truth is guaranteed (up to an twice the observational error and a term dependent on $q$). In particular, if the data becomes dense uniformly in $I \subseteq \inspace$ with a rate of $r(n)$  then, for some $\bar L \in [0,L^*]$ and any $n \in \nat$, we have 
\begin{equation}
\sup_{x \in I} \metric_\outspace\bigl(\predfn(x) , f(x)\bigr) \leq (\bar L+ L^*) r(n)^\hexp +  \frac { q} 2 +2 \obserrbnd  \stackrel{n \to \infty}{\convto}  \frac q 2 + 2  \obserrbnd . \end{equation} 
\label{cor:worstcaseconvhoeldertarget}
\end{cor}

Of course in the absence of observational errors, one can choose $\hestthresh = 0$. In this case, the corollary implies that LACKI will learn the ground-truth arbitrarily well in the limit of infinitely dense data.

\begin{remark} [Curse of dimensionality]
Our bounds rely on the proximity (expressed by the rate functions) of the query input to the previously observed data.
Refer to Thm. \ref{thm:convergenceifboundedconstandsamplecons_LACKI}.
Roughly speaking, for a particular query input $x$, our guarantee in (A) asserts that the closer the query is to the previously seen data, the better the confidence in prediction accuracy. In (B) this is extended to a worst-case statement implying that the smaller the worst-case proximity of the data to any query in $I$, the smaller the worst-case prediction error can be. 
Unfortunately, this worst-case proximity and therefore, the prediction error bound, is subject to the \emph{curse of dimensionality}. That is, the number of samples necessary to guarantee a desired reduction in worst-case prediction uncertainty will inevitably have to scale exponentionally with the dimensionality of the space. A manifestation of this fact can be seen in Sec. \ref{sec:probconv_LACKI} where we give a sample complexity bound for uniformly distributed input samples. 
\end{remark}

Having established that our LACKI rule can learn any H\"older function with any H\"older constant, we will now attend to extend the results to non-H\"older functions. In preparation of the necessary derivations we will first rehearse universality and H\"older properties of radial basis function networks. 

Park and Sandberg derived universal approximation guarantees for radial-basis function networks \cite{Park1991}. In particular, on page 252 in their article the authors make an assertion that translates to our notation as follows: 

\begin{lem}[Expressiveness of RBFNs] \label{lem:RBFNunifapproxcompact} Assume $\inspace \subseteq \Real^d$ is compact. Given parameter  vector $\theta := (w_1,\ldots,w_m,\sigma_1,...,\sigma_m,c_1,\ldots, c_m)$ and kernel function $K: \inspace \to \outspace $ let $\beta(\cdot;\theta ) = \sum_{i=1}^m w_i \, K(\frac{\cdot - c_i}{\sigma_i} )  $ denote a radial basis function network (RBFN). Assume $K: \Real^d \to \Real$ is continuous and has non-vanishing integral, i.e. $\int_{\Real^d} K(x) \d x \neq 0$.
Then, the set $S_K:= \{ \beta(\cdot; \theta) \vert  m \in \nat, \theta \in \Real^{3m}  \}$ of all RBFNs is uniformly dense in the set $C(\inspace)$ of continuous functions on compact domain $\inspace$. That is,  $\forall f \in C(\inspace) \forall \epsilon >0  \exists m, \theta \in \Real^{3m} : \sup_{x \in \inspace}{ \abs{f(\cdot) -\beta(\cdot;\theta)  } } <\epsilon $. 
\end{lem}

\begin{remark} \label{rem:LipconstofRBFN}
We note that, for any finite-dimensional parameter $\theta$, any RBFN $\beta(\cdot;\theta)$ is Lipschitz continuous as long as the kernel $K$ is. This can be seen by applying Lem. \ref{lem:Hoeldarithmetic} which allows us to conclude that the Lipschitz constant of RBFN $\beta(\cdot;\theta ) = \sum_{i=1}^m w_i \, K(\frac{\cdot - c_i}{\sigma_i} ) $ is given by $L_\beta = \sum_{i=1}^m \abs{\frac{w_i}{\sigma_i}} L_{K}$ where $L_K \in \Real_{\geq 0}$ denotes a Lipschitz constant of $K$. By the same Lemma it is easy to see that choosing the Gaussian kernel for $K$ satisfies both the Lipschitz requirement as well as the integrability requirements of Lem. \ref{lem:RBFNunifapproxcompact}. As a by-product this means that on a compact support, any continuous function can be approximated by some Lipschitz function with arbitrarily small, positive worst-case error $\epsilon >0$. Note, it may well be the case that the Lipschitz constant of the approximator grows with decreasing approximation error bound $\epsilon$. We consider this to be inevitable when the approximated function is not Lipschitz.
\end{remark}

Harnessed with these preparatory statements we can move on to show that the LACKI rule can be set up to learn any continuous function up to arbitrary low error.
\begin{thm}[Universality of LACKI]
\label{thm:LACKIuniversality}
Assume we are given a sequence $\seq{\data_n}{n \in \nat}$ of samples with observational errors bounded by $\obserrbnd \in \Real_{\geq 0}$. We set the parameters of the LACKI rule to $\lbf =-\infty,\ubf =\infty, \underline L =0$ and $\hestthresh := 2 \bar e +q $ for some $q >0$. In this theorem, we assume that the set of interest $I \subseteq\inspace$ is compact.
\textbf{Then, we have:}

The LACKI rule as per Def. \ref{def:LACKI} is a universal approximator in the following sense:
If the sequence of input grids $\seq{\grid_n}{n \in \nat}$ relative to $I$ (uniformly) then the sequence of predictors $\seq{\predfn}{n \in \nat}$ computed by the LACKI rule (uniformly) converges to any continuous target $f : \inspace \to \Real$ up to error $2 \obserrbnd + \frac{3q}{2}$.
That is, the following holds true:

\begin{itemize}
\item (I) Let $x \in I$. If $\exists r_x \in o(1): (\grid_n) \stackrel{ r_x}{\bd} x$ then $\exists C \in \Real : \seq{\Metrico{\predfn(x)}{f(x)}}{ } \stackrel{C r_x^\hexp}{\convto}[2\obserrbnd + \frac{3q}{2}]$.
\item (II) If $\exists r \in o(1): (\grid_n) \stackrel{ r}{\bdu} I$ then $\exists C \in \Real : {\errmetric^\infty}{ } \stackrel{C r^\hexp}{\convto}[2\obserrbnd + \frac{3q}{2}]$.
\end{itemize}
\end{thm}
\begin{proof}
We choose any parameter $\hestthresh =  2 \obserrbnd + q$ with $q >0$. As observed in Rem. \ref{rem:LipconstofRBFN}, Lem. \ref{lem:RBFNunifapproxcompact} allows us to infer that there exists a Lipschitz function $h$ that approximates the target with worst-case error of at most $\frac q 2$. That is, $\sup_{x \in \inspace} \metric_\outspace\bigl(h(x), f(x) \bigr) \leq \frac q 2$. (Also, note Lipschitz continuity implies H\"older continuity for any H\"older exponent $\hexp \in (0,1]$, and hence, $h \in \hoelset {L_h}{ }{\hexp}$ for some $ L_h \in \Real_{\geq 0}$.)

Consequently, there exists a function $\phi':\inspace \to \outspace$ with $\sup_x \Metrico{ 0}{\phi'(x) } \leq \frac q 2$ accounting for the discrepancy between the H\"older function $h$ and the target $f$: $f = h+ \phi'$. 

Furthermore,we define $\phi$ to be the bounded observational noise. Hence, we have $\tilde f = f+ \phi$ and $\sup_x \Metrico{0}{\phi(x)} \leq \obserrbnd$.
Combining both functions into $\psi := \phi+\phi'$, we can write $\tilde f = h + \psi$ with $\sup_x \Metrico{0}{\psi(x)}\leq \frac q 2 + \obserrbnd =: \bar \nu$.

This can be interpreted as follows:
Instead of viewing the given sample as being generated by target $f$ (with some observational error $\phi$) we can view the sample as being generated by the H\"older function $h$ corrupted by the extended ``observational noise'' $\psi$ accounting for both the original observational error and the discrepancy between the target and H\"older function $h$.
This gives us a reduction to the case of learning H\"older functions with observational error bounded by $\bar \nu$. Firstly, we note that $\hestthresh = 2 \obserrbnd +q =2 \bar \nu$ (which entails that the sequence $\seq{L(n)}{n \in \nat}$ is bounded by some constant $\bar L  =\sup_{x,x', \metric(x,x') >0} \frac{\metric_\outspace(h(x), h(x'))   -q}{\metric^\hexp(x,x')} \leq L_h$). 
Linking to Thm. \ref{thm:convergenceifboundedconstandsamplecons_LACKI}, we get all the desired statements with regard to learning $h$. These can easily be converted into statements about learning $f$ by adding the worst-case difference $\frac q 2$ between $f$ and $h$ to all error bounds. 
For example, leveraging $\sup_x \Metrico{ 0}{\phi'(x) } \leq \frac q 2$ and $\lambda = 2 \obserrbnd + q$ and going through analogous steps as in the previous theorem we obtain:   
\begin{align}
\metric_\outspace(\predfn(x), f(x) ) &= \metric_\outspace(\predfn(x), h+\phi' (x) ) 
\leq   \metric_\outspace\bigl(\predfn(x) , h(x)\bigr) + \metric_\outspace(0,\phi'(x) ) \\
&\leq (\bar L+ L_h) \metric(x,\xi_n^x )^\hexp + \frac \lambda 2 +\bar \nu + \frac q 2 \\
&\leq (\bar L+ L_h) \metric(x,\xi_n^x )^\hexp  +2 \obserrbnd + \frac {3q}{2} \label{ineq:euyiweeh}
\end{align}
where  $\xi^x_n := \arg\inf_{s \in \grid_n} \metric(x,s)$ denotes a nearest neighbour of $x$ in the input sample $\grid_n$.

So, convergence (pointwise or uniform) of the grid to the input space with a rate of at most $r(n)$ implies that the right-hand side of (\ref{ineq:euyiweeh}) and hence, the prediction error,  
converges (pointwise or uniformly) to the interval $[0,2 \obserrbnd + \frac {3q}{2}] $ with a rate of at most $(\bar L + L_h) r^\hexp(n)$ as $n \to  \infty$.

\end{proof}


\subsubsection{Convergence in probability with uniformly distributed inputs}
\label{sec:probconv_LACKI}
Above we have given guarantees relative to the deterministic convergence rates of the input sample to the domain.
In this subsection, we shall study probabilistic convergence rates as a function of the sample size in situations where the sample is obtained by drawing inputs independently from a uniform probability distribution on $I=\inspace := [0,1]^d$. 

We can show that the worst-case prediction error given by $\sup_{x \in \inspace}\Metrico{\predfn(x)}{ f(x)}$ vanishes (up to the usual worst-case bounds in the presence of observational errors) in probability for canonical input-space metrics:  

%
%

\begin{thm} Let $\inspace = [0,1]^d $ be the domain of target function $f \in \hoelset {L^*} { } \hexp $. Assume the input data $\grid_n = \{s_1,\dots,s_n\}$ contains $n$ data sample inputs which are drawn independently at random from a uniform distribution over $\inspace$. Furthermore, assume there are no observational errors, i.e. $\obserrbnd =0$, and, that $\metric(x,x') = \norm{x-x'}_\infty, \forall x,x' \in \inspace$. The worst-case error of our LACKI predictor vanishes in probability. 

That is, 
$$\forall \epsilon >0 \forall \delta \in (0,1) \exists N \in \nat \forall n \geq N : \Pr[ \sup_{x \in \inspace} \Metrico{\predfn(x)}{f(x)} >\epsilon] \leq \delta.$$
In particular, for all $\delta \in (0,1)$ we have 
$\Pr[ \sup_{x \in \inspace} \Metrico{\predfn(x)}{f(x)} >\epsilon] \leq \delta$
\begin{enumerate}
\item  for any $\epsilon \geq 2 L^*$, provided that $n \geq 1$;
\item for any $\epsilon < 2 L^*$, provided that $n \geq N := \ceil{\frac{ \log(\delta \, 2^{-kd}  )}{\log(1- 2^{-kd})} }$ with $k= \ceil{\frac{\log(\epsilon^{-1}2 L^*)}{\log 2}}$.
\end{enumerate}
\begin{proof}
Let $r_n := \sup_{x \in \inspace} \min_{s \in \grid_n} \metric(x,s) = \sup_{x \in \inspace} \min_{s \in \grid_n} \norm{x-s}_\infty \leq 1$ and let 

$P_n^\epsilon := \Pr[ \sup_{x \in \inspace} \Metrico{\predfn(x)}{f(x)} >\epsilon]$ which we intend to bound from above.
Remember, from Cor. \ref{cor:worstcaseconvhoeldertarget} $\sup_x \Metrico{\predfn(x)}{f(x)} \leq 2 L^*  r_n $. Hence, for $\epsilon \geq 2 L^*$, $P_n^\epsilon \leq 0, \forall n \in \nat$.

So, it suffices to focus on the case where $\epsilon < 2 L^*$. Now, $\sup_x \Metrico{\predfn(x)}{f(x)} \leq \epsilon  $ is implied by $\sup_x \Metrico{\predfn(x)}{f(x)} \leq 2 L^* r_n $ provided that $  r_n  \leq \frac {\epsilon}{2 L^*}$. So, we define an event $E_n$ that ensures $r_n$ satisfies the latter inequality with a probability that grows as $n$ increases.
To this end, we introduce a partition of the domain into $m$ hyper-rectangles $H_1,...,H_m$ of equal size, each having edge length $l_k=\frac 1 {2^k}$ where $k$ is a natural number such that $l_k \leq  \frac {\epsilon}{2 L^*} $. As a valid choice, we set $k:= \ceil{\frac{\log(\epsilon^{-1}2 L^*)}{\log 2}}$. Note, $\Pr[s_i \in H_j] = l_k^d = \frac{1}{2^{dk}}$.  By construction, in the event that each hyper-rectangle ends up containing at least one sample input of $\grid_n$, we have  $r_n \leq \frac {\epsilon}{2 L^*}$. 
We define the complement of this event as $\bar E_n := \{(s_1,...,s_n)  \in \inspace^n | \exists j \in \{1,...,m\} \forall i \in \{1,...,n\}: s_i \notin H_j  \}$. Let  $W:= \{ s=(s_1,...,s_n)  | \sup_x \Metrico{\predfn(x)}{f(x)} > \epsilon \}$ be the event that the sample inputs are located in such a fashion that they give rise to a worst-case error larger than $\epsilon$.  We have: $s \notin \bar E$ implies that $r(n) \leq \frac {\epsilon}{2 L^*} $ which in turn implies $\sup_x \Metrico{\predfn(x)}{f(x)} \leq \epsilon $, i.e. that $s \notin W$. Hence, $W \subseteq \bar E_n$ and thus, $ P_n^\epsilon =\Pr[W] \leq \Pr[\bar E_n]$.
So, to bound $P_n^\epsilon$ from above it suffices to bound $\Pr[\bar E_n]$ from above which we will do next: We can employ the union bound, utilise that $m = 2^{kd}$ and the fact that the $s_i$ are drawn i.i.d. from a uniform to see that $\Pr[\bar E_n] \leq \sum_{j=1}^m \prod_{i=1}^n \Pr[s_i \notin H_j] = 2^{kd}  (1-\frac{1}{2^{dk}})^n \stackrel{n \to \infty}{\convto} 0$ which shows the main statement of the theorem. To find an $n$ sufficently large to ensure $\Pr[W] \leq \delta$ we consider the inequality $2^{kd}  (1-\frac{1}{2^{dk}})^n \leq \delta$. Taking the $\log$ on both sides and rearranging yields the sufficient condition: $n \geq \frac{ \log(\delta \, 2^{-kd}  )}{\log(1- 2^{-kd})}$. 
\end{proof}
\end{thm}

\subsubsection{Some guarantees for online learning}

In the theorems above, we considered the worst-case asymptotics for the case where the data becomes dense in the domain. Here the error was evaluated on the entire input domain. By contrast, we will now consider an online learning setting where we incrementally get to observe samples along the trajectory of inputs $\seq{x_n}{n \in \nat }$ and are interested in the long-term one-step-lookahead prediction errors on this trajectory.
That is, we are interested in the evolution of prediction errors $\Metrico{\predfn(x_n)}{f(x_n) }$
where the predictor $\predf_{n}(\cdot)$ is based on $\data_{n} = \data_{n-1} \cup \{ \bigl(\state_{n-1}, \tilde f(\state_{n-1}) \bigr)\}, \forall n >1 $. 

We will show that this error trajectory vanishes (up to observational errors), provided that the input sequence $\seq{x_n}{n \in \nat}$ is bounded.

In preparation of these considerations, we will establish the following facts:

\begin{lem} 
Assume we are given a trajectory $\seq{x_n}{n \in \nat}$ of inputs with $x_n \in \inspace$ where input space $\inspace$ can be endowed with a shift-invariant measure. Furthermore, assume the sequence  is bounded, i.e.  
$\metric_\inspace(x_n,0) \leq \beta$ for some $\beta \in \Real_+$ and all $n \in \nat$.
Finally assume the inputs of the available data coincide with the complete history of past inputs, i.e. $G_n = \{ x_i | i \in \nat, i < n\}$.
Then we have: \[ \dist(G_n,x_n) = \min\{\metric_\inspace(g,x_n) | \, g \in G_n\} \stackrel{n \to \infty}{\longrightarrow} 0.\]
\begin{proof}
The intuition behind the following proof is that if the distances were not to converge, there was an infinite number of disjoint balls around the input points that summed up to infinite volume. This however, would be a contradiction to the presupposed boundedness of the sequence.
We formalise this intuition as follows:
We can rephrase the desired convergence statement as 
\begin{equation}
\forall \epsilon > 0 \exists n \in \nat \forall m > n : \dist(x_{m}, G_{m}) \leq \epsilon.
\end{equation} 
For contradiction, suppose
\begin{equation}
\exists \epsilon > 0 \forall n \in \nat \exists m(n) > n : \dist(x_{m(n)}, G_{m(n)}) > \epsilon.
\end{equation} 
Hold such an $\epsilon >0$ fixed. For $n \in \nat$, we can define mapping $m: \nat \to \nat$ such that  $ m(n)>n$ and $  \dist(x_{m(n)}, G_{m(n)}) > \epsilon.$
By definition of $G_{m(n)} =\{ x_i | i < m(n)\} $ we have:
\eqn{eq:i34kjjk3}{\forall i < m(n) : \metric_\inspace(x_{m(n)},x_i) > \epsilon.}
This tells us that we can define a subsequence $(\xi_n)$ of $(x_n)$, where $\xi_n = x_{\phi(n)}$ for some strictly monotonically increasing mapping $\phi:\nat \to \nat$, such that 
$ \forall i \neq j: \metric_\inspace(\xi_i,\xi_j) > \epsilon$
and hence, the $\frac \epsilon 2 -$balls around the subsequence members are disjoint: $$\forall i \neq j: \ball{\frac \epsilon 2}{\xi_i}  \cap  \ball{\frac \epsilon 2}{\xi_j} = \emptyset.$$

For a given choice of mapping $m :\nat \to \nat$, we can inductively construct such a subsequence as per: $ \xi_{n} = x_{\phi(n)}$ where $\phi(1) = 1, \phi(n+1) = m(\phi(n)), \, (n>1)$.

Next, let $C_n := \bigcup_{i =1}^n \ball{\frac \epsilon 2}{\xi_i} $ be the union of all $\frac \epsilon 2$-balls around each point of the first $n$ elements of the subsequence. Moreover,  define $\bar I = \bigcup_{n \in \nat} \ball{\frac \epsilon 2}{x_n} $.
By definition, each $x_n$ (and thus, each $\xi_n$) is contained in $\bar I$.
Since sequence $(x_n)_{n \in \nat}$ is bounded, $\bar I $ has a finite volume relative to some positive, shift-invariant measure $\mu$. I.e. $\mu(\bar I) < \infty$ (e.g. choose the Lebesgue measure for $\mu$). 
Since $C_n \subseteq \bar I$ and owing to the disjointness of the $\frac \epsilon 2$- balls of the subsequence,  we have :
 $\mu(C_n) = \sum_{i=1}^n \mu(B_i) \leq \mu(\bar I)< \infty$ where $B_i := \ball{\frac \epsilon 2}{\xi_i}$. Owing to the assumed shift-invariance, we can assign the same measure $M$ to each ball, i.e. $M:=\mu(B_1) =...= \mu(B_n)\forall n \in \nat$. Thus, $\mu(C_n) = n M$.
Define $q:= \ceil{\frac{\mu(\bar I)}{M}} \in \nat$. This is an upper bound on the number of disjoint balls of measure $M$  that can be contained in $\bar I$. Intuitively, since $\mu(\bar I)$ is finite, $\bar I$ cannot contain an infinite number of non-intersecting balls. 
Concretely, for  $n > q$ we have:
\begin{align}
\mu(\bar I) &\geq \mu(C_n) = M \, n  > M \, q = M \,  \ceil{\frac{\mu(\bar I)}{M}} \geq \mu(\bar I).
 \end{align}
Obviously, $\mu(\bar I) > \mu(\bar I)$ is a false statement, establishing the desired contradiction.

\end{proof}
\label{lem:bndseq_entails_distgridvanish}

\end{lem}  

\begin{thm}
Assume that, for some $q\geq0$, we chose $\hestthresh = 2 \obserrbnd + q$ in our LACKI prediction rule. And, assume that the target $f$ is H\"older continuous up to some error level $\bar E_h$. That is, $f = \phi + \psi$ with $\phi \in \hoelset {L^*} { } \hexp$ and a function $\psi$ such that $\sup_x \metric_\outspace\bigl(0,\psi(x)\bigr) \leq \bar E_h \in \Real$.

Assume we are given a trajectory $\seq{x_n}{n \in \nat}$ of inputs that is bounded, i.e. where 
$\metric(x_n,0) \leq \beta$ for some $\beta \in \Real_+$ and all $n \in \nat$.
Furthermore, assume $\data_{n+1} = \data_n \cup \{ \bigl(x_n, \tilde f(x_n)\bigr) \}$ and thus, $\grid_n = \{ x_i | i \in \nat, i < n\}$.
Then the prediction error on the sequence vanishes up to the level of sample-consistency and H\"older continuity in the following sense:
 \[\metric_\outspace\bigl(\predfn(x_n),f(x_n) \bigr) \stackrel{n \to \infty}{\longrightarrow} [0,\frac q 2 + 2  \obserrbnd  + 2 \bar E_h].\]
In particular, in case the observations are error-free ($\tilde f = f$) and assuming the target is H\"older continuous then, when choosing $\hestthresh = 0$, the prediction error is guaranteed to vanish. That is,
\[\metric_\outspace\bigl(\predfn(x_n),f(x_n) \bigr) \stackrel{n \to \infty}{\longrightarrow}0.\]

\begin{proof}

Let $\xi_n  \in \argmin_{g \in \grid_n} \metric(x_n,g)$ denote the nearest neighbour of $x_n$ in $\grid_n = \{x_1,...,x_{n-1}\}$.

Since sequence $(x_n)$ is bounded, Lem. \ref{lem:bndseq_entails_distgridvanish} is applicable and hence: (i) $\lim_{n \to \infty} \metric (x_n,\xi_n) = 0$.

From Lem. \ref{lem:LACKIsampleconsistency} we conclude 
$\metric_\outspace\bigl(\predfn(\xi_n) ,  f(\xi_n)  \bigr) \leq 2\obserrbnd +  \frac q 2$. Hence, appealing to the triangle inequality, we see that 
(ii) $\metric_\outspace\bigl(\predfn(x_n) ,  f(\xi_n)  \bigr) \leq \metric_\outspace\bigl(\predfn(x_n) ,  \predfn(\xi_n) \bigr) + 2 \obserrbnd + \frac q 2$.

Moreover we note that the predictors $\predfn$ have H\"older constants $L(n)$ and that the $L(n)$ are bounded from above by some $\bar L \in \Real$. Thus,  $(iii)$ $\exists \bar L \in \Real\forall n \in \nat : \predfn \in \hoelset {\bar L} { } \hexp$. 
  
In conclusion,
$0\leq\metric_\outspace\bigl(\predfn(x_n) , f(x_n)\bigr) \leq \metric_\outspace\bigl(\predfn(x_n) ,  f(\xi_n)  \bigr) + \metric_\outspace\bigl(f(\xi_n) , f(x_n)\bigr) \stackrel{(ii)}{\leq} \metric_\outspace\bigl(\predfn(x_n) , \predfn(\xi_n) \bigr) + 2 \obserrbnd + \frac q 2 + \metric_\outspace\bigl(f(\xi_n) , f(x_n)\bigr) \leq \metric_\outspace\bigl(\predfn(x_n) , \predfn(\xi_n) \bigr) +2 \obserrbnd + \frac q 2 + \metric_\outspace\bigl(\phi(\xi_n) , \phi(x_n)\bigr) + 2 \bar E_h 
\newline
\stackrel{(iii)}{\leq} (\bar L+ L^* ) \metric(x_n,\xi_n )^\hexp + 2 \obserrbnd + \frac q 2 + 2 \bar E_h  \stackrel{n \to \infty}{\longrightarrow} 2 \obserrbnd + \frac q 2 + 2 \bar E_h $.
\end{proof}
\label{thm:vanisishingseqprederr_LACKI}
\end{thm} 

\subsubsection{Computational complexity}
The computational complexity for computing the parameter update $L(n+1)$ based on $L(n)$, the pre-existing data $\data_n$ of size $N_n$ and a newly arriving sample input is in $\mathcal O(N_n m) $ where $m$ is the effort for evaluating the pseudo-metric. Typically, $m$ will scale linearly with the dimensionality of the input space. Therefore, online updates cost training time that will be linear in the number of existing training data and input dimensionality. In batch training, for a batch of $N$ samples, computation of the estimate will require an effort in $\mathcal O(N^2 D)$. 
Once the parameter $L(n)$ is computed, the effort for evaluating $\predfn(x)$ is linear in the number of samples and, again, typically linear in the input and output space dimensionality. 

However, it should be noted though that generalised nearest neighbor techniques can be utilised to reduce the prediction effort to expected logarithmic effort in the sample size (see \cite{Beliakov2006}). Devised for standard Lipschitz interpolation, this approach could be readily applied to our LACKI inference rule.

\subsection{Tests on an artificial regression problem}
\label{sec:regr_benchmarks}
 \begin{figure*}
        \centering
		    \includegraphics[ height = .25 \textheight,width=.999\textwidth]{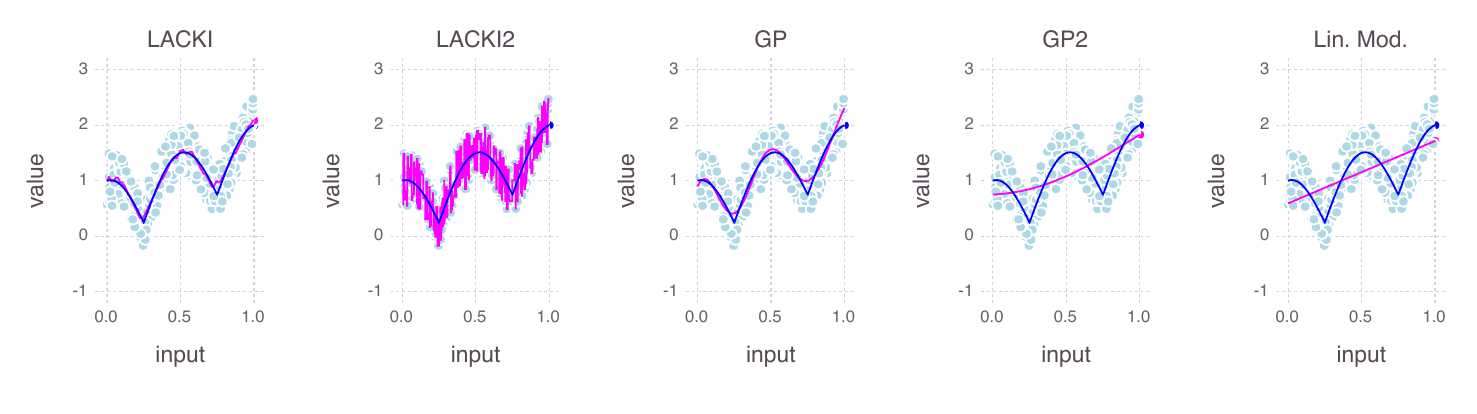}
   \caption{\textbf{Exp. 0.} The predictors of the various regression methods on target function $f_1$ for $d=1$ and training sample size $N_n = 500$. Training examples plotted as light blue dots, the graph of the target function is plotted in dark blue and the predictions are plotted in magenta.}
	 \label{fig:bench011dimregr}
\end{figure*}

 \begin{figure*}
        \centering
  \subfigure{
    \includegraphics[ height = .19 \textheight, width=.32\textwidth]
								{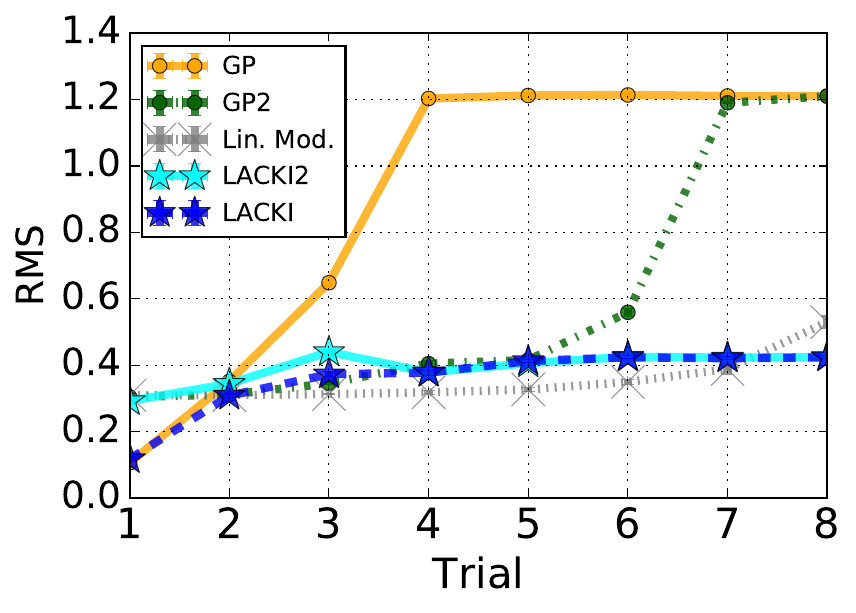}
  } 
  \subfigure{
    \includegraphics[ height = .19 \textheight, width=.31\textwidth]
								{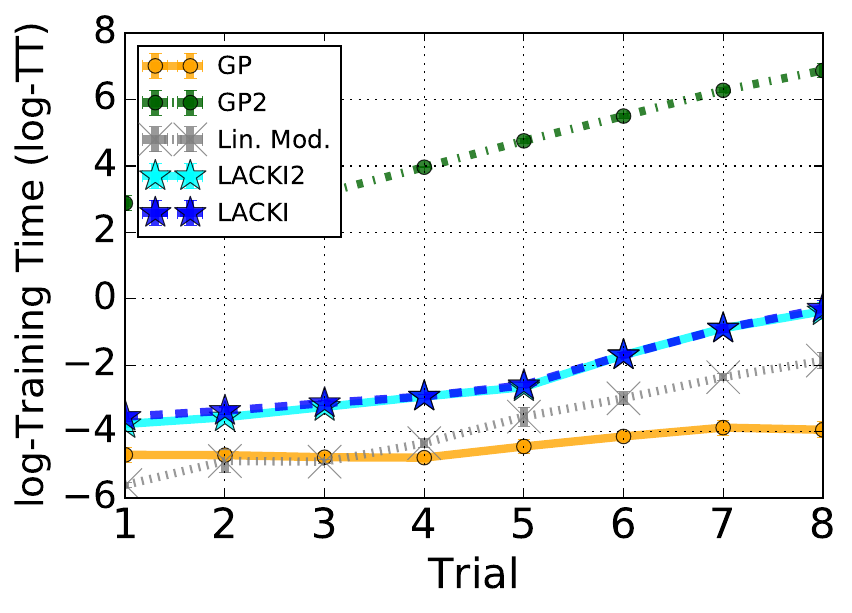}
								  } 
    \subfigure{
    \includegraphics[ height = .19 \textheight, width=.31\textwidth]
								{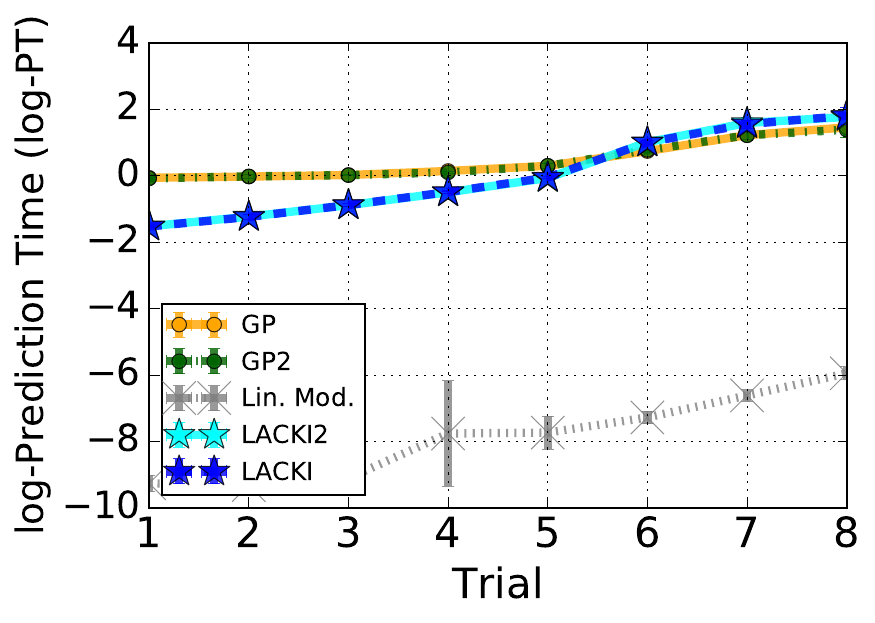}
  } 
  \caption{\textbf{Exp. 0.} Comparisons on target function $f_1$ over a range of different trials. In the $n$th trial, test function 1 was sampled uniformly at random on the domain $2^n$-dimensional input space domain $\inspace =[0,1]^{2^n}$ with \emph{fixed training data size} $\abs{\data_n} = 500$ but varying input space. Depicted are the measured means and standard deviations for each method and each trial. Note how LACKI's predictive performance degrades comparatively mildly with increasing dimensionality on the given test function.}
	 \label{fig:bench00}
\end{figure*}	

Having a established worst-case prediction error bounds of our LACKI method, we now aim to illustrate the benefits and shortcomings of our approach in a number of regression problems. In order to assess the predictive performance accurately, we need access to the ground truth. Therefore, we restrict our benchmarks to artificial data for the time being. 
We conducted three different experiments. In each of the experiments we varied the conditions over the course of randomised trials as follows:

\begin{itemize}
\item \textbf{Exp. 0}: Here, we investigated the regression performance on the target function $f_1: [0,1]^d \to \Real, x \mapsto \abs{\cos(2\pi x_1)}+x_1$ for increasing input-space dimensionality $d$. The observations were perturbed by i.i.d. uniform noise sampled from the interval $[-.5,.5]$.  The training sample size was fixed to $N_n = 1000$ examples. In trial $n$, the input space dimensionality was set to $d=2^n$.
\item \textbf{Exp. 1}: 
A repetition of the setup of Exp.1 but for fixed input space dimensionality $d=2$ and increasing sample size. 
In trial $n$ the training data size was set to $N_n = 2^n+1$.
\item \textbf{Exp. 2}: A repetition of the setup of Exp. 2, but for the benchmark target (considered in \cite{Beliakov2006}) $f_2: x \mapsto \sin(x_1) \sin(x_2)+0.05 \bigl(\sin(5x_1) \sin(5 x_2)\bigr)^3$ under no observational noise.
\end{itemize}
In each trial of each experiment, we recorded the following performance measures:
\begin{enumerate}
	\item  For a predictor $\predf$, on a set $\inspace_{test}$ of test inputs we recorded the empirical root-mean-square error \textbf{RMS}$ =  \sum_{x \in \inspace_{test}} \abs{f(x) -\predf(x)}^2 $ \newline as well as  maximal prediction error 
	 \textbf{ME} $= \max_{x \in \inspace_{test}} \abs{f(x) -\predf(x)}$.
	\item The log of the run time measurement (in sec.) of training the predictor (\textbf{log-TT}).
	\item The log of the run time measurement (in sec.) of computing the predictions of the random test inputs divided by the number of test inputs (\textbf{log-PT}). 
\end{enumerate}
Here, the test sample inputs in $\inspace_{\text{test}}$ were drawn i.i.d. at random from a uniform distribution over the domain $\inspace$.
Therefore, the pertaining performance measures were random variables. In each trial, we have obtained a sample (of size 30) of these random variables and recorded their empirical means and standard deviations for each trial and for the following regression techniques:
\begin{itemize}
\item \textbf{LACKI}: Our LACKI method with parameter choice set to $\hestthresh =1$. 
\item \textbf{LACKI2}: Our LACKI method as above, but with parameter choice $\hestthresh =0$. 
\item \textbf{GP}: A Gaussian process  \cite{GPbook:2006} with fixed covariance function $k(x,x';\theta) = \theta_1 \exp( \frac{\norm{x-x'}^2}{2 \theta_2})$. We determined the parameters manually to give good results on Exp. 1 for $d=1$. This tuning process resulted in the choice of $\theta = (1,\frac 1 4)$ with observational noise variance $\frac 1 {12}$ (the resulting predictor for a 1-dimensional data set is depcited in Fig. \ref{fig:bench011dimregr}). The predictor was chosen to coincide with the mean function of the posterior process.
\item \textbf{GP2}: A GP with hyper-parameters determined by following the standard approach of maximising the marginal log-likelihood of the data \cite{GPbook:2006}. Optimisation was done without restarts employing BFGS. The optimiser was started with initial hyper-parameters set to $\theta = (1,1)$ and observational noise parameter being initialised with $\frac{\obserrbnd}{2}$.
\item \textbf{Lin. Mod.}: A linear regression model fitted with the least-squares method.
\end{itemize}

The code was implemented in pure \emph{Julia 0.4.7} with the GPs making use of the library \emph{GaussianProcesses.jl}. The code was executed on a 2015-MBP furnished with i7 processors and 16 GB RAM running OS X 10.11.6.

\begin{figure*}
        \centering

  \subfigure{
    \includegraphics[ height = .19 \textheight,width=.32\textwidth]
								{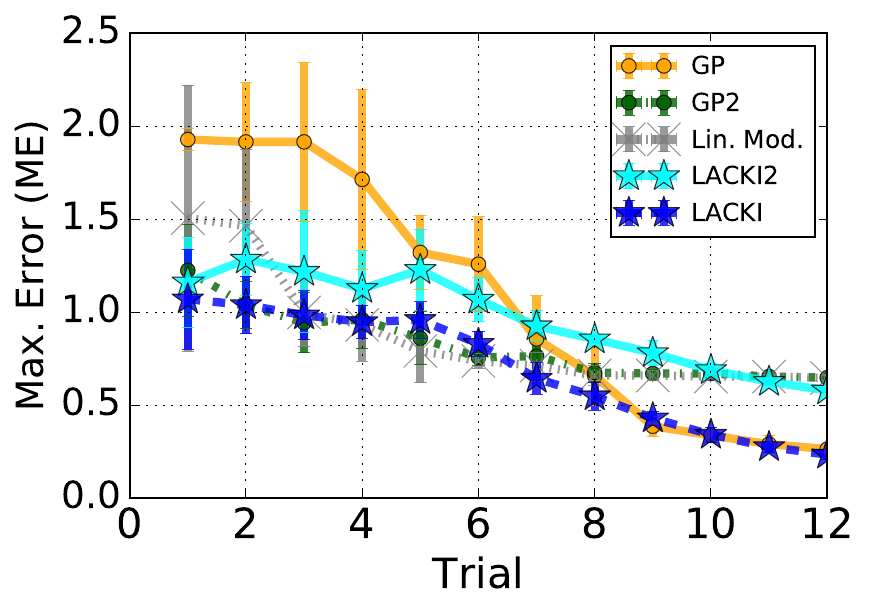}
  } 
  \subfigure{
    \includegraphics[  height = .19 \textheight,width=.31\textwidth]
								{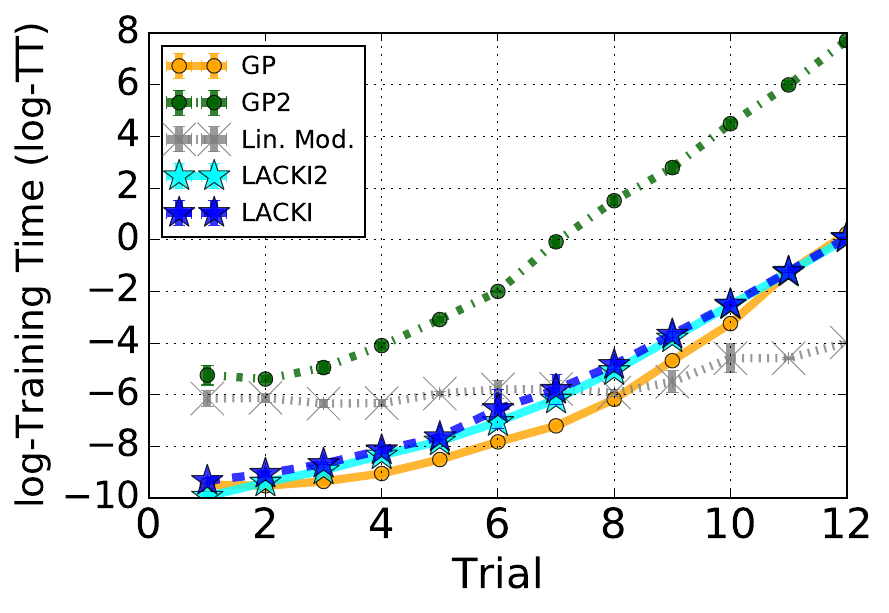}
								  } 
    \subfigure{
    \includegraphics[ height = .19 \textheight,width=.31\textwidth]
								{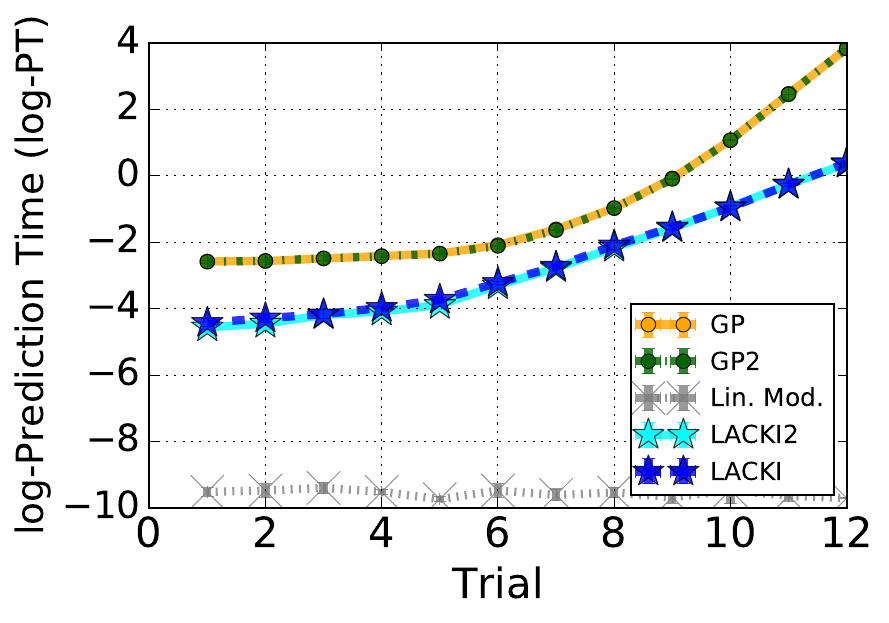}
  } 
   \caption{\textbf{Exp. 1.} Comparisons on target function $f_1$ over a range of different trials. In the $n$th trial, the noise-perturbed target function was sampled uniformly at random on the domain $\inspace =[0,1]^2$ with \emph{varying training data size} $\abs{\data_n} = 2^n+1$. The plots depict measured means and standard deviations over 30 repetitions for each trial.  Prediction error measures were estimated based on $25000$ test samples drawn independently from the domain. Run times for training and prediction are depicted on a log-scale. We note that LACKI, with correct noise parameter $\lambda = 1$ overall outperformed the other methods in terms of prediction accuracy.  }
	 \label{fig:bench01}
\end{figure*}

\begin{figure*}
        \centering

  \subfigure{
    \includegraphics[ height = .19 \textheight,width=.32\textwidth]
								{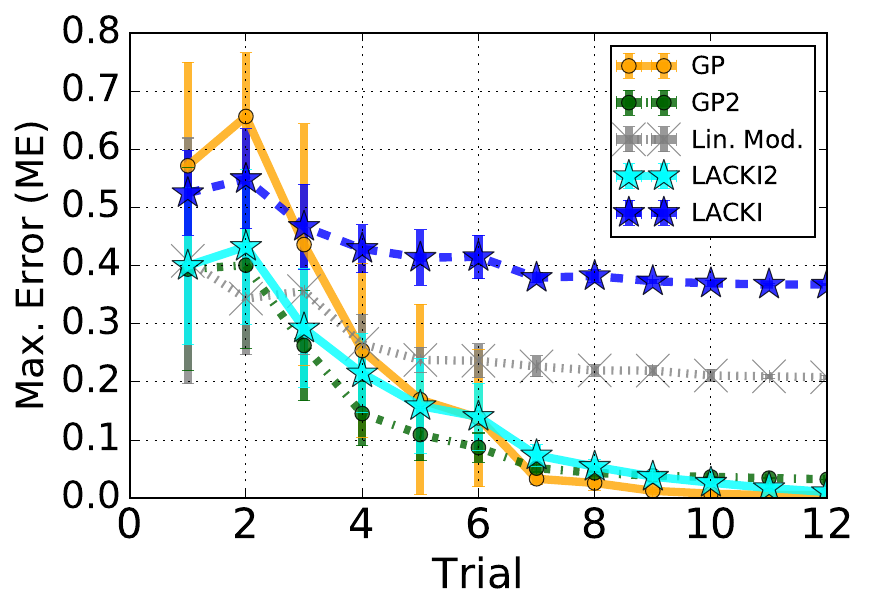}
  } 
  \subfigure{
    \includegraphics[ height = .19 \textheight,width=.31\textwidth]
								{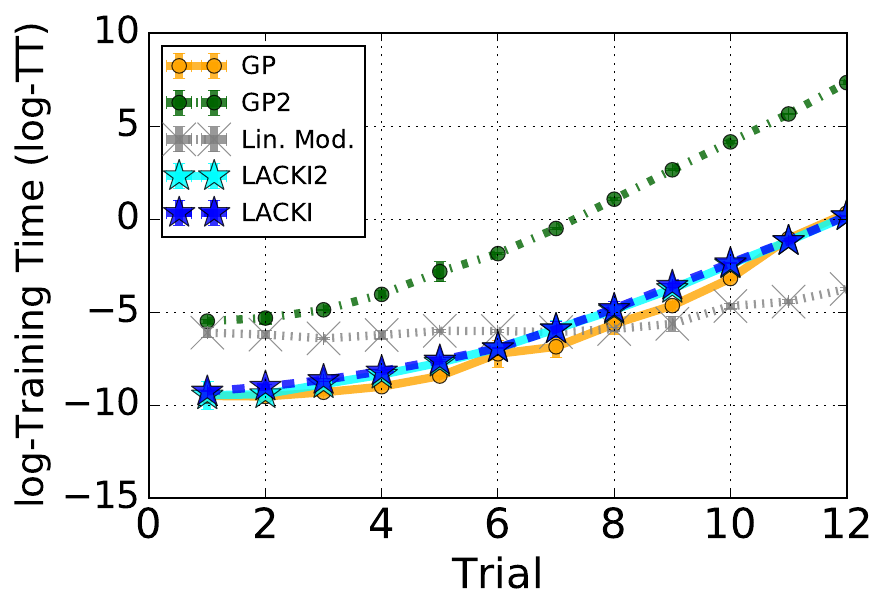}
								  } 
    \subfigure{
    \includegraphics[ height = .19 \textheight,width=.31\textwidth]
								{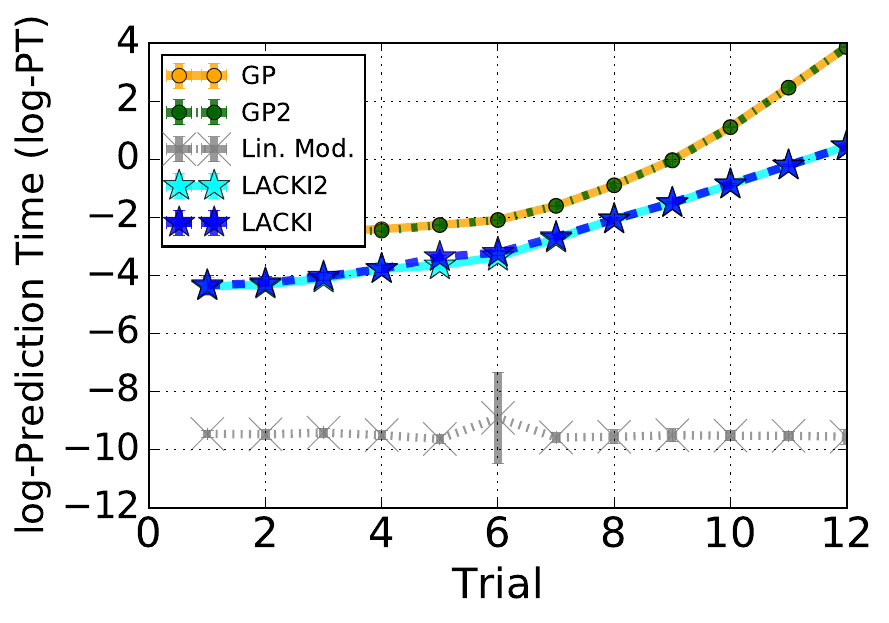}
  } 
   \caption{\textbf{Exp. 2.} Repetition of Exp. 1 but with target function $f_2$ and no observational noise.  In the $n$th trial, the target was sampled uniformly at random on the domain $\inspace =[0,1]^2$ with \emph{varying training data size} $\abs{\data_n} = 2^n+1$.  Prediction errors were estimated based on $25000$ test samples drawn independently from the domain.  Run times for training and prediction times are depicted on a log-scale.  LACKI, having a now falsely set noise parameter $\hestthresh = 1$, did not achieve good predictive results. By contrast, both GPs, as well as LACKI2 (with the correctly set parameter $\hestthresh =0$) managed to learn the target accurately with increasing data. However, in comparison to LACKI2, we note that the GPs exhibited much higher variability in performance for lower sample sizes and that GP2, which tended to have slightly lower prediction error, achieved this at the expense of substantially higher computational training effort.  }
	 \label{fig:bench02}
\end{figure*}

\textbf{Discussion:}
The results of the experiments are depicted in Fig. \ref {fig:bench00} - Fig. \ref{fig:bench02}.
We note that, when the noise hyper-parameter $\hestthresh$ was set correctly to twice the level of observational error (e.g. LACKI in Exp. 0 and Ex. 1, and LACKI 2 in Exp.2), our approach yielded good predictive performance that outperformed the GPs considering both prediction accuracy and computation. Interestingly, even in the presence of stochastic observational noise, LACKI was able to yield  worst-case prediction error below the level of observational error. Furthermore, the prediction errors seemed to vanish in the limit of increasing sample size with a rate matching or outperforming the competing regression methods (ref. Fig. \ref{fig:bench01}). Furthermore, observe that the performance deterioration with increasing dimensionality of the input space (ref. Fig. \ref{fig:bench00}) was less than with the GPs. This is noteworthy since one might expect the GP-based predictors to benefit from stronger regularisation. Of course, we cannot claim that this superior performance over GPs will hold in general, but it is interesting to note that it can hold.

As stated above, the superior performance of the LACKI approach was contingent on setting the observational noise parameter to the correct value of $\hestthresh = 2 \obserrbnd$. In fact, Exp. 2 was designed to expose a shortcoming of our approach: namely the sensitivity to correctly setting $\hestthresh$ to zero in the absence of observational noise ($\obserrbnd =0$) when no information about the Lipschitz constant of the target is known. As we can see in Fig. \ref{fig:bench02}, falsely setting $\hestthresh =1$ resulted in poor predictive performance in Exp. 2. We explain this as follows: the target was confined to the interval $[-1.005,1.005]$. With gradients being very small at values close to the boundary of the interval, setting $\hestthresh=1$ denied an increase of $L(n)$. Thus, the resulting predictor of LACKI remained a constant and hence, gave rise to a relatively large and non-decreasing prediction error.  
Note, this problem would be prevented by either setting $L(n)$ to a valid Lipschitz constant of the target, or, by setting $\lambda =0$ reflecting the absence of observational errors. The positive effect of the latter is also testified by the plot of LACKI2 in Fig. \ref{fig:bench02} showing that our method for this setting was competitive with the best GP method on the noiseless regression task.

While our experiments might suggest that GPs may not always be as sensitive to the noise hyper-parameter settings we would like to emphasise the fact that the GP learners were carefully initialised to give good performance on the problems. And, suitable alterations of these choices could provoke substantially degradation of the GPs predictive performance, even in the limit of dense data.

In summary, the experimental results suggest the following observation: If the parameter $\hestthresh$ is correctly set to $2 \obserrbnd$, our LACKI method can offer a fast and reliable approach to regression under bounded stochastic noise that can outperform GP regression when assessed along performance metrics that reflect prediction accuracy and computational effort. We note that the latter could be further enhanced by applying nearest-neighbor approaches for Lipschitz interpolation in lieu to methods proposed in  \cite{Beliakov2006,calliess2014_thesis}.

\section{Lazily Adaptive Constant Kinky inference and model reference adaptive control}
\label{sec:MRAC_lacki}
So far, we have established some learning guarantees for LACKI as a method for supervised learning. In this section, we utilise our results and discuss the use of LACKI in the context of model-reference adaptive control. We introduce the control framework with a simple example of controlling the simulated roll dynamics of an aircraft under wing rock. In the second half of the remainder of the paper, we use our theoretical guarantees established in the first half in order to provide global asymptotic convergence guarantees of the closed-loop trajectory to the reference. 
Throughout the entire section we simplify our analysis by assuming the pseudo-metrics $\metric$ and $\metric_\outspace$ are in fact canonical norm-induced metrics. For instance, we assume that $\inspace = \Real^d$ is a finite-dimensional vector space and we have $\metric(x,x') =\norm{x-x'}$ for some norm $\norm{\cdot}$ equivalent to the maximum norm $\norm{\cdot}_\infty$.

\subsection{Online learning and tracking control in the presence of wing rock dynamics }
\label{sec:KIMRAC}
As pointed out in \cite{chowdharyacc2013}, modern fighter aircraft designs are susceptible to lightly damped oscillations in roll known as ``wing rock''. Commonly occurring during landing \cite{Saad2000}, removing wing rock from the dynamics is crucial for precision control of such aircraft.
Precision tracking control in the presence of wing rock is a nonlinear problem of practical importance and has served as a test bed for a number nonlinear adaptive control methods \cite{Chowdhary2013,Monahemi1996,chowdharyacc2013}.

For comparison, we replicate the experiments of the recent work of Chowdhary et. al. \cite{Chowdhary2013,ChowdharyCDC2013}.\footnote{We are grateful to the authors for kindly providing the code.}
Here the authors have compared their Gaussian process based approach, called \textit{GP-MRAC}, to the more established adaptive model-reference control approach based on RBF networks \cite{Sanner1992,Kim1998}, referred to as \textit{RBFN-MRAC}. Replacing the Gaussian process learner by our kinky inference learner, we readily obtain an analogous approach which we will refer to as \textit{LACKI-MRAC}. As an additional baseline, we also examine the performance of a simple P-controller.

While with the exact same parameters settings of the experiments in \cite{Chowdhary2013}, performance of our LACKI-MRAC method comes second to GP-MRAC, we also evaluate the performance of all controllers over a range of 555 random parameter settings and initial conditions. As we will see, across this range of problem instances and parameter settings, LACKI-MRAC markedly outperforms all other methods.

\subsubsection{Model reference adaptive control}\label{sec:mrac}
Before proceeding with the wing rock application we will commence with (i) outlining model reference adaptive control (MRAC) \cite{astroemadaptivectrlbook2013} as considered in \cite{Chowdhary2013} and (ii) describe the deployment of kinky inference to this framework. 
We will now rehearse the description of MRAC for second-order systems following \cite{Chowdhary2013}. 

Assume $m \in \nat$ to be the dimensionality of a configuration of the system in question and define $d = 2m$ to be the dimensionality of the pertaining state space $\statespace$.

Let $x = [x_1;x_2] \in \statespace$ denote the state of the plant to be controlled.
Given the control-affine system 
 
\begin{align}
\dot x_1 &= x_2 \\
\dot x_2 &= a(x) + b(x) \, u(x) \label{eq:secorddynctrlaff}
\end{align}

it is desired to find a control law $u(x)$ such that the closed-loop dynamics exhibit a desired reference behaviour:

\begin{align}
\dot \xi_1 &= \xi_2 \\
\dot \xi_2 &= f_{r}(\xi,r)
\end{align}
where $r$ is a reference command, $f_r$ some desired response and $t \mapsto \xi (t)$ is the reference trajectory.

If a priori $a$ and $b$ are believed to coincide with $\hat a_0, \hat b_0$ respectively, the inversion control 
$u = \hat b_0^{-1} (- \hat a_0 +u')$ is applied. This reduces the closed-loop dynamics to 
$\dot x_1 = x_2, \dot x_2 = u' + \tilde a(x,u) $
where $\tilde a(x,u)$ captures the modelling error of the dynamics: 
\begin{equation}
	\tilde a (x,u ) = a(x) - \hat a_0(x) + \bigl(b(x) - \hat b_0(x)\bigr) u.
\end{equation}
 Let $I_d \in \Real^{d \times d}$ denote the identity matrix.  If $b$ is perfectly known, then $b - \hat b_0^{-1} = 0$ and the model error can be written as $\tilde a (x)= a(x) - \hat a_0(x)$. In particular, $\tilde a$ has lost its dependence on the control input.

In this situation \cite{Chowdhary2013,ChowdharyCDC2013} propose to set 
the pseudo control as follows: $u'(x) :=  \nu_{r} + \nu_{pd} - \nu_{ad}$ where $\nu_{r} = f_{r}(\xi,r)$ is a feed-forward reference term,  $\nu_{ad}$ is a yet to be defined output of a learning module \emph{adaptive element} and $\nu_{pd} = [K_1 K_2] e$ is a feedback error term designed to decrease the \textit{tracking error} $e(t) = \xi(t) - x(t)$ by defining $K_1,K_2 \in \Real^{m \times m}$ as described in what is to follow.

Inserting these components, we see that the resulting \textit{error dynamics} are:

\begin{equation}\label{eq:errordynmrac}
	\dot e = \dot \xi - [x_2; \nu_r + \nu_{pd}+ \tilde a(x) ] = M e + B \bigl(\nu_{ad}(x) -  \tilde a(x)\bigr)
\end{equation}

where $M = \left(\begin{array}[h]{cc}
			O_m &  \, I_{m}\\
			-K_1 & -K_2 
					\end{array}\right)$ and $B = \left(\begin{array}[h]{c}
			O_m \\ I_m
					\end{array}\right)$.
If the feedback gain matrices $K_1,K_2$ parametrising $\nu_{pd}$ are chosen such that $M$ is stable then the error dynamics converge to zero as desired, provided the learning error $E_\lambda$ vanishes: $E_\lambda (x(t)) = \norm{\nu_{ad}(x(t)) -  a(x(t))} \stackrel{t \to \infty} {\longrightarrow} 0$. 

It is assumed that the adaptive element is the output of a learning algorithm that is tasked to learn $\tilde a$ online. This is done by continuously feeding it training examples of the form $\bigl(x(t_i), \tilde a(x(t_i)) + \varepsilon_i\bigr)$ where $\varepsilon_i$ is observational noise.  

Intuitively, assuming the learning algorithm is suitable to learn target $\tilde a$ (i.e. $\tilde a$ is close to some element in the hypothesis space \cite{mitchellbook:97} of the learner) and that the controller manages to keep the visited state space bounded, the learning error (as a function of time $t$) should vanish.

Substituting different learning algorithms yields different adaptive controllers. \textit{RBFN-MRAC} \cite{Kim1998} utilises radial basis function neural networks for this purpose whereas \textit{GP-MRAC} 
employs Gaussian process learning \cite{GPbook:2006} to learn $\tilde a$ \cite{Chowdhary2013,ChowdharyCDC2013}. 


In what is to follow, we utilise our LACKI method as the adaptive element. Following the nomenclature of the previous methods we name the resulting adaptive controller \textit{LACKI-MRAC}.

\subsubsection{The wing rock control problem}
The wing rock dynamics control problem considers an aircraft in flight. Denoting $x_1$ to be the roll attitude (angle of the aircraft wings) and $x_2$ the roll rate (measured in angles per second), the controller can set the aileron control input $u$ to influence the state $x := [x_1;x_2]$.

Based on \cite{Monahemi1996}, Chowdhary et. al. \cite{Chowdhary2013,ChowdharyCDC2013} consider the following model of the wing rock dynamics: 

\begin{align}
\dot x_1 &= x_2 \\
\dot x_2 &= a(x) + b \, u 
\end{align}
where $b =3$ is a known constant and 
$a(x) = W_0^* + W_1^* x_1 + W_2^* x_2 + W_3^* \abs{x_1} x_2 + W_4^* \abs{x_2} x_2 + W^*_5 x_2^3$ is an priori unknown nonlinear drift.

Note, the drift is non-smooth but it would be easy to derive a Lipschitz constant on any bounded subset of state space if the parameters $W := (W_0^*,\ldots, W_5^*)$ were known.

To control the system we employ LACKI as the adaptive element $\nu_{ad}$.
In the absence of the knowledge of a Lipschitz constant, we start with a guess of $\underline L=1$ (which will turn out to be too low) and update it following the procedure described in Sec. \ref{sec:lacki}.

In a first instance, we replicated the experiments conducted in \cite{Chowdhary2013,chowdharyacc2013} with the exact same parameter settings. That is, we chose $W_0^* = 0.8, W_1^* = 0.2314, W_2^* = 0.6918, W_3^* = -0.6245, W_4^* = 0.0095, W_5^* = 0.0214$.

The simulation initialised with start state $x = (3,6)^\top$ and simulated forward with a first-order Euler approximation with time increment $\tinc = 0.005 [s]$ over a time interval $\indsett = [t_0,t_f]$ with $t_0 = 0[s]$ and $t_f = 50[s]$. Training examples and control signal were continuously updated every $\Delta_u= \Delta_o = \tinc [s]$. The RBF and GP learning algorithms were initialised with fixed length scales of 0.3 units. The GP was given a training example budget of a maximum of 100 training data points to condition the posterior model on. Our LACKI learner was initialised with $\underline L =\hexp= 1$ and updated online following our lazy update method described above.

The test runs also exemplify the working of the lazy update rule.
The initial guess $\underline L=1$ was too low. However, our lazy update rule successfully picked up on this and had ended up increasing constant to $L=2.6014$ by the end of the online learning process.

\begin{figure*}
        \centering
				  \subfigure[Tracking error (RBF-MRAC).]{
    \includegraphics[ width=.3\textwidth, clip, trim = 2.5cm 8cm 3cm 9cm]
								{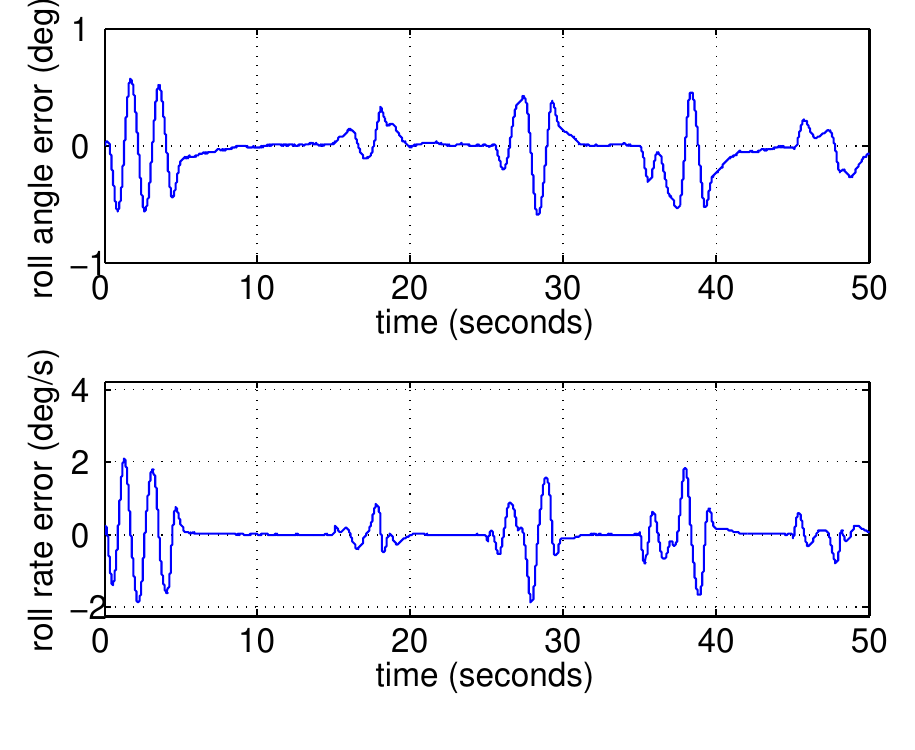}
  } 
					  \subfigure[Tracking error (GP-MRAC).]{
    \includegraphics[width=.3\textwidth, clip, trim = 2.5cm 8cm 3cm 9cm]
								{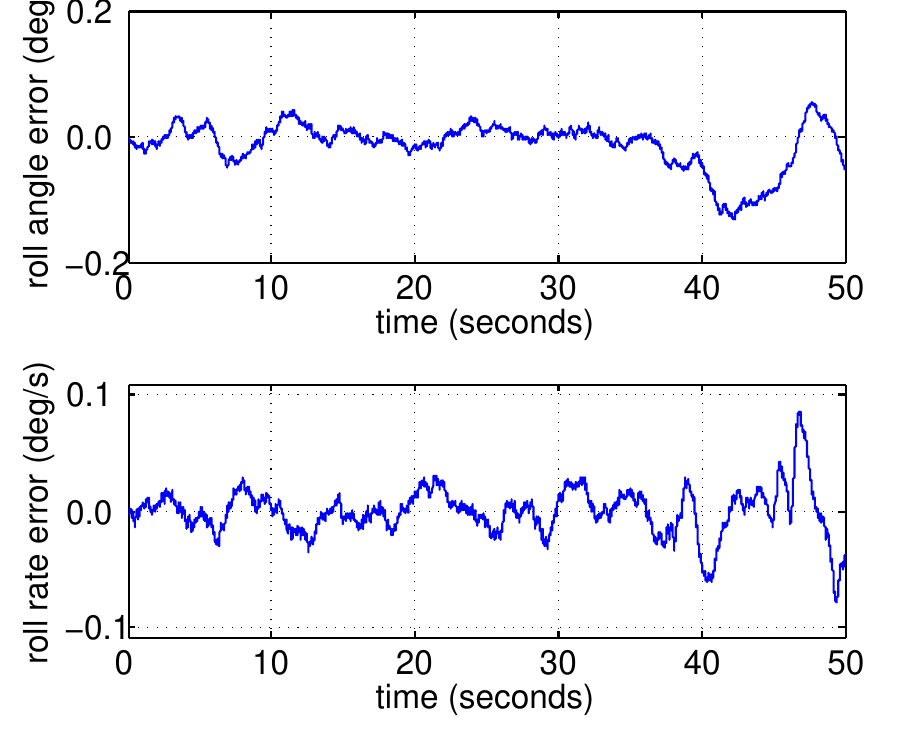}
  } 
					  \subfigure[Tracking error (KI-MRAC).]{
    \includegraphics[width=.3\textwidth, clip, trim = 2.5cm 8cm 3cm 9cm]
								{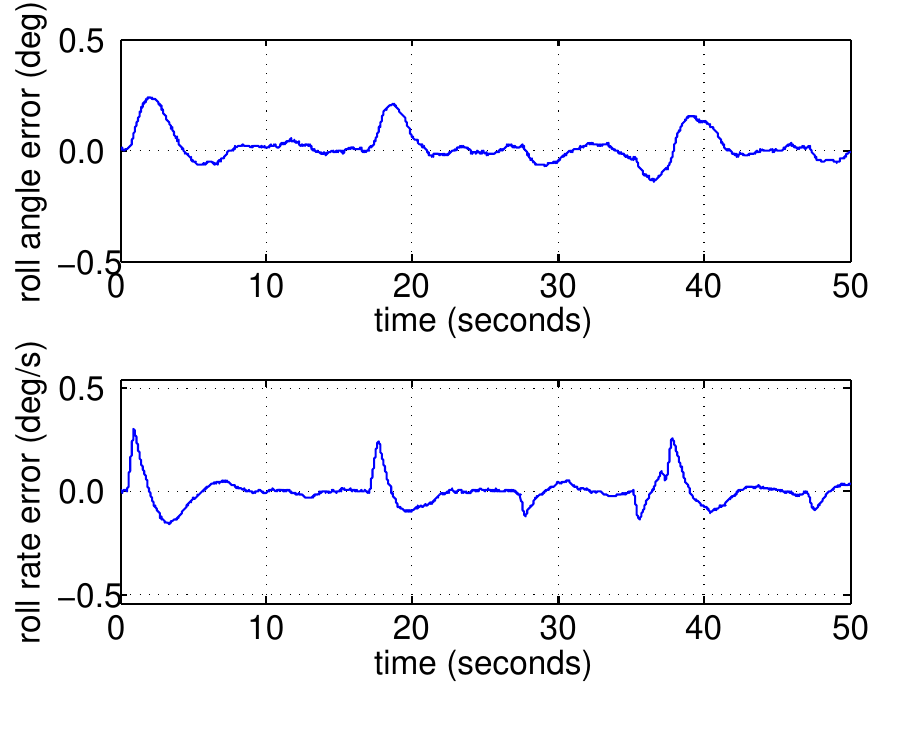}
  }
   \caption{Tracking error comparison of first example.}
	 \label{fig:wrtrackerrorsex1}
\end{figure*}	   

The results are plotted in Fig. \ref{fig:wrtrackerrorsex1}. 
We can see that in terms of tracking error of the reference our LACKI-MRAC outperformed RBF-MRAC and was a close runner-up to GP-MRAC which had the lowest tracking errors.  

\begin{figure*}
        \centering
				  \subfigure[Prediction v.s. ground truth  (RBF-MRAC).]{
    \includegraphics[width = .3\textwidth, clip, trim = 3cm 9cm 3cm 10cm]
								{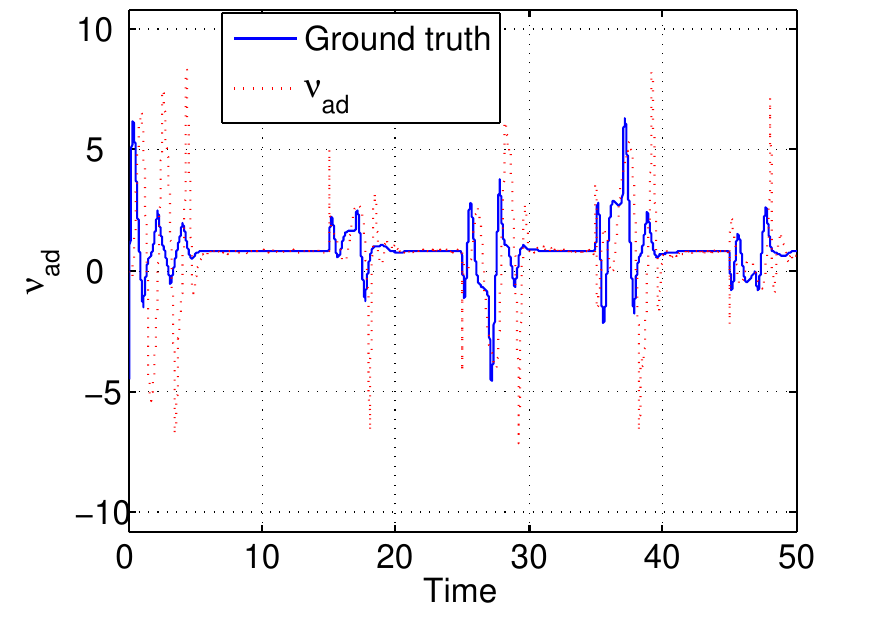}
  } 
					  \subfigure[Prediction v.s. ground truth  (GP-MRAC).]{
    \includegraphics[width=.3\textwidth, clip, trim = 3cm 9cm 3cm 10cm]
								{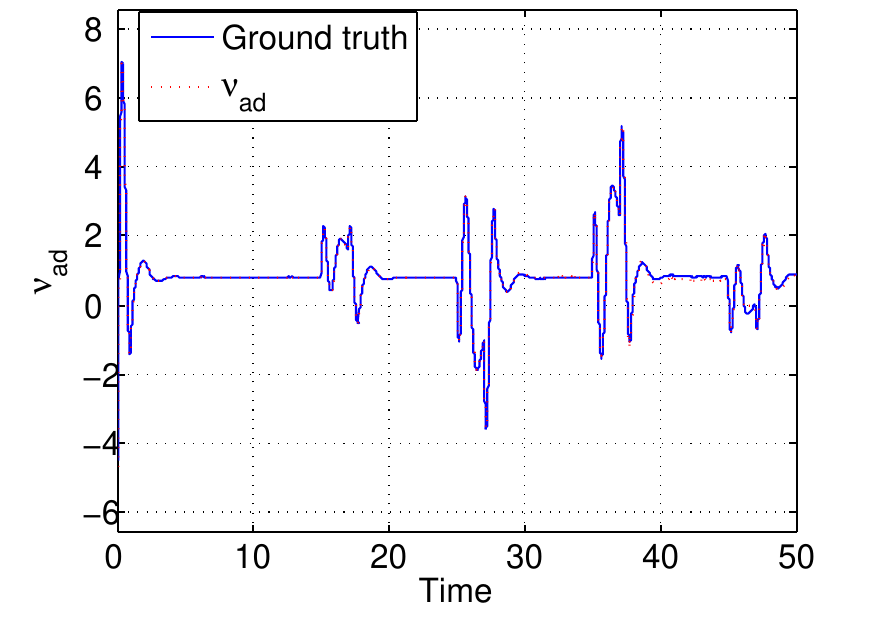}
  } 
					  \subfigure[Prediction v.s. ground truth  (LACKI-MRAC).]{
    \includegraphics[width=.3\textwidth, clip, trim = 3cm 9cm 3cm 10cm]
								{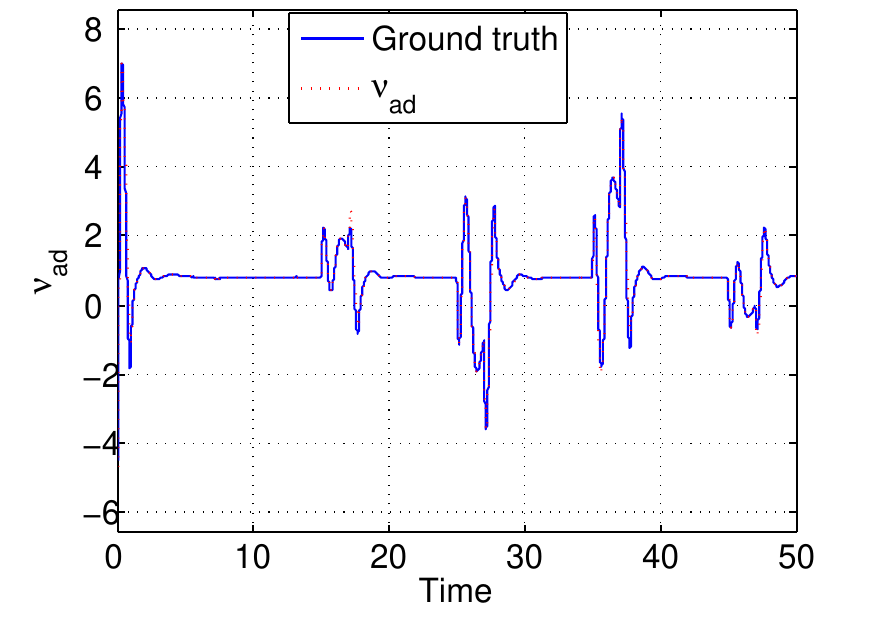}
  } 
   \caption{Prediction vs ground truth comparisons for the first example. Both nonparametric methods accurately predict the true drift and clearly outperform the RBFN learner.}
	\label{fig:wrprederrorsex1}
\end{figure*}

To obtain an impression of the learning performance of the three learning algorithms we also recorded the prediction error histories for this example problem. The results are depicted in Fig. \ref{fig:wrprederrorsex1}. We can see that our kinky inference method and the GP method both succeeded in predicting the drift quite accurately while the RBFN method was somewhat lagging behind.
This is consistent with the observations made in \cite{Chowdhary2013,ChowdharyCDC2013}. The authors explain the relatively poor performance of the radial basis function network method by the fact that the reference trajectory on occasion led outside the region of state space where the centres of the basis function were placed in advance. By contrast, due to the non-parametric nature of the GP, GP-MRAC does not suffer from such a priori limitations. In fact, it can be seen as an RBF method that flexibly places basis functions around all observed data points \cite{GPbook:2006}. We would add that, as a non-parametric method, LACKI-MRAC shares this kind of flexibility, which might explain the fairly similar performance. 

However, being an online method, the authors of GP-MRAC explicitly avoided hyperparameter training via optimising the marginal log-likelihood. The latter is commonly done in GP learning \cite{GPbook:2006} to avoid the impact of an unadjusted prior but is often a computational bottle neck. Therefore, avoiding such hyperparameter optimisation greatly enhances learning and prediction speed in an online setting. However, we would expect the performance of the prediction to be dependent upon the hyperparameter settings. As we have noted above, the Lipschitz constant depends on the part of state space visited at runtime. Similarly, we might expect length scale changes depending on the part of state space the trajectory is in. Unfortunately, \cite{Chowdhary2013,ChowdharyCDC2013,chowdharyacc2013} provide no discussion of the length scale parameter setting and also called the choice of the maximum training corpus size ``arbitrary''. 

Since the point of learning-based and adaptive control is to be able to adapt to various settings, we test the controllers across a range of randomised problem settings, initial conditions and parameter settings.

We created 555 randomised test runs of the wingrock tracking problems and tested each algorithm on each one of them. The initial state $x(t_0)$ was drawn uniformly at random from $[0,7] \times [0,7]$, the initial kernel length scales were drawn uniformly at random from $[0.05,2]$, and used both for RBF-MRAC and GP-MRAC. The initial H\"older constant $\underline L$ for LACKI-MRAC was initialised at random from the same interval but was allowed to be adapted as part of the online learning process. Furthermore, we chose $\hestthresh =0$. The parameter weights $W$ of the system dynamics specified above were multiplied by a constant drawn uniformly at random from the interval $[0,2]$. To allow for better predictive performance of GP-MRAC we doubled the maximal budget to 200 training examples. 
The feedback gains were chosen to be $K_1=K_2=1$. 

In addition to the three adaptive controllers we also tested the performance of a simple $PD$ controller with just these feedback gains (i.e. we executed x-MRAC with adaptive element $\nu_{ad}=0$). This served as a baseline comparison to highlight the benefits of the adaptive element over simple feedback control.

\begin{figure}
        \centering
				  \subfigure[Results over 555 randomised examples.]{
    \includegraphics[width = .5\textwidth, clip, trim = 1.5cm 8cm 1cm 7cm]
								{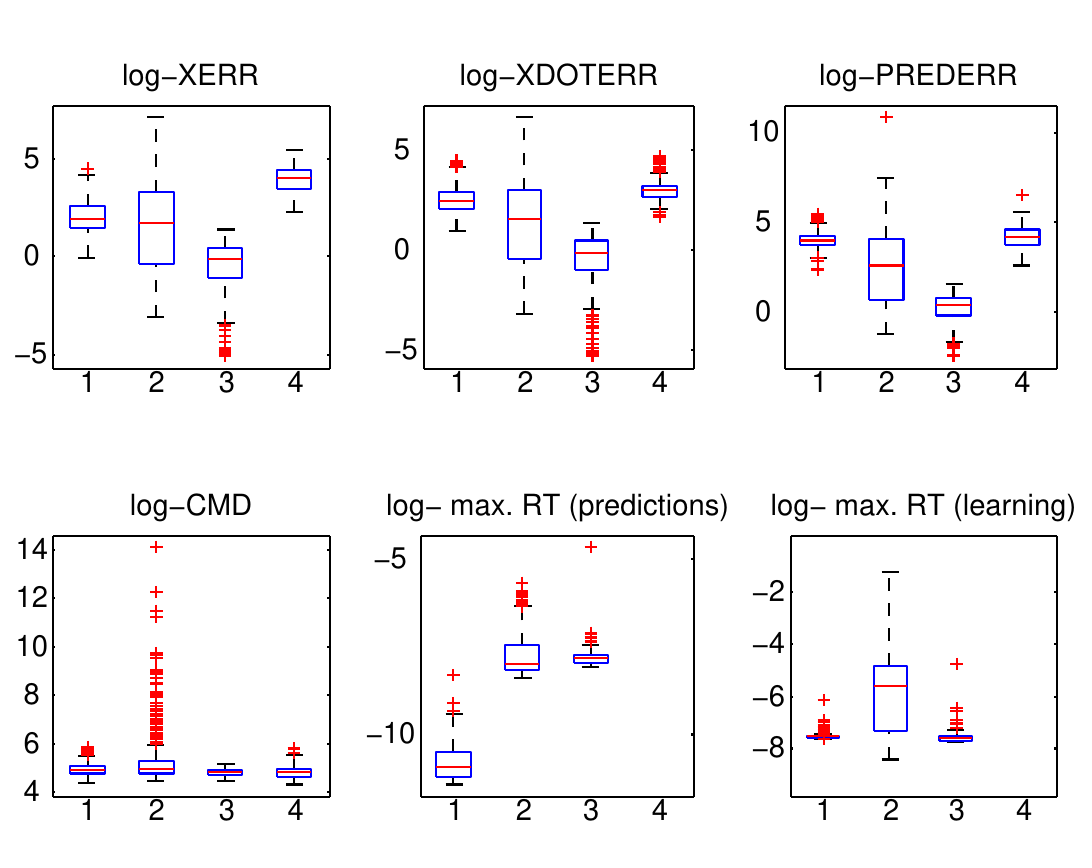}
  } 
   \caption{Performance of the different online controllers over a range of 555 trials with randomised parameter settings and initial conditions. 1: RBF-MRAC, 2: GP-MRAC, 3:LACKI-MRAC, 4: P-Controller. LACKI-MRAC outperforms all other methods with respect to all performance measures, except for prediction runtime (where the parametric learner RBF-MRAC performs best).} 
       \label{fig:wingrockresultsbp}
\end{figure}

The performance of all controllers across these randomised trials is depicted in Fig. \ref{fig:wingrockresultsbp}. Each data point of each boxplot represent a performance measurement for one particular trial.

	For each method, the figures show the boxplots of the following recorded quantities: 
	\begin{itemize}
		\item \textit{log-XERR}: cummulative angular position error (log-deg), i.e. $\log(\int_{t_0}^{t_f} \norm{\xi_1(t) - x_1 (t)} \dt )$.
		\item \textit{log-XDOTERR}:  cummulative roll rate error (log-deg/sec.), i.e. $\log(\int_{t_0}^{t_f} \norm{\xi_2(t) - x_2 (t)} \dt )$.
		\item \textit{log-PREDERR}: log-prediction error, i.e. 
		
		$\log(\int_{t_0}^{t_f} \norm{\nu_{ad}(x(t)) - \tilde a(x(t))} \dt )$.
		\item \textit{log-CMD}: cummulative control magnitude (log-scale), i.e. $\log(\int_{t_0}^{t_f} \norm{u(t)} \dt )$.
		\item \textit{log-max. RT (predictions)}: the log of the maximal runtime (within time span $[t_0,t_f]$) each method took to generate a prediction $\nu_{ad}$ within the time span.
		\item \textit{log-max. RT (learning)}: the log of the maximal runtime (within time span $[t_0,t_f]$) it took each method to incorporate a new training example of the drift $\tilde a$.
	\end{itemize}
	
	As can be seen from Fig. \ref{fig:wingrockresultsbp}, all three adaptive methods outperformed the simple $\hexp$ controller in terms of tracking error. 
	
	In terms of prediction runtime, the RBF-MRAC outperformed both GP-MRAC and LACKI-MRAC. This is hardly surprising. After all, RBF-MRAC is a parametric method with constant prediction time. By contrast, both non-parametric methods will have prediction times growing with the number of training examples.
That is, it would be the case if GP-MRAC were given an infinite training size budget. Indeed one might argue whether GP-MRAC, if operated with a finite budget, actually is a parametric approximation where the parameter consists of the hyperparameters along with the fixed-size training data matrix. When comparing the (maximum) prediction and learning runtimes one should also bear in mind that GP-MRAC predicted with up to 200 examples in the training data set. By contrast, LACKI-MRAC undiscerningly had incorporated all 10001 training points by the end of each trial.

Across the remaining metrics, LACKI-MRAC markedly outperformed all other methods.

Note, we have also attempted to test all methods across a greater range of problem settings, including larger initial states, more varied hyper-parameter settings, lower feedback gains and more varied choices of dynamics coefficients $W$. However, this resulted in GP-MRAC to often run into conditioning problems. This is a common issue in GP learning due to the necessity of matrix inversion or Cholesky decompositions of the covariance matrix (it seems to be common practice to address this by hand-tuning the observational noise parameters). Similar behaviour ensued when setting the training size budget to large values. All these changes often resulted in long learning runtimes, spiky control outputs and thus, poor overall performance. Similarly, code execution of our RBF-MRAC implementation was frequently interrupted with error messages when the state was initialised to positions outside of the rectangle $[0,7] \times [0,7]$.

We have not investigated the root cause of these issues in greater detail yet. However, it might be worth exploring whether the great robustness of LACKI might be an additional selling point that sets it apart from other recent adaptive control methods. Such robustness is of course important in control settings such as flight control where failure or erratic behaviour of the adaptive element may result in critical incidents. 

An example where GP-MRAC failed to track the reference occurred when repeating our first experiment  with the following modifications: The initial state was chosen to be $x(t_0) = (-90,40)^\top$ corresponding to a rapidly rotating aircraft. Furthermore, the wing rock coefficients $W$ were multiplied by a factor of $5$, amplifying the non-linearities of the drift field. 

\begin{figure*}
        \centering
				  \subfigure[Position (GP-MRAC).]{
    \includegraphics[width = 5.5cm, clip, trim = 4cm 9cm 4cm 9cm]
								{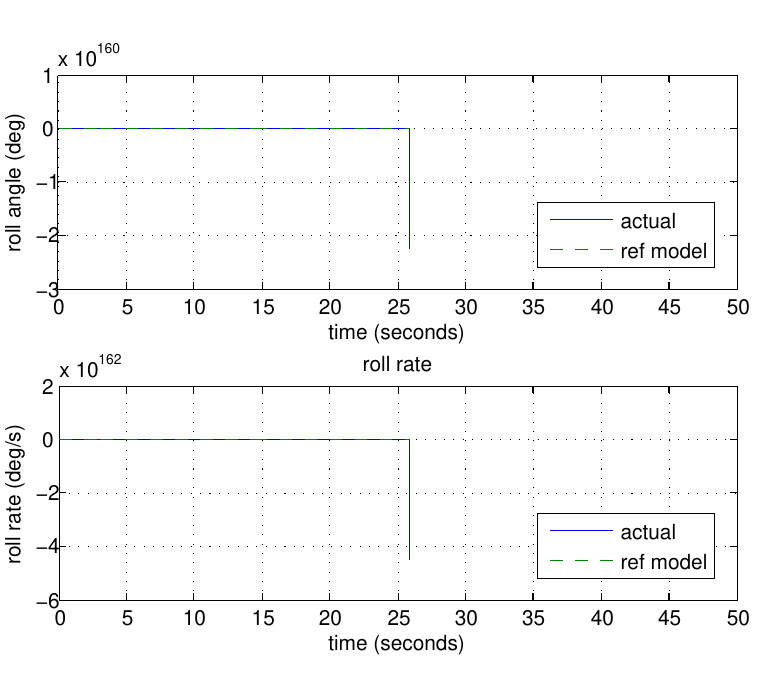}
   \label{fig:gpfailGPWR1}
  } 
					  \subfigure[Tracking error (GP-MRAC).]{
    \includegraphics[width = 5.5cm, clip, trim = 4cm 9cm 4cm 9cm]
								{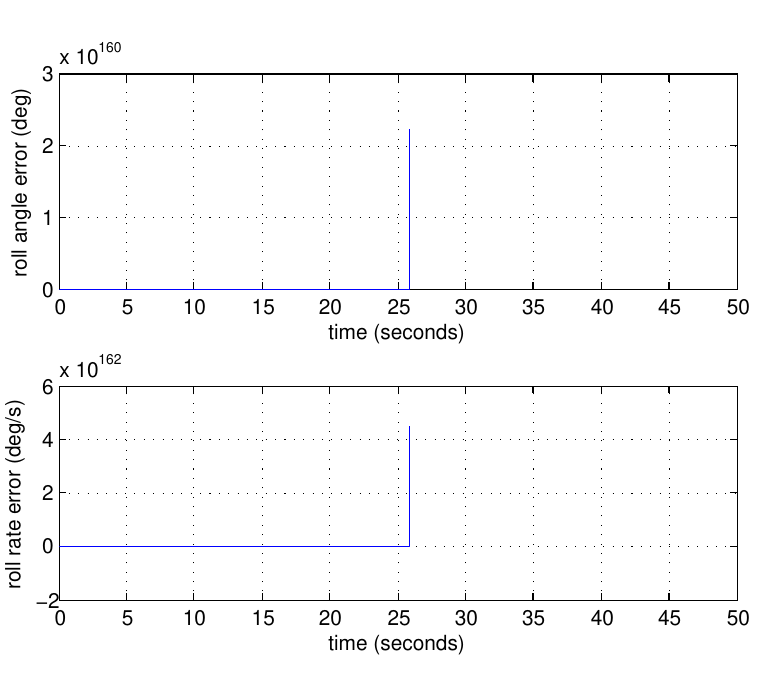}
   \label{fig:gpfailGPWR2}
  } 
							  \subfigure[Log - prediction error (GP-MRAC).]{
    \includegraphics[scale =.35]
								{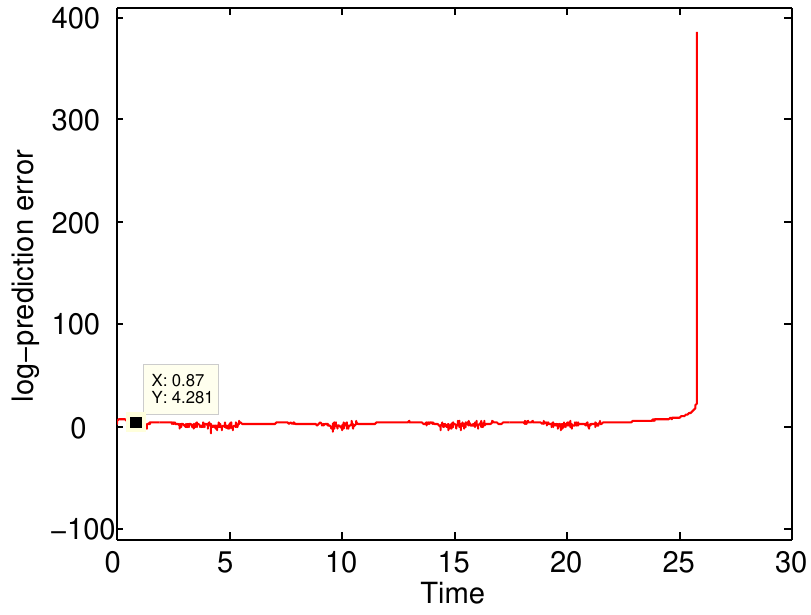}
    \label{fig:gpfailGPWR3}
  } 
					  \subfigure[Position (LACKI-MRAC).]{
    \includegraphics[width = 5.5cm, clip, trim = 3.5cm 9cm 4cm 10cm]
								{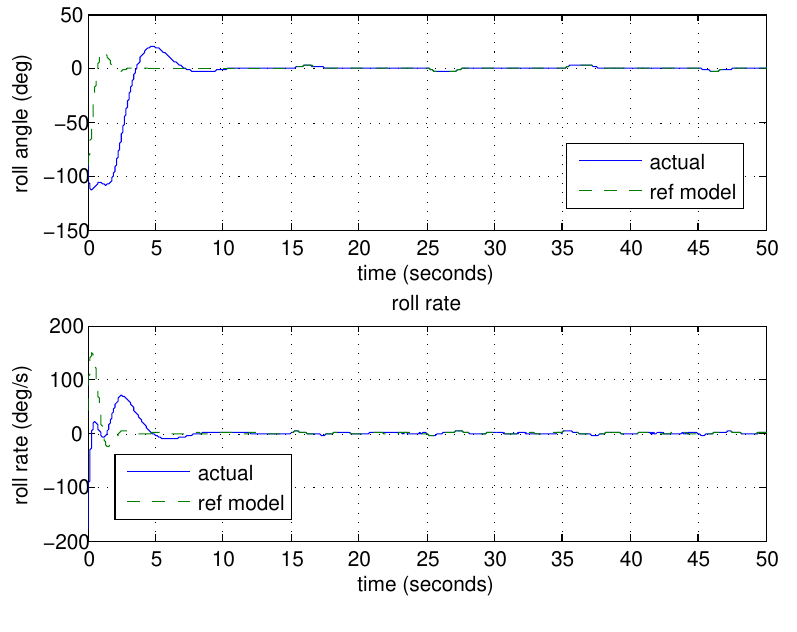}
    \label{fig:gpfailKIWR1}
  } 
					  \subfigure[Tracking error (LACKI-MRAC).]{
    \includegraphics[width = 5.5cm, clip, trim = 3.5cm 9cm 4cm 10cm]
								{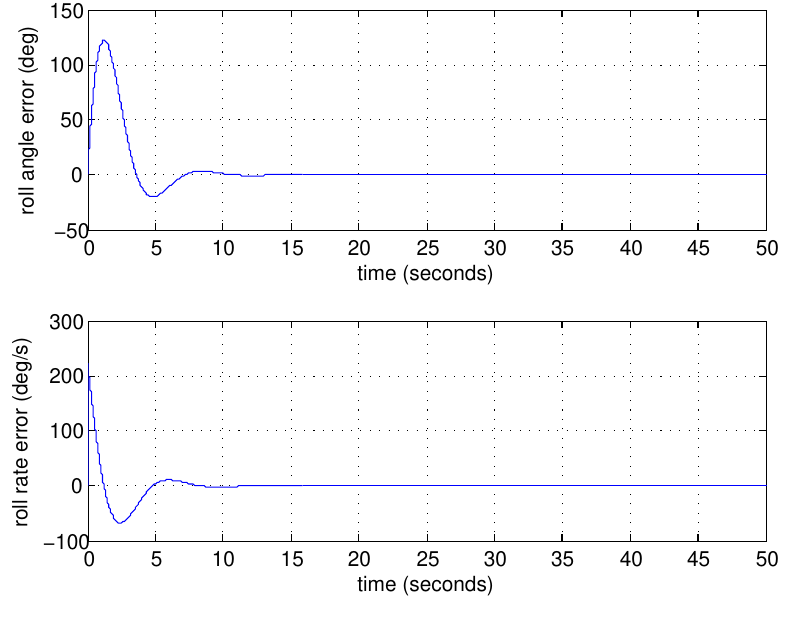}
    \label{fig:gpfailKIWR2}
  } 
					  \subfigure[Log - prediction error (LACKI-MRAC).]{
    \includegraphics[scale =.35,clip, trim = 0cm 0cm 0cm .1cm]
								{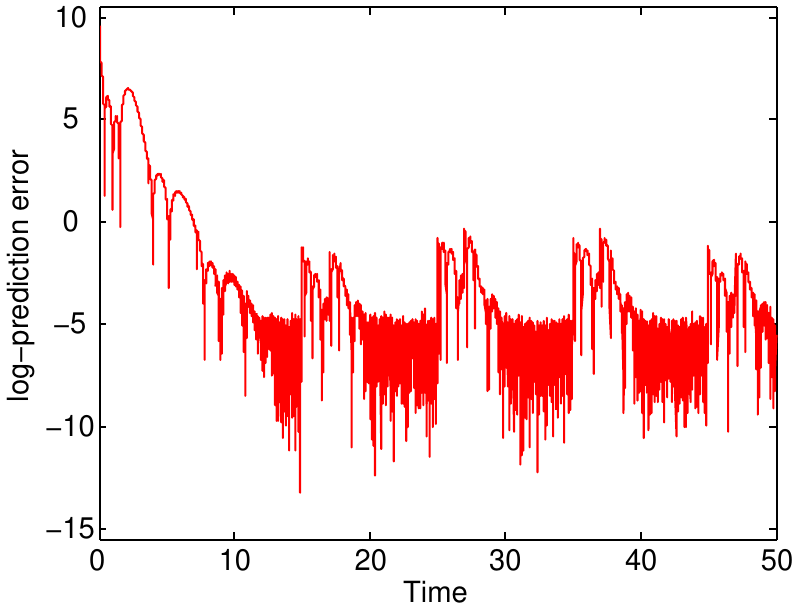}
    \label{fig:gpfailKIWR3}
  } 
	%
   \caption{Example where GP-MRAC fails. By contrast, LACKI-MRAC manages to adapt and direct the system back to the desired trajectory.}
	\label{fig:gpfail}
\end{figure*}

When initialised with a length scale parameter of 0.3, the GP ran into conditioning problems and caused the output of the adaptive element in GP-MRAC to produce spikes of very large magnitude and thus, further destabilised the system. We tried the problem with various kernel length scale settings ranging from $0.3$ to $20$. Increasing the length scale parameter to length scale of at least 1 seemed to fix the conditioning problem. Nonetheless, GP-MRAC still did not manage to learn and stabilise the system in any of these settings. A record of GP-MRAC's performance in this example (for length scale of 1) is depicted in Fig.  \ref{fig:gpfailGPWR1} -  \ref{fig:gpfailGPWR3}. As the plots show, GP-MRAC starts with relatively high tracking and prediction error from which it could not recover. At about 26 seconds into the simulation the state rapidly diverged.

\begin{figure*}
        \centering
				  \subfigure[State path (LACKI-MRAC).]{
    \includegraphics[width = .4\textwidth, clip, trim = 3cm 9cm 3cm 9cm]
								{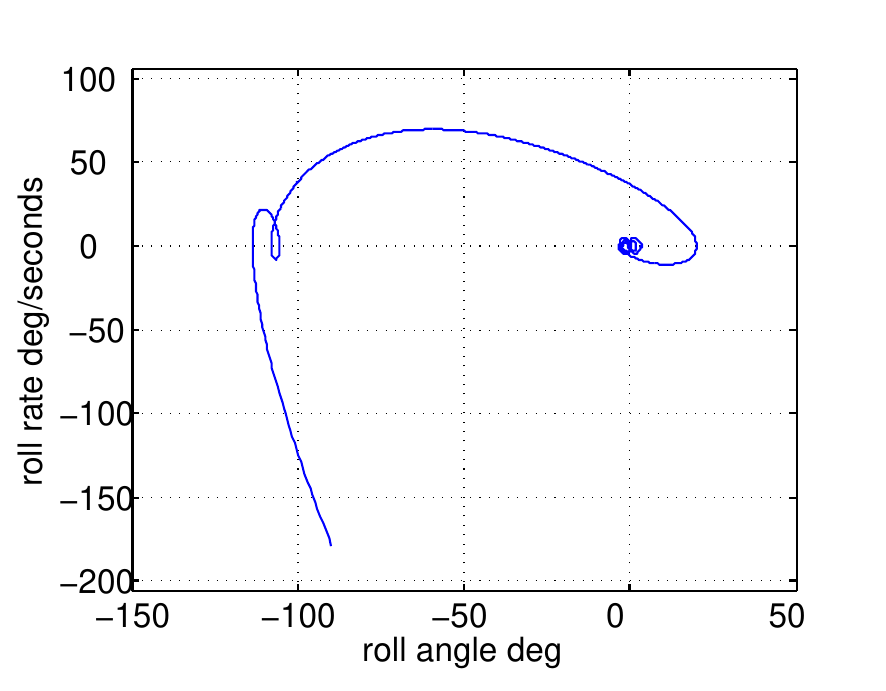}
  } 	
	  \subfigure[Learned drift model (LACKI-MRAC).]{
    \includegraphics[width = .4\textwidth, clip, trim = 4.3cm 9cm 4cm 9cm]
								{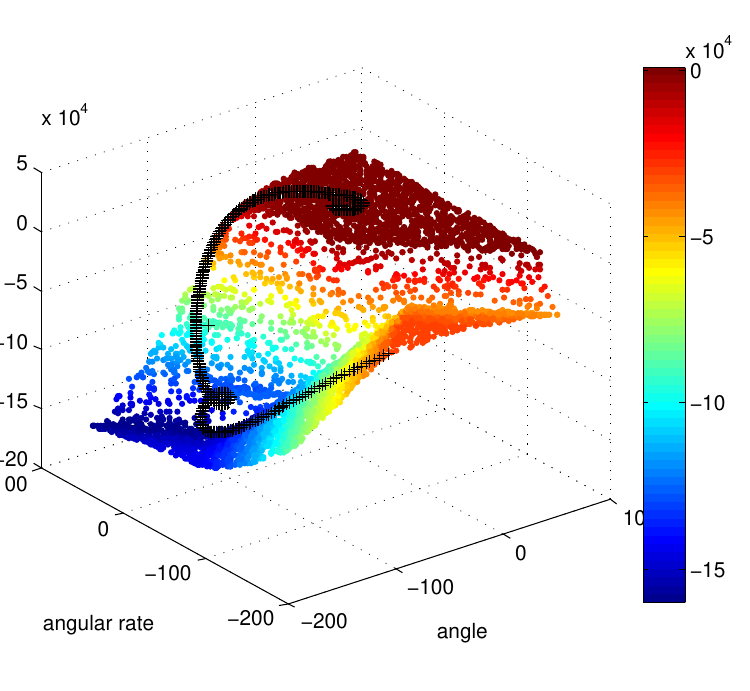}
  } 	
   \caption{Depicted are the state path and the drift model learned online by LACKI-MRAC.}
	\label{fig:gpfail2}
\end{figure*}

For comparison, we also tried LACKI-MRAC on the same problem, starting with initial $L=1$ as before. Starting out with a relatively large tracking and prediction error, LACKI-MRAC nonetheless managed to recover and successfully track the system (see  Fig.  \ref{fig:gpfailKIWR1} -  \ref{fig:gpfailKIWR3}). The state path and learned drift model obtained by LACKI-MRAC are depicted in Fig. \ref{fig:gpfail2}.

\subsection{Convergence guarantees in discrete-time systems}
\label{sec:KIMRACstabbounds}

In the previous section, we gave an illustration LACKI-MRAC -- a combination of a feedback-linearising controller with our KI learning method. The results are encouraging -- our adaptive control law managed to learn a dynamic system online and to track a reference in the presence of wing-rock dynamics where other methods failed. 
To simulate the dynamics we relied on a first-order Euler approximation resulting in a discrete-time dynamical system.

In this section, we study the convergence of LACKI-MRAC in such discrete-time systems.
In these we will consider both the offline and the online learning setting. In the former, the LACKI-learner receives a sample set once and builds a controller that remains unaltered during execution time.
For this case, we will provide robustness bounds on the control success (as quantified by a bound on the norm of the error dynamics) as a function of the remaining uncertainty of the trained LACKI model. 

In the online learning setting, which we considered in the wing-rock control simulations, the LACKI-learner gets updated with the most recent observation after each time step. Provided the initial uncertainty is bounded on the given state space, we will be able to guarantee that LACKI-MRAC leads to a closed-loop system that achieves the tracking objective with increasing time and learning experience. 

\subsubsection{Tracking error bounds for the offline learning setting} 
In the offline learning setting, the predictor sequence $\seq{\predfn}{n \in \nat}$ of adaptive elements is based on only one fixed data set $\data_0$ at time $0$ which is not updated subsequently.
That is $\data_n = \data_0, \forall n \in \nat$.

Given some data set $\data_n =\data_0$ at time $n$ and the resulting predictor $\predfn(\cdot)$, the model error is given by $F_n(\cdot) := f(\cdot) - \predfn(\cdot)$. Since we assume constant data, the error does not change either. That is, we have $F_n(x) = \ldots = F_0(x) = F(x), \forall n,\forall x$.
 
In our analysis, we consider discrete-time dynamical systems. For example the dynamics might be first-order Euler approximations of the control-affine dynamics of Sec. \ref{sec:KIMRAC}.
Consequently, the error dynamics as per Eq. \ref{eq:errordynmrac} translate to the recurrence relation:
\begin{equation}\label{eq:errordynmrac_nIMRACdiscrete}
	e_{n+1}  = M e_n + \tinc F(x_n)
\end{equation}
where $\tinc \in \Real_+$ is a positive time increment and  $n$ is the time step index. Furthermore,
\begin{equation}
	F(x_n) = f(x_n) - \predfn(x_n) = B \bigl(\nu_{ad}(x_n) -  \tilde a(x_n)\bigr) 
\end{equation}
is an uncertain increment due to the model error of the learner (cf. Eq. \ref{eq:errordynmrac}), 						$B = \left(\begin{array}[h]{c}
				O_m \\ \tinc I_m
						\end{array}\right)$ and 

\begin{equation}
	M = \left(\begin{array}[h]{cc}
				I_m &  \, \tinc I_{m}\\
				-\tinc K_1 & I_m- \tinc K_2 
						\end{array}\right) 
\end{equation}
					is the (error state) transition matrix. 	
Here, $m = \frac d 2$ is half the dimensionality of the state space, $I_m$ denotes the $m \times m$ identity matrix and $K_1,K_2$ are gain matrices that can be freely chosen by the designer of the linear pseudo controller.		
By induction, it is easy to show that the recurrence can be converted into the closed-form expression:
\begin{align*}
	\vc e_n &= M \vc e_{n-1} + \tinc \vc F(x_{n+1})  \\ 
	&= \dots= M^k \, \vc e_0 + \tinc \sum_{i=0}^{n-1} M^{n-1-i} \, \vc F(x_i)	.	
\end{align*}	
For vector norm $\norm{\cdot}$, let $\matnorm{\cdot}$ denote a matrix norm such that $\norm{Mv} \leq \matnorm{M} \norm{v}$ for all suitable matrices $M$ and vectors $v$. For instance, for the Euclidean norm $\norm{\cdot} =\norm{\cdot}_2$, we can choose the spectral norm $\matnorm{\cdot}= \specnorm{\cdot}$ as a matrix norm. Or, for  $\norm{\cdot} =\norm{\cdot}_\infty$, if the vector space is $d$-dimensional, we can choose the matrix norm $\matnorm{\cdot}= \sqrt{d} \specnorm{\cdot}$. We desire to bound the  norm of the error. To this end, we leverage that the norms are sub-additive and sub-multiplicative to deduce: 
\begin{align}
	\norm{\vc e_n} &\leq   \matnorm{M^{n}} \, \norm{\vc e_0} + \tinc \sum_{i=0}^{n-1}  \matnorm{M^{n-1-i}} \, \norm{\vc F(x_i)}		\\
	&\leq \matnorm{M^{n}} \, \norm{\vc e_0} + \tinc \maxerrn_n	 \sum_{i=0}^{n-1}  \matnorm{M^{n-1-i}} =: \varrho[n]
	\label{ineq:errnormbasic}
\end{align}
where  $\maxerrn_n$ is chosen such that we can guarantee that $\maxerrn_n	\geq \max_{i=1,...,n-1} \norm{F(x_i)}$.
For instance, we could choose $\maxerrn_n := \sup_{s \in \iaspace_n} \norm{F(s)} $
	where $\iaspace_n = \bigcup_{t< k} \{x \in \iaspace | \norm{x - \xi[t]} \leq \varrho[t] \}$ is the union of the possible states around the reference trajectory $\xi[\cdot]$ at previous time steps.
\begin{remark}	[Bounded error innovations] 
We assume there exists a maximum model error norm $\maxerrn$, i.e. $\forall i: \norm{F_i(x_i)} \leq \maxerrn$ for some bound $\maxerrn \in \Real$.
This is a realistic assumption in any physically plausible systems where the drift forces $a(x)$ are inevitably bounded. Since also, for any finite data, our LACKI learner has a bounded prediction function (cf. Lem. \ref{lem:LACKIpredbounded}), the discrepancy given by $F$ is bounded as well. 
\label{rem:bnd_err_innovations}
\end{remark}

Given the boundedness of $F$ we have: \vspace{-1em}
\begin{equation}
	\norm{\vc e_n} \leq   \matnorm{M^{n}} \, \norm{\vc e_0} + \tinc \maxerrn \sum_{i=0}^{n-1}  \matnorm{M^{n-1-i}}.
\end{equation}

The right-hand side is convergent provided the gains $K_1, K_2$ are chosen such that $M$ is stable, i.e. $\specrad(M) <1$, and provided $\maxerrn_n$ is bounded (see e.g. \cite{calliess2014_thesis}). 

Whether or not $M$ is stable, in low dimensions, the sums can be computed offline and in advance. This is of great benefit if 
 the controller that is building on the error bounds is utilising optimisation-based control with a finite time-horizon.   

To obtain a (conservative) closed-form bound on the error norms \cite{calliess2014_thesis}
contains a derivation and discussion of the following result:

\begin{thm}\label{thm:normboundsboundeddisturbmainbody}
Let $(F_n)_{n \in \nat_0}$ be a norm-bounded sequence of vectors with $\maxerrn_n :=\max_{i \in \{0,\ldots,n-1\}} \norm{F_i(x_i)} \leq \maxerrn \in \Real$. 
For error sequence $(e_n)_{n \in \nat_0}$ defined by the linear recurrence 
	$e_{n+1}  = M e_n + \tinc F_n(x_n) \,\,\,\, (n \in \nat_0)$, we can define the following bounds:
	
	\begin{enumerate}
	\item $\norm{\vc e_n} \leq \matnorm{M^{n}} \, \norm{\vc e_0} + \tinc \maxerrn_{n} \sum_{i=0}^{n-1}  \matnorm{M^{i}}$. If $\specrad(M) <1$ and $\exists \maxerrn$ $:\maxerrn\geq \maxerrn_{n-1} \geq 0, \forall n$, the right-hand side converges as $n \to \infty$.
	\item Let $k_0 \in \nat, k_0 >1$ such that $\matnorm{M^n} <1, \forall n \geq k_0$, let $\varphi := \matnorm{M^{k_0}} <1$ and let $\delta_n := \floor{n/k_0}$. If $r:=\specrad(M) < 1$, for $n > k_0$,  we also have:
	\begin{align*} \norm{\vc e_n} &\leq c \, \varphi^{\delta_n} \, \norm{\vc e_0} + \tinc \maxerrn_{n} \Bigl(\sum_{j=0}^{n_0-1} \matnorm{M^j} + c \, k_0 \frac{\varphi-\varphi^{ \floor{\frac{n}{k_0}}+1 }}{1-\varphi} \Bigr)\\
	& \stackrel{n \to \infty}{\longrightarrow} C \leq \tinc \maxerrn \sum_{j=0}^{n_0-1} \matnorm{M^j} 
	+ \frac{\tinc \maxerrn\,  c\, k_0 \,\varphi}{1-\varphi} 
	\end{align*} 
	for some constants $C,c \in \Real$.
	Here, possible choices for $c \in \Real$ are: 
	
(i) $c= \max\{\matnorm{M^i} | i=1,\ldots,k_0-1 \}$ 

or (ii) $c = \frac{1}{(d-1)!} \Bigl(\frac{1 - d}{\log r}\Bigr)^{d-1}	\matnorm{M}^{d-1}\, \, \, r^{\frac{1 - d}{\log r}-d+1}  $.
Since $\matnorm{M} \neq 1$, one can also choose (iii) $c:= \matnorm{M}^{n_0}$. 

\item If $\matnorm{M} \neq 1$, we have: \newline
$\norm{\vc e_n} \leq \matnorm{M}^{n} \, \norm{\vc e_0} + \tinc \maxerrn_{n}   \frac{1-\matnorm{M}^n}{1-\matnorm{M}}. $
	\end{enumerate}
%
\end{thm}

\subsubsection{Convergence guarantees in the online learning setting}
In this subsection, we lift the assumption that the available sample is static. Instead we assume that at time step $n+1$ we get to see an additional sample of the uncertain drift at the state visited in the previous time step $n$. 
That is, the predictor $\predf_{n+1}(\cdot)$ is based on $\data_{n+1} = \data_n \cup \{ \bigl(\state_n, \tilde f(\state_n) \bigr)\}, \forall n $.
Therefore, we also need to index the innovations vector field by time. That is, $F_n := f - \predfn$ denotes the prediction error function (or innovation) due to the data available at time $n \in \nat$. 
As pointed out in Rem. \ref{rem:bnd_err_innovations}, the error innovations $F_n$ are assumed to be bounded for all finite sample sets $\data_n$.

Above, we have seen that any continuous function can be approximated by some H\"older continuous LACKI predictor up to an arbitrarily small error. 
For convenience, we will establish the following definition:

\begin{defn}
We say that a continuous function $f$ is $L^*-p$-H\"older \emph{up to} error $\bar E_h \in \Real$ on domain $\inspace$ iff there is an $L^*-p$-H\"older function $\phi \in \hoelset {L^*} { } p$ and a function $\psi$ such that $\forall x: f(x) = \phi(x)+\psi(x), \, \sup_{x \in \inspace} \metric_\outspace\bigl(0,\psi(x)\bigr) \leq \bar E_h$.
\end{defn}

\begin{thm}[Tracking error convergence]
\label{thm:trackerrconv}
Assume that, for some $q\geq0$, we choose $\hestthresh = 2 \obserrbnd + q$ in our LACKI prediction rule and that the sequence of innovations $\seq{F_n(x_n)}{n\in \nat} $ as well as the reference $\seq{\xi_n}{n \in \nat}$ are bounded. 
If the initial error innovation function is bounded, i.e. if $\exists b \in \Real \forall x: \norm{F_0(x)}_\infty \leq b $, and, if $M$ is a stable matrix, i.e. if $\specrad(M) <1$, then the tracking error converges to the interval $[0,\frac q 2 + 2  \obserrbnd  + 2 \bar E_h]$. That is,
\[\norm{e_n}_\infty \stackrel{n \to \infty}{\longrightarrow} [0, \sigma\bigl(\frac q 2 + 2  \obserrbnd  + 2 \bar E_h\bigr )] \]
where $\sigma := \tinc  \sum_{i=0}^\infty \matnorm{M^i} < \infty$.

\end{thm}
\begin{proof} Let $\norm{\cdot} := \norm{\cdot}_\infty$ with accociated matrix norm $\matnorm{\cdot} := \sqrt{d} \specnorm{\cdot}$.

Let $\epsilon >0$. We desire to show: 
\begin{equation}\label{}
\exists N \in \nat \forall n\geq N: \norm{e_n} \leq \epsilon+\frac q 2 + 2  \obserrbnd  +2 \bar E_h.
\end{equation}

If sequence $\seq{F_n(x_n)}{n\in \nat} $ is bounded then, by Thm. \ref{thm:normboundsboundeddisturbmainbody}, the error sequence \seq{e_n}{n\in \nat} is bounded.  That is, $\exists b \in \Real \forall n: \norm{e_n} \leq \beta$.
Knowing that the error dynamics are bounded by some $\beta \geq 0$ we see that $\matnorm{M^{k}} \, \norm{\vc e_n} \leq \matnorm{M^{k}} \beta \stackrel{k \to \infty}{\longrightarrow} 0$. Here, the convergence to zero follows from the assumption that $M$ is a stable matrix. Hence,
we have:

$$(I)  \, \, \forall n \exists k_0(n) \in \nat \forall k \geq k_0(n): \matnorm{M^{k}} \, \norm{\vc e_n} \leq \frac\epsilon 2. $$

If in addition, the reference is bounded this implies that the sequence $(x_n)$ is bounded, too. Thm. \ref{thm:vanisishingseqprederr_LACKI} implies convergence of the innovations and hence, assuming $\metric_\outspace (f,f') = \norm{f-f'}$, we have:
\begin{equation}\label{eq:Fnconv}
\forall \varepsilon >0 \exists n_0 \forall n \geq n_0 : \norm{F_n(x_n)}\leq  \varepsilon +\frac q 2 + 2  \obserrbnd  +2 \bar E_h.
\end{equation}

Referring to (i) in Thm. \ref{thm:normboundsboundeddisturbmainbody}, With a change of variables we can follow analogous steps to convert Ineq. (\ref{ineq:errnormbasic}) to state that for all $k \in \nat, n \in \nat_0$ we have: 
\begin{align}
	\norm{\vc e_{n+k}} 	&\leq \matnorm{M^{k}} \, \norm{\vc e_n} + \tinc Q_{n:n+k}	 \sum_{i=0}^{k-1}  \matnorm{M^{k-1-i}} 	\label{ineq:errnormgeneral}
\end{align}
$Q_{n:n+k}:=\max\{\norm{F_n(x_n)},\ldots,\norm{F_{k+n-1}(x_{k+n-1})} \}$. 
With Gelfand's formula and the standard root test for series it is easy to establish convergence of the series: That is, $\sigma = \lim_{k \to \infty} \tinc \sum_{i=0}^{k-1}  \matnorm{M^{k-1-i}} < \infty$. And, we have $\tinc \sum_{i=0}^{k-1}  \matnorm{M^{k-1-i}} \leq \sigma, \forall k$.
Hence,
\begin{equation}
\norm{\vc e_{n+k}} \leq \matnorm{M^{k}} \, \norm{\vc e_n} + \sigma \, Q_{k:n+k}, \forall n \in \nat_0, k \in \nat.
\label{eq:errdyn320894}
\end{equation}

With (\ref{eq:Fnconv}) follows that there exists $n_0 \in \nat_0$ such that we have: $$ (II) \,\,  \forall k \in \nat : Q_{n_0:n_0+k}\leq  \frac{\epsilon}{2 \sigma} +\frac q 2 + 2  \obserrbnd  +2 \bar E_h.
 $$

Combining (I) and (II) with Eq. \ref{eq:errdyn320894} allows us to conclude that for any $n \geq N:= n_0 +k_0(n_0)$  we have 
\begin{equation*}
\norm{\vc e_{n}} \leq \frac \epsilon 2 + \sigma \Bigl (\frac{\epsilon}{2  \sigma} +\frac q 2 + 2  \obserrbnd  +2 \bar E_h\Bigr) =\epsilon +\sigma (\frac q 2 + 2  \obserrbnd  +2 \bar E_h \Bigr). 
\label{eq:errdyn3208942}
\end{equation*}
\end{proof}
Note, since the error converges to a bounded set the state will converge to the target trajectory. So, if the target trajectory is bounded, the continuity of the control law (as a function of state) implies that the control is bounded as well.

\begin{cor}
In the special case of error-free observations of a H\"older continuous target function, choosing a parameter $\hestthresh =0$ implies that the tracking error vanishes, i.e. :
\[\norm{e_n}_\infty \stackrel{n \to \infty}{\longrightarrow} 0. \]
Furthermore, the control action sequence $\seq{u(x_n)}{n \in \nat}$ converges, provided the reference trajectory $\seq{\xi_n}{n \in \nat}$ is bounded.
\end{cor}

\begin{proof} The convergence statement is an immediate consequence of the preceding theorem. 
Remember from Sec. \ref{sec:mrac} that the control action at time $n$ is of the form $u_n := u(x_n) = -\predfn(x_n) - K e_n +c $ for some constant $c$. We show that $(u_n)$ is a Cauchy sequence, provided that the reference sequence $\xi_n$ is. Since $\inspace$ is a Hilbert space, the desired convergence result follows.

So, let $\epsilon >0$. Since $(e_n),(\xi_n)$ converge, also the state sequence $(x_n)$ converges. Hence, all three are convergent Cauchy sequences. In particular, there is $N$ such that for all $n,m >N$: $\norm{e_n-e_m}< \frac{\epsilon} {2 \matnorm{K}}$ and $\norm{x_n-x_m}< \frac{\epsilon} {2\bar L}$. 
Hence, utilising the definition of the control law and the fact that all predictors are H\"older continuous with H\"older constant $\bar L$, for all $m,n >N:$ $\norm{u_n-u_m} \leq \matnorm{K} \norm{e_n -e_m} + \norm{\predfn(x_n) - \predf_m(x_m) } \leq \frac{\epsilon}{2} + \bar L \norm{x_n -x_m} \leq \epsilon $. Therefore, $(u_n)$ is a Cauchy sequence and hence, convergent.

\end{proof}


\section{Conclusions}
\label{sec:concl}
In this paper, we have introduced \emph{Lazily Adapted Constant Kinky Inference (LACKI)} as an approach to nonparametric machine learning. Our method was built on the framework of \emph{Kinky Inference} (which in turn is a generalisation of well-known approaches such as \emph{Lipschitz Interpolation} \cite{Sukharev1978,Zabinsky2003,Beliakov2006} and \emph{Nonlinear Set Membership (NSM)}  methods \cite{Milanese2004}). 
 Our approach inherits the numerical simplicity of these methods. On top of this, it can deal with bounded additive observational errors and does not require a priori knowledge about a H\"older constant of the underlying target function--  instead it estimates the constant online from the data. This of course is of great practical interest since this endows LACKI with superior black-box learning capabilities while still allowing us to give theoretical guarantees on learning and control success.
 
To avoid the need to specify the H\"older constant, LACKI adapts the parameter $L(n)$ to reflect a modification of the empirical estimate of the H\"older constant of the underlying target function. 
The adapted parameters were carefully defined to be bounded even if the target function is not H\"older continuous and the data is subject to bounded observational uncertainty.
This allowed us to establish several theoretical guarantees of worst-case consistency. That is, we provided asymptotic guarantees on the ability to learn any H\"older (and non-H\"older) continuous target function as well as convergence rates. Our derivations focussed on worst-case prediction error bounds.

Future work will investigate in how far the bounds can be tightened further (albeit we do not expect that worst-case guarantees could be given that avoid the curse of dimensionality without imposing more confining assumptions on the target functions and the nature of the observational uncertainties). 
However, if we were to shift our attention away from worst-case error analysis under general (possibly systematic) observational uncertainties towards standard mean-square risk analysis and assumptions prevalent in probability theory, we believe less conservative consistency might be establishable. 
To this end, we believe it is possible to modify our proofs to translate mean-square consistency derivations for nearest-neighbor regression methods (e.g. as discussed in \cite{Gyoerfi2002}) to our LACKI approach. This might provide a theoretical underpinning for the smoothing properties of our approach observed in the presence of i.i.d. stochastic noise (refer to Fig. \ref{fig:LACKInoise} and Fig. \ref{fig:LACKInoise2}).

Our learning-theoretic considerations were supplemented by an application of LACKI to 
online learning-based model-reference adaptive control. In a simulated aircraft control problem with nonlinear model uncertainty, we compared our LACKI-based controller against other learning-based methods that were recently proposed in the control literature. Across a range of performance metrics and randomised problem instances, LACKI-MRAC exhibited consistent and robust performance and outperformed its competitors on the majority of randomised test cases.

For discrete-time systems with additive, bounded, nonlinear uncertainty, we provided theoretical guarantees on the tracking success of our LACKI-MRAC controller. For the online learning case where the LACKI learner was assumed to be continuously updated with the most recently visited state observation, we proved tracking success up to a term dependent on the observational error.

In future work, we would like to apply our LACKI learning method to more challenging control tasks that require planning.
To this end, it might be beneficial to link our results to recent work on NSM-based model-predictive control (e.g. \cite{Canale2014}). We believe that the worst-case analysis we focus on in this work is key to establish the necessary links to results existing in the robust MPC literature.

In this work, we have assumed that the observational noise was bounded; we have addressed the issue of unbounded noise in a companion paper \cite{POKIdraft2016} where the H\"older constant parameter $L(n)$ is found as the minimiser of a prediction loss estimator.  We also believe that these estimators could also be used to estimate the noise bound if this exists but is  unknown a priori.

A final suggestion for future work pertains to the question of speeding up prediction time. In the context of Lipschitz interpolation, Beliakov \cite{Beliakov2006} proposed a way to organise the sample in a tree-structure (in lieu to KD-trees utilised in nearest-neighbor search) with the aim to reduce the prediction time to be logarithmic in the sample size. While applicable to our LACKI approach in the batch learning setting, it is not clear to us how this tree structure could be efficiently updated in an online learning setting. Future work could explore avenues of connecting his idea of appealing to notions of generalised nearest-neighbor search to existing efficient approaches of \emph{online} nearest neighbor search.

%

\section{Acknowledgements}
I am grateful for useful feedback from Jan Maciejowski, Carl Rasmussen, Stephen Roberts and Daniel Limon. I would also like to thank Girish Chowdhary and Hassan Kingravi who generously supplied me with their code that allowed me to most closely reproduce their work and use it for my comparisons in Sec. \ref{sec:KIMRAC}. I also gratefully acknowledge funding via the AIS project, NMZR/031 RG64733.
\bibliographystyle{plain}
\bibliography{content/lit}
\appendix 
\section{Supplementaries}
\subsection{H\"older and Lipschitz continuity for inference}
\label{sec:Hoelder_brief}
\subsubsection{Introduction and related work}

H\"older continuous functions are uniformly continuous functions that may exhibit infinitely many points of non-differentiability and yet are sufficiently regular to facilitate inference. That is, they have properties that make it possible to make assertions of global properties on the basis of a finite function sample. 

H\"older continuity is a generalisation of Lipschitz continuity.  Lipschitz properties are widely used in applied mathematics to establish error bounds and, among many other, find application in optimisation \cite{Shubert:72,direct:93} and quadrature \cite{Baran2008,curbera1998,dereich2006} and are a key property to establish convergence properties of approximation rules in (stochastic) differential equations \cite{kloedenandplaten1992,Gardiner2009}. 
Furthermore, most machine learning methods for function interpolation seem to impose smoothness (and thereby, H\"older continuity) on the function. For instance, with our Lem. \ref{lem:Hoeldarithmetic} derived below, it would be possible to show that any finite \textit{radial basis function neural network} \cite{Broomhead1988} with a smooth basis function is H\"older continuous on a compact domain. Or, a \textit{Gaussian process} with a smooth covariance function also has a smooth mean function and a.s. smooth sample paths  \cite{GPbook:2006,grimmetbook2001}. Therefore, posterior inference over functions on compact support made with such a Gaussian process on the basis of a finite sample is H\"older continuous.

Recently, we have become aware of related work published in mathematical and operations research journals \cite{Cooper2006,Cooper1995,Zabinsky2003,Beliakov2006,Beliakovsmoothing2007}. For instance, Zabinsky et. al. \cite{Zabinsky2003} consider the problem of estimating a one-dimensional Lipschitz function (with respect to the canonical norm-induced metric). Similar to the analysis we employ to establish our guarantees, they use a pair of bounding functions and make predictions by taking the average of these functions. While we have developed our kinky inference rule independently, it can be seen as a generalisation of their approach. Our method provides extensions to H\"older continuous multi-output functions over general, multi-dimensional (pseudo-) metric spaces, can cope with with erroneous observations and inputs, can fold in additional knowledge about boundedness, learn parameters from data and provides different guarantees such as (uniform) convergence of the prediction uncertainty. 
As part of the analysis of our method, we construct delimiting functions we refer to as \textit{ceiling} and \textit{floor} functions. The construction of similar functions is a recurring theme that, in the standard Lipschitz context, can be found in global  optimisation \cite{Shubert:72}, quadrature \cite{Baran2008}, interpolation \cite{Beliakov2006,Beliakovsmoothing2007}, as well as in the analysis of linear splines for function estimation \cite{Cooper1995}. Cooper \cite{Cooper2006,Cooper1995} utilises such upper and lower bound functions in a multi-dimensional setting to derive probabilistic PAC-type error bounds \cite{Valiant1984} for a linear interpolation rule. He assumes the data is sampled uniformly at random on a hypercube domain. This precludes the application of his results to much of our control applications where the data normally is collected along continuous trajectories visited by a controlled plant. Our inference rule is different from his and our guarantees do not rely on distributional assumptions. This of course is important in control settings where the common assumption that the input data was drawn independently from a fixed distribution typically is not met.
In this thread of works, perhaps the work that is most closely related to ours is the function interpolation method of Beliakov \cite{Beliakov2006} that is a special case of a kinky inference rule: For a single-output function that is Lipschitz with respect to a special input space norm and where the data is error-free, the authour provides an algorithm that promises logarithmic prediction time. Unfortunately, many of his assumptions are unrealistic in a control setting. And, the improved prediction time is achieved by constructing a data structure from batch data which precludes its use in an online learning setting. However, future work might explore in how far his ideas can be 
converted into an online learning rule. Furthermore, in learning situations where Beliakov's interpolation method is applicable, our theoretical results extend to his work. For instance, our results show how Lipschitz constant estimation can be harnessed to render his approach a universal approximator.

Other work of relevance can be found in analysis. For instance, Miculescu \cite{Miculescu2000} presents work proving that any continuous function on a metric space is a uniform limit of a sequence of \emph{locally} Lipschitz functions and also mentions that the stronger statement, that every function is a limit of a sequence of globally Lipschitz functions, is not true in general. However, he cites earlier work \cite{Georganopoulos1967} that does show that every real-valued continuous function on \emph{compact} domain is a uniform limit of a sequence of Lipschitz functions. In some sense, our work develops a related statement as a by-product. 
From our convergence guarantee of the LACKI rule, we have derived constructive method for showing 
that any continuous function on compact domain is the uniform limit of a sequence of H\"older functions up to an arbitrarily small error.

Finally, in the context of control, Milanese and Novara \cite{Milanese2004} considered NSM methods for interpolation. For a fixed Lipschitz constant, their prediction rule can be
seen as a special case of ours without the $\lbf$ and $\ubf$ parameters and with special choices of metrics. Similar to us, they do consider the problem of estimating the Lipschitz constant from the data and consider bounded noise. However, they obtain the Lipschitz constant estimate via the maximum partial derivative of an arbitrarily chosen fitted parametric model of a bounded input set. And, they give no guarantees on the quality of the resulting estimator that is fitted to the data like this nor do they discuss the impact of the choice of parametric model or the fitting method on the quality of the estimator.
 


\subsubsection{Basic facts and derivations}
In preparation of subsequent parts of the work that take advantage of H\"older properties this section will proceed to establish essential prerequisites.
The remainder of this section is structured as follows: Firstly, we will go over basic definitions and engage in some preliminary derivations that will be of importance throughout this work.
While we do not claim novelty on any of the results we provide proofs for in this section, we had not found them in the literature and hence, had to derive them on our own.

Firstly, we commence with introducing the notions of (pseudo-) metric spaces.

\begin{defn}[(Pseudo-) metrics]
Let $\inspace$ be a set. A mapping $\metric_\inspace: \inspace^2 \to \Real$ is called a \emph{pseudo-metric} if it positive ($\forall x,x' \in \inspace: \metric_\inspace(x,x') \geq 0$) and satisfies the triangle inequality ($\forall x,x',x'' \in \inspace: \metric(x,x') \leq \metric_\inspace(x,x'') + \metric (x'',x')$). If furthermore the pseudo-metric $\metric$ is definite $(i.e. \forall x,x' \in \inspace_\inspace: \metric_\inspace(x,x') =0 \Leftrightarrow x=x')$ then 
the mapping $\metric$ is called a \emph{metric}. The set $\inspace$ endowed with a (pseudo-) metric $\metric_\inspace: \inspace^2 \to \Real$ or the pair $(\inspace, \metric_\inspace)$ are called \emph{(pseudo-) metric space}.
 \end{defn}

\begin{defn} 
Let $(\inspace ,\metric_\inspace ), (\outspace , \metric_\outspace )$ be two (pseudo-) metric spaces and 
$I \subset \inspace$ be an open set. A function $f: \inspace \to \outspace $ is called (L-p-) \emph{H\"older} 
(continuous) on $I \subset \inspace$ if there exists a \emph{(H\"older) constant} $L \geq 0$ and \emph{(H\"older) 
exponent} $p\geq 0$ such that 
\[\forall x,x' \in I : \metric_\outspace \bigl(f(x),f(x')\bigr) \leq L \, \bigl( \metric_\inspace (x,x') \bigr)^p. \]
We denote the space of all L-p- H\"older functions by $\hoelset{L}{}{p}$.
\end{defn}

H\"older functions are known to be uniformly continuous. 
A special case of importance is the class of $L$-\textit{Lipschitz} functions. These are H\"older continuous 
functions with exponent $p=1$. In this context, coefficient $L$ is referred to as\textit{ Lipschitz constant} or \textit{Lipschitz number}.

\begin{ex}[Square root function]\label{ex:sqrtfctHoelder}
As an example of a H\"older function that is not Lipschitz we can consider $x \mapsto \sqrt x$ on domain $I = [0+\epsilon,c]$ where 
$c >\epsilon \geq 0 $. For $\epsilon >0 $ the function is Lipschitz with $L = \sup_{x \in I} \frac{1} {2 \sqrt{x}}$. We can see that the 
coefficient grows infinitely large as $\epsilon \to 0$. By contrast, the function is H\"older continuous 
with H\"older coefficient $L=1$ and exponent $p=\frac 1 2 $ for any bounded $I \subset \Real$.
We can see this as follows: Let $\epsilon =0,$ $x,y \in I$ and, without loss of generality,  $y \geq x$. Let $\xi := \sqrt{x}, \gamma := \sqrt{y}$ and thus, $\gamma \geq \xi$. We have:
$\xi \leq \gamma $ $\Leftrightarrow 2 \xi^2 \leq 2\xi\gamma$ $\Leftrightarrow \gamma^2 - 2 \xi\gamma + \xi^2  \leq \gamma^2 - \xi^2$ $\Leftrightarrow (\gamma-\xi)^2  \leq \gamma^2-\xi^2$ $\Leftrightarrow \abs{\gamma-\xi}^2  \leq \abs{\gamma^2-\xi^2}$
$\Leftrightarrow \abs{\gamma-\xi}  \leq \sqrt{\abs{\gamma^2-\xi^2}}$  $\Leftrightarrow \abs{\sqrt{x}-\sqrt{y}} \leq \abs{y-x}^{\frac{1}{2}}$. Since, $x,y$ were chosen arbitrarily, we have shown H\"older continuity as desired.
\end{ex}

Most commonly, one considers H\"older continuity for the special case of the standard metric induced by a norm, i.e.  $\metric(x,x') = \norm{x-x'}$.
For a function $f: \inspace \to \outspace$, the H\"older condition becomes:
\[\forall x,x' \in I : \norm{f(x)-f(x')}_\outspace \leq L \, \norm{x-x'}_\inspace^p. \]

Similarly, we can consider H\"older continuity for each output component: 

\begin{defn} \label{def:outputwisehoelder}
Let $\outspace \subseteq \Real^m $ and $\inspace$ be a space endowed with a metric (or indeed a semi-metric) $\metric_\inspace$. Then, the function $f: \inspace \to \outspace$ is output-component-wise H\"older continuous with exponent $p$ and constant $L \in \Real^m_{\geq 0}$ if $f \in \hoelset L { } p$ where
$\hoelset L {\metric_\inspace} p:= \bigl\{ \phi: \inspace \to \outspace \,\ | \, \forall j \in \{1,...,m \}, \forall x,x' \in \inspace: \abs{\phi_j(x) - \phi_j(x')} \leq L_j \metric^p_\inspace(x,x') \bigr\}$
is the set of all functions whose component functions are H\"older continuous with respect to input space metric $\metric_\inspace$ and an output space metric that is induced by the canonical norm $\metric_{\outspace} (x,x')= \abs{x-x'}$.  
\end{defn}

\begin{remark}[Best H\"older constant]
Note for $p \in (0,1], 0 \leq L_1 \leq L_2$ we have $\hoelset {L_1} {\metric_\inspace} p \subseteq \hoelset {L_2} { } p$. The smallest $L^* \geq 0$ such that function if is $L^*-p-$ H\"older, $f \in \hoelset {L^*} {\metric_\inspace} p$, is called the \emph{best} H\"older constant of $f$.
\end{remark}

Generally, it is obviously true that $\hoelset{L}{}{p} \subseteq \hoelset{L'}{}{p}$ for $L' \geq L$.
With regard to the H\"older exponent, we will now show that smaller exponents are less restrictive than larger ones.

\begin{lem} \label{lem:hoeldexpprop2}
Let $\inspace$ be the input space (not necessarily bounded). For some $p \in (0,1], L \geq 0$ assume that  and $f:\inspace \to \outspace$ is locally $L-p$-H\"older continuous.
Then we have: (i) for any $ q \in (0,p]$, f is also locally $L-q$-H\"older. 
(ii) If $f:\inspace \to \outspace$ is bounded with $\sup_{x,x' \in \inspace} \metric_\outspace (f(x),f(x')) \leq B \in \Real$ and globally $L-p$ H\"older then $f$ is globally $L^*-q$-H\"older, where $L^* := \max\{L,B \}$ and $q \in (0,p]$. In particular, on compact support, Lipschitz continuity entails \emph{H\"older} continuity for any H\"older exponent $p \in [0,1)$.
\begin{proof}
(i) Let $p \in (0,1], f \in \hoelset L { } p $ and $p = q+r, r \in [0,1)$. Let $\xi \in \inspace$ and $I$ denote the intersection of the domain with an $\epsilon$-ball around $\xi$ such that $f$ satisfies the H\"older condition on $I$ and such that $\sup_{x,x' \in I} \metric_\inspace(x,x') \leq 1$.  For all $x,x' \in I$ we have $\metric_\outspace (f(x),f(x')) \leq L\metric_\inspace^p (x,x') = L \metric_\inspace^q(x,x') \metric_\inspace^r(x,x') \leq L \metric_\inspace^q(x,x')$ where the last inequality holds since $r \in [0,1)$ and $\sup_{x,x' \in I} \metric_\inspace(x,x') \leq 1$.

(ii) Let $x,x' \in \inspace$. If $\metric_\inspace(x,x') \leq 1$ we can show $\metric_\outspace (f(x),f(x'))  \leq L \metric_\inspace^q(x,x')$ following through the same sequence of inequalities as above in the proof of (i). Now, let $\metric(x',x) >1$. We have $\metric_\outspace (f(x),f(x')) \leq B \leq B  \metric_\inspace(x,x')^q$.
\end{proof}
\end{lem}

\begin{thm} \label{thm:hoelderconcat}
Let $(\statespace, \d)$ be a metric space and $f,g : \statespace \to \statespace$ be two H\"older continuous mappings with H\"older constants $L(f), L(g)$ and H\"older exponents $p_f,p_g$, respectively.
Then, the concatenation $h=f \circ g: \statespace \to \statespace $ is also H\"older continuous with H\"older constant $L(h):= L(f) L(g)^{p_f}$ and exponent $p_h:=p_g \, p_f$.
That is, 
\[\forall \state,\state' \in \statespace: \d\bigl(h(\state),h(\state')\bigr) \leq L(h) \, \bigl(\d(\state,\state')\bigr)^{p_h}.\]
\begin{proof}
$\d\bigl(f \circ g(\state),f\circ g(\state')\bigr)$ $\leq L(f)\,  \Bigl(\d(g(\state),g(\state'))\Bigr)^{p_f}$
$\leq L(f)\,  \Bigl(L(g)\, \d(\state,\state')^{p_g}\Bigr)^{p_f}$ 

$= L(f)\, L(g)\,^{p_f}   \Bigl(\d(\state,\state')\Bigr)^{p_g\, p_f} $ where in the first step we were using H\"older-continuity of $f$ and in the second, we were using H\"older continuity of $g$ combined with the fact that $(\cdot)^{p_f}$ is a monotonically increasing  function. 
\end{proof}
\end{thm}


Several numerical methods, such as Lipschitz optimisation \cite{Shubert:72}, rely on the knowledge of a Lipschitz constant. In the more general case of H\"older continuous functions this will correspond to the need of knowing a H\"older constant and exponent. To avoid having to derive these for each new function from first principles, we establish the following collection of facts that allows us determine bounds on H\"older constants of combinations of functions with known H\"older parameters.
While we provide proofs for a restatement in the H\"older continuous setting, a number of the following statements have also been proven in \cite{Weaver1999} in the context of Lipschitz algebras.

\begin{lem}[H\"older arithmetic] \label{lem:Hoeldarithmetic}
Let, $I,J \subset \inspace$ where $\inspace$ is a metric space endowed with metric $\metric$. Let $f : \inspace \to \Real$ be H\"older on $I$ with constant $L_I (f) \in \Real_+$ 
and $g :\inspace \to \Real$ be H\"older on $J$ with constant $L_J (g) \in \Real_+$. Assume both functions have the same H\"older exponent $p \in (0,1]$. That is, $\forall x, x' \in \inspace: \abs{f(x)-f(x')} \leq L(f) \metric(x,x')^p$ and  $\forall x, x' \in \inspace: \abs{g(x)-g(x')} \leq L(g)  \metric(x,x')^p$.
We have:

\begin{enumerate}
	\item Mapping $x \mapsto |f(x)|$ is H\"older on $I$ with constant $L_I(f)$ and exponent $p_f$.
	\item If $g$ is H\"older on all of $J=f(I)$ the concatenation $g \circ f: t \mapsto g(f(t))$ is H\"older on $I$ with constant 
	      $L_I(g \circ f) \leq$ $L_J (g) \, L_I^p(f)$ and exponent $p^2$.
	\item Let $r \in \Real$. $r \, f: x \mapsto r \, f(x)$ is H\"older on $I$ having a constant $L_I (r \,f) = |r| \, L_I(f)$.
	\item $f+g: x \mapsto f(x) + g(x)$ has H\"older constant at most $L_I(f) + L_J(g)$.
	\item Let $m_f = \sup_{x\in \inspace } f(x)$ and $m_g = \sup_{x \in \inspace } g(x)$. Product function $f\cdot g: x \mapsto f(x) \, g(x)$ has H\"older exponent $p$ and a H\"older constant on $I$ which is at most $(m_f \, L_J(g)+ m_g \, L_I(f))$.
	\item For some countable index set $\indsett$, let the sequence of functions $f_i$ be H\"older with exponent $p$ and constant $L(f_i)$ each. Furthermore, let $H(x) =\sup_{i \in \indsett} f_i(x) $ and $h(x) := \inf_{i \in \indsett} f_i(x)$ be finite for all $x$. Then $H,h$ are also H\"older with exponent $p$ and have a H\"older constant which is at most $\sup_{i \in \indsett} L(f_i)$.
	\item Let $b := \inf_{x \in \inspace }| f(x)| > 0$ and let 
	$\phi(x) = \frac{1}{f(x)}, \forall x \in \inspace$ be well-defined.  
	      Then $L_I(\phi) \leq b^{-2} \, L_I(f)$.  
	\item Let $p=1$ (that is we consider the Lipschitz case), let $I$ be convex and $\metric(x,x') = \norm{x-x'}$ where $\norm{\cdot}$ is a norm that induces a sub-multiplicative matrix norm (e.g. all $p-$ norms are valid). $f$ cont. differentiable on I $\Rightarrow$ $L_I(f) \leq \sup_{x \in I } \norm{\nabla f(x)}. $ 
	For one-dimensional input space, $\inspace = \Real$, $L_I(f) = \sup_{x \in I } \abs{\nabla f(x)}$ is the smallest Lipschitz number. 
	 \item Let $c \in \Real$, $f( t) = c, \forall x \in I $. Then $f$ is H\"older continuous with constant $L_I(f) =0$ and for any coefficient $p_f \in \Real$.  
	\item $L_I(f^2) \leq 2 \, L_I(f)\, \sup_{t \in I} f\,$.
	\item With conditions as in 8), and input space dimension one, we have $\forall q \in \mathbb Q : L_I(f^q) = |q| \,\sup_{\tau \in \indset } |f^{q-1}(\tau) \, \dot f(\tau)| $.
\end{enumerate}
\begin{proof}

\textbf{1)}  We show $|f|$ has the same constant and exponent as $f$. Let $X,X' \in \inspace $ arbitrary. 
We enter a case differentiation:

\textit{1st case: $f(x), f(x') \geq 0$}. 

Hence, $\bigl| \abs{f(x)}- \abs{f(x')} \bigr| = \bigl| f(x) - f(x') \bigr|  \stackrel{f \,Hoelder}{\leq} L_I(f) \metric(x,x')^{p}$.\\

\textit{2nd case: $f(x) \geq 0, f(x') \leq 0$.} 

Note, $|y| = - y$, iff $y \leq 0$. Hence,  $\bigl| |f(x)| - |f(x')| \bigr| \leq \bigl| |f(x)| + |f(x')| \bigr| $
$= \bigl| |f(x)| - f(x') \bigr|  =  \bigl| f(x) - f(x') \bigr| \leq L_I(f) \, \metric(x,x')^{p}$.\\

\textit{3rd case: $f(x) \leq 0, f(x') \geq 0$.} Completely analogous to the second case.\\

\textit{4th case: $f(x), f(x') \leq 0$}. 

$\bigl| |f(x)| - |f(x')| \bigr| = \bigl| f(x') - f(x) \bigr|= \bigl| f(x) - f(x') \bigr|  \stackrel{f \, Hoelder}{\leq} L_I(f) \metric(x,x')^{p}$.\\

The remaining points are also proven in \cite{Weaver1999} in the context of Lipschitz functions.

\textbf{2)} Special case of Thm. \ref{thm:hoelderconcat}.
%

\textbf{3)}  For arbitrary $x,x' \in \inspace , r \in \Real$ we have:

$\bigl| r \, f(x) - r \, f(x')| \bigr| = |r|\, |f(x) - f(x')| \leq |r| \,L_I(f)\,  \metric(x,x')^{p}$.\\ 

\textbf{4)}  For arbitrary $x,x' \in \inspace , r \in \Real$ we have:

$\bigl| g(x) + f(x) - (g(x') + f(x'))| \bigr| = \bigl| g(x)  - g(x') + f(x)- f(x')| \bigr|$ 
$\leq \bigl| g(x)  - g(x')\bigr|  + \bigl| f(x)- f(x')| \bigr|$ $\leq (L_J(g)+L_I(f))\,  \metric(x,x')^{p}$.\\

\textbf{5)}  Let  $x,x' \in \inspace $, $d := f(t) - f(t')$.

$\bigl| g(x) f(x) - g(x')  f(x') \bigr| = \bigl| g(x) (f(x') +d) - g(x') f(x') \bigr|$ 
$= \bigl|\bigl( g(x) - g(x') \bigl)  f(x')+ g(x)  d \bigr|  $

$\leq \bigl| g(x) - g(x') \bigr|  |f(x')|   + \bigl|g(x)\bigr| \,  |d|  $
$\leq L_I(g) \metric(x,x')^p  |f(x')|   + \bigl|g(x)\bigr| \,  L_I(f) \metric(x,x')^p  $

$\leq L_I(g) \metric(x,x')^p  \sup_{x' \in \inspace } \{|f(x')|\}  \\ + \sup_{x \in \inspace }\{\bigl|g(x)\bigr|\} \,  L_I(f) \metric(x,x')^p  $

$= \Bigl(L_I(g)  \sup_{x' \in \inspace } \{|f(x')|\}   \\+ \sup_{x \in \inspace }\{\bigl|g(x)\bigr|\} \,  L_I(f)\Bigl) \metric(x,x')^p  $.\\

\textbf{6)}  The proof of Proposition 1.5.5 in \cite{Weaver1999} proves our statement if one replaces their 
metric $\rho$ by $\metric^p$.

\textbf{7)}  Let  $x,x' \in \inspace $.
$\bigl| \frac{1}{f(x)} - \frac{1}{f(x')} \bigr|$ 
$=\bigl| \frac{f(x')}{f(x') f(x)} -\frac{f(x)}{f(x') f(x)} \bigr|$ 
$= \frac{\bigl|f(x')-f(x) \bigr|}{\bigl|f(x')\bigr| \bigr| f(x)\bigr|}$ 
$\leq \frac{L_I(f) \metric(x,x')^p}{\inf_{x \in \inspace } |f(x)|^2}$.\\

\textbf{8)} Let $p=1$ and $\metric(x,x') = \norm{x-x'}$ be a norm that induces a sub-multiplicative matrix norm. Define $\ell := \sup_{x \in I } \norm{ \nabla f(x) } = L_I(f)$. 
Firstly, we show that it is a Lipschitz constant: Let $x,x' \in I $ and 
$\overline{xx'} := \{tx + (1-t) x' \, | \, t \in [0,1]\}$. 
Owing to convexity of I, $\overline{xx'} \subset I$. Due to the mean value theorem $\exists \xi \in \overline{xx'} \subset I: |f(x) - f(x')|=  T_\xi (x-x')$. where $T_\xi (x) = \SP{\nabla f(\xi)}{ x}$ is a linear OP. Assuming the derivative of $f$ is bounded, $T_\xi$ is a bounded OP and we have $\abs{T_\xi (x-x') } \leq \matnorm{T_\xi} \norm{x-x'}$ where 
$\matnorm{T_\xi} = \sup_{\norm{x} = 1} \abs{\SP{\nabla f(\xi)}{x}} \leq \norm{\nabla f(\xi)}$. In conjunction,
$|f(x) - f(x')| \leq \norm{\nabla f(\xi)} \norm{x-x'}$. 

Secondly, we show that $\ell$ is the smallest Lipschitz constant in the one-dimensional case: Let $\bar \ell$ be another Lipschitz constant on $I$ such that $\bar \ell \leq \ell$. Of course, here $\norm{\cdot} = \abs{\cdot}$. Since $I$ is compact and $\norm{\nabla f(\cdot) }$ is continuous, there is some $\xi \in I$ such that $\norm{\nabla f(\xi)} = \sup_{x \in I} \norm{\nabla f(\xi)} = \ell$. Pick any sequence $(y_k)_{k=1}^\infty$ contained in $I$ such that $y_k \stackrel{k \to \infty}{\longrightarrow} \xi$ and $y_k \neq \xi$.
$\forall k: y_k \in I $ and $\bar \ell$ is a Lipschitz constant on $I$. Hence, $ \bar \ell \geq \frac{|f(y_k) - f(\xi)|}{\norm{y_k-\xi}}\stackrel{k \to \infty}{\longrightarrow} \norm{ \nabla f(\xi)} = \ell$. Thus, $\bar \ell = \ell$.

\textbf{9)} Trivial. 

\textbf{10)} Special case of property 5).

\textbf{11)} $L(f^q) \stackrel{8)}{=} \sup_{\tau \in \indset } |\dif{}{t} f^q(\tau)| = |q| \,\sup_{\tau \in \indset } |f^{q-1}(\tau) \, \dot f(\tau)| $. 
\end{proof}
\end{lem} 

As a simple illustration, assume we desire to establish that $f(t) = \max\{ 1- 3 \sin(t), \exp\bigl(- \sin(t) \bigr)\}$ is Lipschitz and to find a Lipschitz number on $\Real$. Application of 8. shows that $t \mapsto \sin(t)$ and $t \mapsto \exp(- t)$ have a Lipschitz number of $1$. Application, of 2., 9. 1. and 6. then show that $L(f) =3$ is a Lipschitz number of $f$.




\subsection{H\"older continuity of the exponentiated map }

The main objective of this subsection is to show that, for each $s \in \inspace$, the function $\phi_s: x \mapsto L \metric(x,s)^p$ is H\"older continuous with coefficient $L$ and exponent $p$ with respect to (pseudo-) metric $\metric: \inspace^2 \to \Real$. This is important in our derivations since these maps are essential building blocks of the kinky inference procedure. Therefore it is easy to employ Lem. \ref{lem:Hoeldarithmetic} to show that the full kinky inference rule is H\"older continuous.

To establish the H\"older regularity result, we will first show that $(x,y) \mapsto \metric(x,y)^p$, for $p \in (0,1] $, is a metric.  This can be utilised to show that
$\phi_s \in \hoelset L \metric p$.

Before proceeding we need to establish a few facts. 
Firstly, we remind ourselves of the following well-known fact:
\begin{lem} \label{lem:pd_n_concave_subadditive}
A nonnegative, concave function $g:\Real_{\geq 0} \to \Real$ with $g(0) = 0$ is subadditive. 
That is, $\forall x,y \in \Real_{\geq 0}: g(x+y) \leq g(x) + g(y)$. 
 \begin{proof}
If $x = y = 0$ then subadditivity trivally holds:  $0=g(x+y) \leq g(x) + g(y) = 0$.
So, let $x, y \in \Real_+$ such that $x >0 \vee y >0$.
We have, $g( x +y) = \frac{x}{x+y} g(x+y) + \frac{y}{x+y} g(x+y) \leq g(\frac{x}{x+y} (x+y) ) +  g(\frac{y}{x+y}(x+y)) = g(x) + g(y)$.
Taking into account that $\frac{x}{x+y}, \frac{x}{x+y} \in [0,1]$, the last inequality can be seen as follows:

 Since $g$ is concave we have 
$\forall p \in [0,1], x \in \Real: g(p x) =  g(px + (1-p) 0) \geq p g(x) + (1-p) g(0) \geq p g(x)  $. The last inequality follows from $g(0) \geq 0$.
\end{proof}
\end{lem}

\begin{lem} \label{lem:x2p_pdNsubadd}
 Let $h: \begin{cases} \Real_{\geq 0} \to \Real_{\geq 0},\\ \, x \mapsto x^p\end{cases}$, for $p \in (0,1] $. $h$ is positive definite and subadditive. 
 That is, (i) $h(0) = 0 $ and  $h(x) > 0, \forall x \neq 0$ and (ii) $\forall x,y \in \Real_{\geq 0}: h(x+y) \leq h(x) + h(y)  $.
 Moreover, $h$ is concave. If $p \in (0,1)$, h is strictly concave. 
 
\begin{proof}
\textit{Pos. def. :} $h(0) = 0$.  Also $\lim_{x \to 0_+} h(x) =0$. Since $\nabla h (x) = p h^{p-1}(x) >0 $ for $x >0$, $h$ is strictly monotonically increasing on $\Real_+$. Hence, $h(x) > 0, \forall x >0$. 

\textit{Concavity:} If $p =1$, $h$ is linear and hence, concave. If $p \in (0,1)$, $\nabla^2 h(x) = p (p-1) h(x)^{p-2} > 0$ regardless of $x$.

\textit{Subadditivity:} Follows directly with Lem. \ref{lem:pd_n_concave_subadditive} on the basis of established positive definiteness and concavity.
\end{proof}
\end{lem}

\begin{lem}\label{lem:hoeldererror_metric}
Let $p \in (0,1]$.
With definitions as above, we assume that set $\inspace$ is endowed with a pseudo- metric $\metric$. Function
$\metricp: \begin{cases} \inspace \times \inspace \to \Real_{\geq 0} \\ (x,y) \mapsto \Bigl(\metric(x,y)\Bigr)^p \end{cases}$ is a pseudo-metric on $\inspace$.
If $\metric$ is a metric then so is $\metricp$.
\begin{proof}
 By Lem. \ref{lem:x2p_pdNsubadd}, $x\mapsto x^p$ is pos. def. and sub-additive. Therefore, positive definiteness and the triangle inequality of $\metric$ readily extend to $\metricp$ as follows: 

\textit{Pos. def.:}
Let $x=0$. $\metricp(x,x) = \metric(x,x)^p = 0^p = 0$. If $x \neq 0$ then $k :=\metric(x,x) \neq 0$. Hence   $\metricp(x,x) = d(x,x)^p = k^p \neq 0$.

\textit{Triangle inequality:}
Choose arbitrary $x,y,z \in \inspace $. We have $\metricp(x,z) = \metric(x,z)^p \leq (\metric(x,y) + \metric(y,z) \bigr)^p \leq \metric(x,y)^p + \metric(y,z)^p = \metricp(x,y) + \metricp(y,z)$. Here, the first inequality followed from the triangle inequality of pseudo-metric $\metric$ and the second from subadditivity properties established in Lem. \ref{lem:x2p_pdNsubadd}.

\textit{Symmetry:} If $\metric$ is a metric it is symmetric. Hence, $\metricp(x,y) = \metric(x,y)^p = \metric(y,x)^p = \metricp(y,x), \forall x,y \in \inspace $ in which case $\metric$ also is symmetric.
\end{proof}
\end{lem}

Before proceeding we establish a slight generalisation of the well-known \textit{reverse triangle inequality}:

\begin{lem}[Reverse Triangle Inequality] \label{thm:revtriangle}
Let  $\inspace$ be a set and $\metric : \inspace^2 \to \Real$ a symmetric function that satisfies the triangle inequality. 
That is,  $\forall x,y,z \in \inspace: \metric (x,y) = \metric(y,x) \wedge \metric(x,z) \leq \metric(x,y) + \metric(y,z)$.

Then \[\forall x,y,z \in \inspace: \abs{\metric(x,y) - \metric(z,y)} \leq \metric(x,z).\]
\begin{proof}
For contradiction, assume 
$\abs{\metric(x,y) - \metric(z,y)}>\metric(x,z)$ for some $x,y,z \in \inspace$.
This is implies  
$
(i)\,\,\, \metric(x,y) - \metric(z,y)>\metric(x,z) 
\, \, $  or  $
\,\,  (ii) \,\,\,\, \metric(z,y)-\metric(x,y) >\metric(x,z)$.
Both inequalities contradict the triangle inequality:
$(i) \Leftrightarrow \metric(x,y)  >\metric(x,z) +\metric(z,y) $ and 
$(ii) \Leftrightarrow  \metric(z,y) > \metric(x,z) + \metric(x,y) =  \metric(z,x) + \metric(x,y)$.

\end{proof}
\end{lem}

\begin{lem}
Let $\metric$ be a (pseudo-) metric on $\inspace$. For arbitrary $s \in \inspace $ we define $\phi_s: \inspace \to \Real $ as $\phi_s: x \mapsto \metric (x,s) $.
$\phi_s $ is Lipschitz with respect to metric $\metric$. That is, \[\forall x,y \in \inspace : \abs{\phi_s(x) - \phi_s(y) } \leq \metric (x,y). \]
\begin{proof}
$\abs{\phi_s(x) - \phi_s(y)} = \abs{\metric(x,s) - \metric(y,s)} \stackrel{Thm. \ref{thm:revtriangle}} {\leq} \metric(x,y), \forall x,y,s \in \inspace $.
\end{proof}
\end{lem}
Finally, combining this lemma with Lem. \ref{lem:hoeldererror_metric} immediately establishes that mappings of the form $\metric(\cdot,s)^p$ are H\"older continuous with exponent $p$ ($\in (0,1]$):

\begin{thm} \label{thm:d2pmapishoelder}
Let $\metric$ be a (pseudo-) metric on set $\inspace$ and let $p \in (0,1], L \geq 0$. For arbitrary $s \in \inspace $ we define $\phi_s:\begin{cases}  \inspace \to \Real \\ x \mapsto L \, \bigl(\metric (x,s) \bigr)^p\end{cases}$.
$\phi_s $ is H\"older with exponent $p$. That is, \[\forall x,y \in \inspace : \abs{\phi_s(x) - \phi_s(y) } \leq L \, \metric (x,y)^p. \]
In particular, for any norm $\norm \cdot$ on $G$ and $s \in \inspace $, mapping $x \mapsto L \norm{x-s}^p$ is H\"older with constant $L$ and exponent $p$.

\begin{proof}
By Lem. \ref{lem:hoeldererror_metric}, $\metric^p(\cdot,\cdot)$ is a (pseudo-) metric on $\inspace$. Hence, Lem. \ref{thm:revtriangle} is applicable yielding:
$\abs{\phi_s(x) - \phi_s(y)} = L \, \abs{\metric(x,s)^p - \metric(y,s)^p} \stackrel{Lem. \ref{thm:revtriangle} } {\leq} \metric(x,y)^p, \forall x,y,s \in \inspace$. The last sentence concerning the norms follows from the fact that a mapping $(x,y) \mapsto \norm{x-y}$ defines a metric.
\end{proof}
\end{thm}

\bibliographystyle{plain}

\end{document}